\documentclass[final,1p]{elsarticle}
\graphicspath{{./Figures/}}
\usepackage{amssymb,latexsym,amsmath,amsthm}
\usepackage[colorlinks,linkcolor=blue, anchorcolor=blue, citecolor=blue]{hyperref}
\usepackage{color,caption,cases}
\usepackage{bm}
\usepackage{graphicx,graphics}
\usepackage{booktabs}
\usepackage{hyperref}
\usepackage{tabularx}
\usepackage[T1]{fontenc}
\usepackage{epstopdf}
\usepackage{multirow, makecell}
\biboptions{numbers,sort&compress} 
\allowdisplaybreaks
\newtheorem{theorem}{Theorem}[section]

\newtheorem{lemma}[theorem]{Lemma}

\theoremstyle{definition}

\newtheorem{example}[theorem]{Example}

\theoremstyle{remark}
\newtheorem{remark}[theorem]{Remark}

\newcommand{\f}{\frac}
\newcommand{\p}{\partial}

\newcommand{\mR}{\mathrm{R}}
\newcommand{\mF}{\mathrm{F}}
\newcommand{\mS}{\mathrm{S}}
\newcommand{\mo}{\mathcal{O}}
\numberwithin{equation}{section}

\begin{document}    

\begin{frontmatter}

\title{A decoupled linear, mass-conservative block-centered finite difference method for the Keller--Segel chemotaxis system}

\author[OUC]{Jie Xu} \ead{jxu129@163.com}
\author[OUC]{Hongfei Fu\corref{Fu}}\ead{fhf@ouc.edu.cn}

\address[OUC]{School of Mathematical Sciences \& Laboratory of Marine Mathematics, Ocean University of China, Qingdao 266100, P.R. China}

\cortext[Fu]{Corresponding author.}

\begin{abstract}
As a class of nonlinear partial differential equations, the Keller--Segel system is widely used to model chemotaxis in biology. In this paper, we present the construction and analysis of a decoupled linear, mass-conservative, block-centered finite difference method for the classical Keller--Segel chemotaxis system. We show that the scheme is mass conservative for the cell density at the discrete level. In addition, second-order temporal and spatial convergence for both the cell density and the chemoattractant concentration are rigorously discussed, using the mathematical induction method, the discrete energy method and detailed analysis of the truncation errors. Our scheme is proposed and analyzed on non-uniform spatial grids, which leads to more accurate and efficient modeling results for the chemotaxis system with rapid blow-up phenomenon. Furthermore, the existence and uniqueness of solutions to the Keller--Segel chemotaxis system are also discussed. Numerical experiments are presented to verify the theoretical results and to show the robustness and accuracy of the scheme.
\end{abstract}


\end{frontmatter}

\section{ Introduction}
In 1970s, Keller and Segel \cite{KS'70, KS'71} made groundbreaking contributions by developing a set of partial differential equations (PDEs) to model chemotaxis - a crucial biological process wherein a cell (or organism) migrates in reaction to a chemical stimulus, either attractive or repulsive. Although their pioneering work has been acknowledged to contain some biological inaccuracies, particularly as their model can lead to unbounded solutions within a finite time (a scenario not observed in nature), the classical Keller--Segel system continues to hold substantial mathematical significance.

Mathematically, the representation of the Keller--Segel chemotaxis involves identifying two real functions, denoted as $\rho = \rho(\bm{x}, t)$ and $c = c(\bm{x}, t)$,  subject to the following dimensionless system
\begin{subequations}\label{model:ks}
	\begin{numcases}{}
  \p_t \rho  = \Delta \rho-\lambda \nabla \cdot(\rho \nabla c), & in $\Omega \times(0, T]$, \label{model:ks:a}\\
    \p_t c  = \Delta c- c+  \rho, & in $\Omega \times(0, T]$, \label{model:ks:b}
\end{numcases}
\end{subequations}
with homogeneous Neumann boundary conditions
\begin{equation}\label{model:ks:bc}
\begin{aligned}
	&\p_{\bm{n}} \rho  :=\nabla \rho \cdot \bm{n}=0, \quad \p_{\bm{n}} c=0, & & \text { on } \p \Omega \times(0, T], 
\end{aligned}
\end{equation}
and initial conditions
\begin{equation}\label{model:ks:ic}
\begin{aligned}
	&\rho(\bm{x},0)  =\rho^o(\bm{x}), \quad c(\bm{x},0)=c^o(\bm{x}),  & & \text { in } \Omega.
\end{aligned}
\end{equation}
Here $\Omega \subset \mR^2$ is a two-dimensional convex, bounded and open domain, $\bm{n}$ represents the unit outer normal vector onto the boundary, and $\lambda>0$ is a chemotactic sensitivity constant. 

In Biology, $\rho$ denotes the cell (or organism) density at position $\bm{x} \in \Omega$ and time $t \in [0, T]$, and $c$ represents the concentration of the chemoattractant, a chemical signal that induces cell migration. Both cells and chemoattractant undergo diffusion within the spatial domain, auto-diffusion phenomena observed in cells and chemoattractants are characterized by the terms $- \Delta \rho$ and $- \Delta c$, respectively, while the migration mechanism is represented by the nonlinear cross-diffusion term $- \lambda\nabla \cdot (\rho \nabla c)$. Actually, the nonlinear term presents a significant challenge for both the theoretical analysis and the numerical modeling of the system described by \eqref{model:ks}. In addition, the degradation and production of the chemoattractant are represented by the terms $- c$ and $ \rho$, respectively.
It is important to note that the production of chemoattractant by cells, to which the cells are attracted, can lead to the phenomenon known as chemotactic collapse. This involves uncontrolled aggregation that can result in a sudden increase, or explosion, within a finite time. This phenomenon is recognized as a key feature of the classical Keller--Segel model, presenting a significant challenge, particularly for numerical methods. The classical parabolic-parabolic type Keller--Segel model \eqref{model:ks} obeys the mass conservation law under the homogeneous Neumann boundary conditions, i.e.,
\begin{equation}\label{model:ks:MC}
\begin{aligned}
Mass_\rho(t) := \int_\Omega \rho(\bm{x}, t)  d\bm{x} = \int_\Omega \rho^o(\bm{x})  d\bm{x}=Mass_\rho(0).
\end{aligned}
\end{equation}
When the initial data $\rho^o$ has a total mass strictly below $8 \pi /\lambda$, a continuous solution is guaranteed for all times in a bounded domain $\Omega$. Conversely, if the initial data has a total mass exceeding $8 \pi /\lambda$ while maintaining a finite second moment, the solution exhibits a finite-time blow-up phenomenon, as proved in \cite{CC'08, LWZ'18}.

Due to  difficulties in obtaining the analytic solutions to the Keller–Segel equations, extensive investigations have been conducted on numerical modeling and analysis. For instance, Li et al. \cite{LSY'17} employed a local discontinuous Galerkin method to obtain optimal convergence rates based on a special finite element space before the blow-up occurs, and a positivity-preserving limiter is developed to model the blow-up time; 
Xiao et al. \cite{XFH'19} developed a semi-implicit characteristic finite element method (FEM) to simulate the blow-up scenario of the chemotaxis system on surface, as well as pattern formation and aggregation phenomena of the bacteria. Their method can achieve second-order accuracy in both $L^2$- and $H^1$-norm errors; In 2019, Sulman and Nguyen \cite{SN'19} introduced an adaptive moving mesh implicit-explicit FEM, and they demonstrated the effectiveness of the non-uniform spatial grids; More recently, Shen et al. \cite{SX'20} and Huang et al. \cite{HS'21} focused on developing positivity/bound preserving and unconditionally energy stable temporal schemes for the Keller--Segel system, while it is not obvious to give a concrete error analysis; Additionally, some other numerical methods, such as operator splitting methods \cite{TSL'19, M'03}, hybrid finite-volume-finite-difference methods \cite{CE'18}, finite difference methods \cite{S'09, EK'12, HZ'23}, and meshless methods \cite{BGG'20, DA'19}, have also been developed for the chemotaxis system. However, up to now, there has been no consideration of finite difference methods specifically designed on non-uniform spatial grids for the chemotaxis model, let alone any theoretical numerical analysis.  

It is well-known that the usage of non-uniform grid numerical methods can enhance the accuracy and efficiency of simulations, by distributing more grid points at positions where the solution has a large deformation and less grid points at positions where the solution changes rather slowly. Thus, 
the possibility of finite-time blow-up in the chemotaxis model \eqref{model:ks} motivates us to develop an efficient finite difference method on non-uniform grids to better simulate this phenomenon. As is well known, the block-centered finite difference (BCFD) method, sometimes called cell-centered finite difference method \cite{AWY'97}, can gain second-order spatial accuracy on non-uniform grids, without any accuracy lost compared to standard finite difference methods. Thus, it is widely considered in the literature for various models, see \cite{RP'12, XXF'22, WXF'24, SXLF'21, LSR'19, RL'15}. As far as we know, there is no published work on BCFD method for the Keller--Segel model.
Inspired by \cite{LSR'19}, where a SAV-BCFD method is considered for gradient flows for the first time, and an $L^2$-norm error analysis is rigorously proved on uniform spatial grids, we aim to introduce a mass-conservative Crank-Nicolson type BCFD scheme on non-uniform spatial grids, which is linearly and decoupled from the computational point of view. However, the analysis techniques in \cite{LSR'19} cannot be directly applied to non-uniform grids and the chemotaxis system. Incorporating high-order small perturbation techniques, along with the mathematical induction and discrete energy methods, a rigorous optimal second-order error analysis for both the cell density and chemoattractant concentration is conducted on non-uniform spatial grids, and mass conservation of the scheme is also proved. To the best of the authors’ knowledge, this seems to be the first paper with rigorous proof of second-order convergence for a linear decoupled finite difference scheme on general non-uniform grids for the Keller--Segel system \eqref{model:ks}. In summary, our new scheme enjoys the following remarkable advantages:
\begin{itemize}
    \item \emph{Mass conservation}:  The proposed scheme can ensure the mass conservation law, thereby can accurately enhance the physical realism of the simulations.
    
    \item \emph{High efficiency}: The decoupling of solutions for the cell density and the chemoattractant concentration with non-uniform spatial grids can greatly improve the efficiency of numerical simulations.

    \item \emph{Optimal-order convergence on non-uniform spatial grids}: A significant advantage of this method is its ability to achieve second-order accuracy in both time and space, even on non-uniform spatial grids. Rigorous theoretical analysis and numerical simulations on non-uniform spatial grids are carried out. 
    
    \item \emph{Blow-up simulation}: This method can accurately simulate various blow-up phenomena by using specified non-uniform spatial grids. Although it cannot be proved that the positivity of the cell density is strictly preserved for the developed scheme, well-chosen non-uniform spatial grids help reducing any non-positivity violations.
\end{itemize}
The remainder of this paper is organized as follows: In Section \ref{sec:not}, we introduce some key notations and preliminary results. Section \ref{sec:mc-bcfd} presents the construction of a mass-conservative linearly decoupled BCFD scheme along with a rigorous error analysis on non-uniform spatial grids. In Section \ref{sec:num}, we conduct some numerical experiments to support our theoretical findings, and a three-dimensional Keller–Segel system is also tested. Finally, we summarize our results and discuss future research plans. Throughout this paper, we denote by $K$ a generic positive constant that is independent of the grid parameters, but may have different values in different occurrences.

\section{Notations and preliminaries}\label{sec:not}
 Let $N$ be a positive integer and $t^n=n\tau$ ($0\le n\le N$) with $\tau:=T/N$ be a given sequence.
 For temporal grid function  $\{\phi^n\}$, define
\begin{align*}
  d_t\phi^{n} = \f{\phi^n-\phi^{n-1}}{\tau},~~ \phi^{n+1/2} = \f{\phi^{n+1} + \phi^{n}}{2}.
\end{align*}

For simplicity, we consider $\Omega :=(a^x, b^x)\times (a^y, b^y)$. Let $N_x$ and $N_y$ be the number of grids along the $x$ and $y$ coordinates, respectively. Similar to those used in \cite{WW'88,RP'13}, staggered spatial grids are introduced, where the primal grid points $\Pi_x \times \Pi_y$ are denoted by
\begin{align*}
\Pi_x: \quad a^x=x_{1 / 2}<x_{3 / 2}< \ldots < x_{i-1 / 2}<x_{i+1 / 2} < \ldots <  x_{N_x+1 / 2}=b^x,\\
\Pi_y: \quad a^y=y_{1 / 2}<y_{3 / 2}< \ldots < y_{j-1 / 2}<y_{j+1 / 2} < \ldots <  y_{N_y+1 / 2}=b^y,
\end{align*}
with grid sizes $\Delta x_i=x_{i+1/2}-x_{i-1/2},\, i=1, \ldots, N_x$ and 
$\Delta y_j=y_{j+1/2}-y_{j-1/2},\, j=1, \ldots, N_y$, and the auxiliary grid points are denoted by
\begin{align*}
	 \Pi_x^*: \quad x_i= (x_{i-1/2}+x_{i+1/2} ) / 2,  ~~ i=1, \ldots, N_x,\\
     \Pi_y^*: \quad y_j= (y_{j-1/2}+ y_{j+1/2} ) / 2,  ~~ j=1, \ldots, N_y,
\end{align*}
with grid sizes  
$\Delta x_{i+1/2}=x_{i+1}-x_i= (\Delta x_{i+1} +\Delta x_i)/2$,\, $i=1,\ldots, N_x-1$ and 
$\Delta y_{j+1/2}=y_{j+1}-y_j= (\Delta y_{j+1} +\Delta y_j)/2$,\, $j=1,\ldots, N_y-1$. 
We use $\Delta x:=\max_{i} \Delta x_i$ and $\Delta y:=\max_{j} \Delta y_j$ to represent the maximum grid sizes along $x$ and $y$ directions. Let $h:=\max \{\Delta x, \Delta y\}$ and assume that the grid partition is regular,
i.e., there exists a positive constant $\sigma$ such that
$ 	\max \{\f{h}{\min_{i} \Delta x_i}, \f{h}{\min_{j} \Delta y_j} \} \le \sigma.
$ 
Moreover, given spatial grid functions $g=\{g_{i, j}\}$, $\hat{g}=\{g_{i+1/2, j}\}$ and $\check{g}=\{g_{i, j+1/2}\}$ respectively defined on $\Pi_x^*\times \Pi_y^*$, $\Pi_x\times \Pi_y^*$ and $\Pi_x^*\times \Pi_y$,  we define
\begin{align*}
	& [d_x g]_{i+1/2, j}=\f{g_{i+1, j}-g_{i, j}}{\Delta x_{i+1/2}}, ~~	
        && [d_y g]_{i, j+1/2}=\f{g_{i, j+1}-g_{i, j}} {\Delta y_{j+1/2}}, \\
	& [D_x \hat{g}]_{i, j}=\f{\hat{g}_{i+1/2, j}-\hat{g}_{i-1/2, j}}{\Delta x_i}, ~~
	&&  [D_y \check{g}]_{i, j}=\f{\check{g}_{i, j+1/2}-\check{g}_{i, j-1/2}}{\Delta y_j}, 
\end{align*}
and denote $\bm{d}g = \left(d_xg, d_yg\right)$. Besides, we introduce the discrete inner products and norms on $\Pi_x^*\times \Pi_y^*$, $\Pi_x\times \Pi_y^*$ and $\Pi_x^*\times \Pi_y$ as follows:
\begin{align*}
	& (f, g)_{\rm M}=\sum_{i=1}^{N_x} \sum_{j=1}^{N_y} \Delta x_i \Delta y_j f_{i, j} g_{i, j},  &&\|f\|_{\rm M}^2 = (f, f)_{\rm M}, \\
	& (f, g)_x=\sum_{i=1}^{N_x-1} \sum_{j=1}^{N_y} \Delta x_{i+1/2} \Delta y_j f_{i+1/2, j} g_{i+1/2, j},  &&\|f\|_x^2 = (f, f)_x, \\
	& (f, g)_y=\sum_{i=1}^{N_x} \sum_{j=1}^{N_y-1} \Delta x_i \Delta y_{j+1/2} f_{i, j+1/2} g_{i, j+1/2},  &&\|f\|_y^2 = (f, f)_y,\\
	& (\bm {f}, \bm {g})_{\rm TM}  = (f^x, g^x )_x+ (f^y, g^y)_y,    &&\|\bm {f}\|_{\rm TM}^2 = (\bm{f}, \bm{f})_{\rm TM},\\
        & \|d_x g\|_\infty =\max _{i, j} |[d_x g]_{i+1/2, j}|,~\|d_y g\|_\infty =\max _{i, j} |[d_y g]_{i, j+1/2}|, &&\|\bm{d}g\|_\infty= \|d_x g\|_\infty +\|d_y g\|_\infty.
\end{align*}

Next, we present some useful lemmas that will play an important role in the subsequent error analysis. 
\begin{lemma}\label{lem:trunerr:t} 
\begin{enumerate}[(i)] \item If $f'(t) \in L^2\left([t^0, t^1]\right)$,   then we have
\begin{equation*}
		f(t^1)
		= f(t^0) + \Upsilon_1, \quad |\Upsilon_1| \leq K\|f'\|_{L^{2}([t^0, t^1])}\tau^{1/2}.
\end{equation*}

   \item If $f''(t) \in L^{2}([t^{n-1}, t^n])$,   then we have
\begin{align*}
		f(t^{n-1/2})
		&= \f{1}{2}\left( f(t^{n}) + f(t^{n-1}) \right)+ \Upsilon_2, \quad |\Upsilon_2| \leq K\|f''\|_{L^{2}([t^{n-1}, t^n])}\tau^{3/2}.
\end{align*}

   \item If $f''(t) \in L^2([t^{n-2}, t^n])$, then we have
\begin{align*}
	f(t^{n-1/2})
	= \f{1}{2}\left( 3f(t^{n-1}) - f(t^{n-2}) \right)  +\Upsilon_3, \quad |\Upsilon_3| \leq K\|f''\|_{L^{2}([t^{n-2}, t^n])}\tau^{3/2}.
\end{align*}
\end{enumerate}
\end{lemma}
\begin{lemma}\label{lem:trunerr:dt} 
\begin{enumerate}[(i)]
     \item If $f''(t) \in L^2([t^{0}, t^1])$,  then  we have
\begin{equation*}
	f^{\prime}(t^1)
	= \f{f(t^1)-f(t^0)}{\tau} + \Theta_1, \quad  |\Theta_1| \leq K\|f''\|_{L^{2}([t^{0}, t^1])}\tau^{1/2}.
\end{equation*}
 
     \item If $f'''(t) \in L^2([t^{n-1}, t^n])$,  then  we have
\begin{align*}
	&f^{\prime}(t^{n-1/2})
	= \f{f(t^{n})-f(t^{n-1})}{\tau}  + \Theta_2, \quad |\Theta_2| \leq K\|f'''\|_{L^{2}([t^{n-1}, t^n])}\tau^{3/2}.
\end{align*}
\end{enumerate}
\end{lemma}
\begin{remark} The conclusions in Lemmas \ref{lem:trunerr:t}–\ref{lem:trunerr:dt} can be similarly proved as those in Lemma 1.2 of Ref. \cite{SZG'23}, where the Lagrange remainders of Taylor expansion are replaced by the integral form, and the Cauchy-Schwarz inequlity is applied. Detailed discussions are omitted. 
\end{remark}
\begin{lemma}  \label{trun:Dp}
Let $p\in H^3(\Omega)$, then there holds
\begin{equation*}
\left\{
 \begin{aligned} 
     \p_x p (x_{i}, y_j)= [D_x p ]_{i, j}      + \vartheta^x_{i,j}(p),\\
     \p_y p (x_i, y_j) = [D_y p ]_{i, j}   +\vartheta^y_{i,j}(p),
\end{aligned}\right.
\end{equation*}
such that
\begin{equation*}
\|\vartheta^x(p)\|_{\rm M}  \leq K  \|p\|_{H^{3}(\Omega)}h^2,
\quad \|\vartheta^y(p)\|_{\rm M}  \leq K \|p\|_{H^{3}(\Omega)} h^2.
\end{equation*}
\end{lemma}
 \begin{proof}
 The conclusion is a direct result of the application of Bramble-Hilbert lemma \cite{BH'71} to the truncation errors $\vartheta^x(p)$ and $\vartheta^y(p)$, see also Lemma 3.4 of Ref. \cite{BGM'07}  for details.
 \end{proof}
\begin{lemma} \label{trun:ap-dp} 
Let  $p\in W^{3,\infty}(\Omega)$, then there holds
\begin{equation*}
\left\{
 \begin{aligned} 
     \p_x p (x_{i+1 / 2}, y_j)= \big[d_x(p-\delta(p))\big]_{i+1/2, j}
                                 +\epsilon_{i+1/2, j}^x(p), \\
     \p_y p (x_i, y_{j+1 / 2}) =\big[d_y(p-\delta(p))\big]_{i, j+1 / 2}
                                 +\epsilon_{i, j+1 / 2}^y(p),
\end{aligned}
\right.
\end{equation*}
with
\begin{align*}
	\delta_{i, j}(p)& :=
	\f{\Delta x_i^2}{8}  \p_{xx} p  (x_i, y_j)+\f{\Delta y_j^2}{8} \p_{yy} p (x_i, y_j),\\
	\epsilon_{i+1 / 2, j}^x(p)& 
    :=  \f{1}{2 \Delta x_{i+1 / 2}} \int_{x_{i+1 / 2}}^{x_{i+1}}\left(\f{  \Delta x_{i+1}^2}{4}-\left(x-x_{i+1}\right)^2\right) \p_{xxx} p(x, y_j) d x \\
	& \qquad  +\f{1}{2 \Delta x_{i+1 / 2}} \int^{x_{i+1 / 2}}_{x_i}\left(\f{\Delta x_i^2}{4}-\left(x-x_i\right)^2\right) \p_{xxx} p(x, y_j) d x \\
	& \qquad +\f{\Delta y_j^2}{8 \Delta x_{i+1 / 2}} \int_{x_i}^{x_{i+1}} \p_{xyy}p(x, y_j) dx,\\
	\epsilon_{i, j+1 / 2}^y(p)& :=  \f{1}{2 \Delta y_{j+1 / 2}} \int_{y_{j+1 / 2}}^{y_{j+1}}\left(\f{\Delta y_{j+1}^2}{4}-\left(y-y_{j+1}\right)^2\right) \p_{yyy} p\left(x_i, y\right) dy \\
	& \qquad +\f{1}{2 \Delta y_{j+1 / 2}} \int^{y_{j+1 / 2}}_{y_j}\left(\f{\Delta y_j^2}{4}-\left(y-y_j\right)^2\right) \p_{yyy} p\left(x_i, y\right) d y \\
	& \qquad +\f{\Delta x_i^2}{8 \Delta y_{j+1 / 2}} \int_{y_j}^{y_{j+1}} \p_{xxy} p(x_i, y) dy,
\end{align*}
such that 
\begin{equation*}
\begin{aligned}
\|\delta(p)\|_{\rm M} \leq K \|p\|_{W^{2,\infty}(\Omega)} h^2=: K_\delta h^2,  \quad
\|\epsilon^x(p)\|_x +\|\epsilon^y(p)\|_{y} \leq K  \|p\|_{W^{3,\infty}(\Omega)} h^2=:K_\epsilon h^2.
\end{aligned}
\end{equation*}
\end{lemma}
 \begin{proof} 
 	The lemma is proved in Ref. \cite{RP'13} by using the Taylor expansion under assumption
 $p\in W^{3,\infty}(\Omega)$, see Lemmas 4.1--4.2 of \cite{RP'13} for details. In fact, similar to Lemma \ref{trun:Dp}, the estimates in Lemma \ref{trun:ap-dp} can also be obtained via the Bramble-Hilbert lemma under a weaker regularity assumption $p\in H^{3}(\Omega)$. 
 \end{proof}
\begin{lemma} \label{lemma:Dd} 
 Let $\{q_{i, j}\}, \{v_{i+1/2, j}\}$ and $\{w_{i, j+1 / 2}\}$ be any grid functions defined on $\Pi_x^*\times \Pi_y^*$, $\Pi_x\times \Pi_y^*$ and $\Pi_x^*\times \Pi_y$, such that $v_{1/2, j}=$ $v_{N_x+1/2, j}=w_{i, 1/2}=w_{i, N_y+1 / 2}=0$. Then there holds
$$
\begin{aligned}
 (q, D_x v)_{\rm M}=-(d_x q, v)_x, ~~ (q, D_y w)_{\rm M}=-(d_y q, w)_y .
\end{aligned}
$$
\end{lemma}
\begin{proof}
  Please refer to Lemma 4.2 of Ref. \cite{WW'88} for detail.
\end{proof}
Finally, for given function values $\{p_{i,j} = p(x_i,y_j)\}$, we introduce the piecewise bilinear interpolant operator $\ell_h p$, which is similar to that used in \cite{DWW'98}. For any points $(x, y) \in [x_i, x_{i+1}]\times [y_j, y_{j+1}]$, $i=1,\ldots, N_x$, $j=1,\ldots, N_y$, we define  $\ell_h p(x, y)$ for approximating $p(x, y)$ by
 \begin{equation}\label{Operator:bi}
\begin{aligned}
	\ell_h p(x, y)& =    \f{(x_{i+1}-x)(y_{j+1}-y)}{(x_{i+1}-x_i)(y_{j+1}-y_j)}p_{i, j} 
           + \f{(x-x_i)(y_{j+1}-y)}{(x_{i+1}-x_i)(y_{j+1}-y_j)} p_{i+1, j} \\
	& \quad + \f{(x_{i+1}-x)(y-y_j)}{(x_{i+1}-x_i)(y_{j+1}-y_j)}  p_{i, j+1} 
         + \f{(x-x_i)(y-y_j)}{(x_{i+1}-x_i)(y_{j+1}-y_j)} p_{i+1, j+1}.
\end{aligned}
\end{equation}

\begin{lemma}\label{lem:interp} 
Assume that $p \in W^{2, \infty}(\Omega)$, then we have 
\begin{align*}
\left\|\ell_h p-p\right\|_{L^\infty} \leq K  \|p\|_{W^{2, \infty}(\Omega)}h^2.
\end{align*}
\end{lemma}
\section{Numerical scheme and its analysis}\label{sec:mc-bcfd}
In this section, a decoupled and mass-conservative Crank-Nicolson type BCFD algorithm, abbreviated as DeC-MC-BCFD, is proposed for the Keller--Segel chemotaxis system \eqref{model:ks}, in which the Crank-Nicolson formulae and linear extrapolation technique are used in time discretization and block-centered finite difference method is considered in space discretization on general non-uniform staggered spatial grids. We shall provide detailed analysis on demonstrating the mass conservation law and error estimates for the developed algorithm.

\subsection{The DeC-MC-BCFD scheme and its mass conservation}

Throughout the paper, we make the following smoothness assumptions
\begin{equation}\label{Regularity}
\begin{aligned}
&\rho\in L^{\infty}([0,T] ; H^{4}(\Omega)\cap W^{3,\infty}(\Omega)), \quad c \in L^{\infty}([0,T] ; H^{4}(\Omega)), \\
 & \rho_t\in L^{4}([0,T] ; W^{1,\infty}(\Omega)), \quad   c_{t}\in L^{4}([0,T] ; W^{2,\infty}(\Omega)),\\
 &\rho_{tt},~c_{tt}\in L^{2}([0,T] ; W^{2,\infty}(\Omega)),
 \quad \rho_{ttt},~c_{ttt}\in L^{2}([0,T] ; L^{\infty}(\Omega)).
\end{aligned}
\end{equation}
Furthermore, we assume there exists a positive constant $K_*$ such that
\begin{equation}\label{bound}
  \| \nabla c \|_{L^{\infty}([0,T]\times\Omega)}\leq K_*.
\end{equation}

At each time level $t^n=n\tau$ ($0\le n\le N$), let us denote the approximations of $\left\{\rho^n, c^n\right\}$ by $\left\{U^n, Z^n\right\}$, in which $U^n=\{U_{i,j}^n\}$ and $Z^n=\{Z_{i,j}^n\}$ are defined on $\Pi_x^*\times \Pi_y^*$. Noting that the main difficulty in the construction of the DeC-MC-BCFD scheme lies in the discretization of the nonlinear cross-diffusion term $- \lambda\nabla \cdot (\rho \nabla c)$, as the two components of $\nabla c$ are approximated on different staggered spatial grids $\Pi_x\times \Pi_y^*$ and $\Pi_x^*\times \Pi_y$, respectively, while the variable $\rho$ is discretized on $\Pi_x^*\times \Pi_y^*$. Thus, we shall use the piecewise bilinear interpolant operator $\ell_h$ defined by \eqref{Operator:bi} to make $\rho$ and $\nabla c$ match on the same grid points. Then, the DeC-MC-BCFD scheme for the Keller--Segel system \eqref{model:ks} is proposed as follows. 
 
 {\bf Step 1}: Solve $(Z^1,U^1)$ via the prediction-correction BCFD scheme
 \begin{subequations}\label{2D:PC-BCFD}
	\begin{numcases}{}
		d_t \bar{U}_{i,j}^{1} = [D_x (d_x \bar{U})]_{i,j}^{1} + [D_y (d_y \bar{U})]_{i,j}^{1} \nonumber\\
		\qquad\qquad\qquad
             -\lambda \left( [D_x ([\ell_h \bar{U}^1] [d_x Z^0])]_{i,j} 
                           + [D_y ([\ell_h \bar{U}^1] [d_y Z^0])]_{i,j}\right) ,  \label{2D:PC-BCFD:rho}	\\
  d_t Z_{i,j}^{1} =  [D_x (d_x Z)]_{i,j}^{1/2} + \f{1}{2} [D_y (d_y Z^1)]_{i,j}^{1/2}
 	-  Z_{i,j}^{1/2} +  \bar{U}_{i,j}^{1/2} , \label{2D:PC-BCFD:c}\\
  		d_t U_{i,j}^{1} = [D_x (d_x U)]_{i,j}^{1/2} + [D_y (d_y U)]_{i,j}^{1/2} \nonumber\\
           \qquad\qquad\qquad -\lambda \left( [D_x ([\ell_hU] [d_x Z])]_{i,j}^{1/2} 
                   + [D_y ([\ell_hU] [d_y Z])]_{i,j}^{1/2}\right),   \label{2D:PC-BCFD:rho:e2}
  \end{numcases}
\end{subequations}
 for $1 \leq i \leq N_x, 1 \leq j \leq N_y$, where $d_t \bar{U}^{1}= (\bar{U}^1-U^0)/ \tau$, $\bar{U}^{1/2}=(\bar{U}^{1}+ U^{0})/2$.

 {\bf Step 2}: Solve $\{Z^{n+1}, U^{n+1}\}_{n=1}^{N-1}$ via the Crank-Nicolson BCFD scheme
\begin{subequations}\label{2D:CN-BCFD}
	\begin{numcases}{}
 	d_t Z_{i,j}^{n+1} = [D_x (d_x Z)]_{i,j}^{n+1/2} + [D_y (d_y Z)]_{i,j}^{n+1/2}
 	-  Z_{i,j}^{n+1/2} + U_{i,j}^{*,n+1}. \label{2D:CN-BCFD:c} \\
		d_t U_{i,j}^{n+1} = [D_x (d_x U)]_{i,j}^{n+1/2} + [D_y (d_y U)]_{i,j}^{n+1/2} \nonumber\\
		\qquad\qquad\qquad -\lambda \left( [D_x ([\ell_hU] [d_x Z])]_{i,j}^{n+1/2} + [D_y ([\ell_hU] [d_y Z])]_{i,j}^{n+1/2}\right) ,  \label{2D:CN-BCFD:rho}	
  \end{numcases}
\end{subequations}
for $1 \leq i \leq N_x, 1 \leq j \leq N_y$, where $U^{*,n+1} = 3U^n/2-U^{n-1}/2$ is an explicit second-order temporal approximation of $\rho^{n+1/2}$.

Steps 1--2 are enclosed with the following boundary and initial conditions 
 \begin{equation}\label{2D:CN-BCFD:IBc}
     \left\{
 	\begin{aligned}
		&\left[d_x \bar{U}\right]_{i+1 / 2, j}^1 = \left[d_x Z\right]_{i+1/2, j}^n=\left[d_x U\right]_{i+1 / 2, j}^n
              =0, \quad i=\{0, N_x\}, ~ 1 \leq j \leq N_y, \\
		&\left[d_y \bar{U}\right]_{i+1 / 2, j}^1 = \left[d_y Z\right]_{i, j+1 / 2}^n =\left[d_y U\right]_{i, j+1 / 2}^n
             =0, \quad 1 \leq i \leq N_x, ~ j=\{0, N_y\},\\	
		&Z_{i,j}^0 = c^o(x_i,y_j) - \delta_{i, j}(c^o),  \quad 	U_{i,j}^0 = \rho^o(x_i,y_j), \quad 1 \leq i \leq N_x, 1 \leq j \leq N_y.
	\end{aligned}
 \right.
\end{equation}

\begin{remark}  Note that the DeC-MC-BCFD scheme \eqref{2D:PC-BCFD}--\eqref{2D:CN-BCFD:IBc} is linear and decoupled.
In practical computation, on the first time level, we first compute $\bar{U}^1$ through the first-order Euler prediction step \eqref{2D:PC-BCFD:rho}, then determine $Z^1$ via the second-order semi-implicit Crank-Nicolson scheme \eqref{2D:PC-BCFD:c} with known $\bar{U}^1$, and finally perform a correction step to obtain $U^1$ through the second-order fully-implicit Crank-Nicolson scheme \eqref{2D:PC-BCFD:rho:e2}. Then, for $n \geq 1$, as $U^n$ and $U^{n-1}$ are known, we solve $Z^{n+1}$ via the linear scheme \eqref{2D:CN-BCFD:c}, and subsequently, with the obtained $Z^{n+1}$, we solve $U^{n+1}$ via the linear scheme \eqref{2D:CN-BCFD:rho}. This sequential approach facilitates linearization and decoupling, and thereby can greatly simplify the computation and thus enhance the computational efficiency of the method, in which each variable is solved independently.
\end{remark}

\begin{remark} Even without the correction step \eqref{2D:PC-BCFD:rho:e2} in Step I, the scheme \eqref{2D:PC-BCFD}--\eqref{2D:CN-BCFD:IBc} can also be theoretically proven and numerically validated to achieve second-order convergence in both time and space. However, this step can further enhance the accuracy of the numerical solutions, without spending too much computational cost.
\end{remark}

\begin{theorem}[Discrete Mass Conservation]\label{thm:MassConserve}
Let $\left\{Z^{n+1}, U^{n+1}\right\}$ be the solutions to the DeC-MC-BCFD scheme \eqref{2D:PC-BCFD}--\eqref{2D:CN-BCFD:IBc}. Then, there holds
\begin{equation*}
    \left(U^{n+1}, \bm 1 \right)_{\mathrm{M}} = \left(U^0,  \bm 1\right)_{\mathrm{M}}.
\end{equation*}
\end{theorem}
\begin{proof}
By taking inner product on both sides of equations \eqref{2D:PC-BCFD:rho:e2} and \eqref{2D:CN-BCFD:rho}  with $\bm 1$, for $n\geq0$, we have
 \begin{align*}
	\left(d_t U^{n+1}, \bm 1\right)_{\mathrm{M}} &= \left([D_x (d_x U)]^{n+1/2} + [D_y (d_y U)]^{n+1/2},  \bm 1\right)_{\mathrm{M}} \\
	&\qquad -\lambda \left( \left[D_x ([\ell_hU] [d_x Z])\right]^{n+1/2} + \left[D_y ([\ell_hU] [d_y Z])\right]^{n+1/2}, \bm 1\right)_{\mathrm{M}}.  	
\end{align*}

Note that Lemma \ref{lemma:Dd} and the discrete boundary conditions \eqref{2D:CN-BCFD:IBc} show that
\begin{equation}\label{mass：laplace}
\begin{aligned}
 & \left( [D_x (d_x U)]^{n+1/2} + [D_y (d_y U)]^{n+1/2},  \bm 1\right)_{\mathrm{M}} \\
&\qquad =-\left(d_x U^{n+1/2} , d_x  \bm 1\right)_x-\left(d_y U^{n+1/2}, d_y \bm 1\right)_y = 0,
\end{aligned}
\end{equation}
and furthermore, we similarly get
\begin{equation}\label{mass：cross}
\begin{aligned}
    \left( \left[D_x ([\ell_hU] [d_x Z])\right]^{n+1/2} + \left[D_y ([\ell_hU] [d_y Z])\right]^{n+1/2},  \bm 1\right)_{\mathrm{M}} = 0.
\end{aligned}
\end{equation}
Thus, we obtain
\begin{equation}\label{mass：covt}
\left(d_t U^{n+1},  \bm 1\right)_{\mathrm{M}} = 0 \Longrightarrow \left(U^{n+1},  \bm 1\right)_{\mathrm{M}}=\left(U^{n},  \bm 1\right)_{\mathrm{M}},
\end{equation}
which implies the conclusion. 
\end{proof}

\subsection{Truncation errors}\label{subsec:Truc}
In this subsection, we analyze the truncation errors of the developed DeC-MC-BCFD scheme \eqref{2D:PC-BCFD}--\eqref{2D:CN-BCFD:IBc}. For simplicity, let
\begin{align*}
&\rho_{i,j}(t) := \rho(x_i,y_j, t), ~~c_{i,j}(t) :=c(x_i,y_j, t), 
 ~~\rho_{i,j}^n:=\rho_{i,j} (t^n ), ~~c_{i,j}^n:=c_{i,j} (t^n ), 
\end{align*}
for $1 \leq i \leq N_x, 1 \leq j \leq N_y$ and  $0 \leq n \leq N$. Whenever no confusion caused, we usually omit the subscripts.

\paragraph{\indent \bf Analysis of the truncation errors for \eqref{2D:PC-BCFD}} In the first step,  a prediction-correction approach is proposed to construct the decoupled, second-order accurate in time scheme for the first time level, and therefore, we shall analyze the truncation errors separately. 

First, \eqref{2D:PC-BCFD:rho} can be viewed as discretization of \eqref{model:ks:a} using the semi-implicit Euler-BCFD method.
It then follows from Lemma \ref{trun:ap-dp} that the exact solution $\rho^1$ satisfies 
\begin{equation}\label{exact:PC-BCFD:rho1} 
\begin{aligned}
		 d_t \rho_{i,j}^{1} 
		&= \big [ D_x  ( d_x (\rho-\delta(\rho)) +  \epsilon^x(\rho) ) \big]_{i,j}^{1} 
		+ \big[ D_y ( d_y (\rho-\delta(\rho)) +  \epsilon^y(\rho) ) \big]_{i,j}^{1} \\
		&\quad  -\lambda  \left(\big[D_x ( [\ell_h\rho]^1 [d_x Z^0 +  \epsilon^x(c^o)] ) \big]_{i,j} 
	          + \big[D_y ( [\ell_h\rho]^1 [d_y Z^0 +  \epsilon^y(c^o)] ) \big]_{i,j}  \right)
        + \mR^1_{i,j},  
\end{aligned}
\end{equation}
where $\delta(\cdot)$, $\epsilon^x(\cdot)$ and $\epsilon^y(\cdot)$ are defined in Lemma \ref{trun:ap-dp}, and the truncation error 
\begin{equation}\label{exact:PC-BCFD:rho1:e1} 
\begin{aligned}
	 \mR^1 & :=\big(d_t \rho^{1}  - \p_t \rho(t^{1}) \big)
                    + \big(\Delta \rho(t^{1}) - \left[ D_x\p_x \rho + D_y\p_y \rho\right]^{1}\big)\\
                &\quad  + \lambda \left(  D_x ( [\ell_h\rho]^1 \p_x c^o ) -  \left[ \p_x( \rho\p_x c )\right](t^{1})
	           +   D_y ( [\ell_h\rho]^1  \p_y c^o ) -  \left[ \p_y ( \rho \p_y c )\right](t^{1})\right)\\
           & =\big(d_t \rho^{1}  - \p_t \rho(t^{1}) \big)
                    + \big(\Delta \rho(t^{1}) - \left[ D_x\p_x \rho + D_y\p_y \rho\right]^{1}\big)\\
         &\qquad + \lambda \big( D_x (\rho^1 \p_x c^o) - \p_x (\rho^1 \p_x c^o) + D_y (\rho^1 \p_x c^o) - \p_y (\rho^1 \p_y c^o)\big) \\
         &\qquad + \lambda \big( \p_x (\rho^1 (\p_x c^o - \p_x c^1 )) + \p_y (\rho^1 (\p_y c^o - \p_y c^1 ) \big)\\
          &\qquad +  \lambda \big(D_x ([\ell_h\rho-\rho]^1 \p_x c^o) +  D_y ([\ell_h\rho-\rho]^1 \p_y c^o) \big).
\end{aligned}
\end{equation}
Below we use
$
\mR_{II}^1  :=   \lambda \big(D_x ([\ell_h\rho-\rho]^1 \p_x c^o) +  D_y ([\ell_h\rho-\rho]^1 \p_y c^o) \big)
$
to denote the last part of $\mR^1$, and other terms in $\mR^1$ are denoted as $\mR^1_{I}$.
Then, following from Lemmas \ref{lem:trunerr:t}--\ref{trun:Dp},  $\mR^1_{I}$ can be estimated as 
\begin{equation}\label{TruncErr:PC:e11} 
	\begin{aligned}
	\|\mR^1_{I}\|_{\rm M}	&\leq K \big(\|\rho_{tt}\|_{L^{2}([t^0,t^1] ; L^{\infty}(\Omega))} + \|\rho\|_{L^{\infty}([t^0,t^1]; W^{1, \infty}(\Omega) )}\|c_t\|_{L^{2}([t^0,t^1] ; {W^{2, \infty}(\Omega)} )} \big) \tau^{1/2} \\
       &\qquad  + K \big( \|\rho\|_{L^{\infty}([t^0,t^1];H^{4}(\Omega) )} + \|\rho\|_{ L^{\infty}([t^0,t^1];W^{3,\infty}(\Omega) )} \| c\|_{L^{\infty}([t^0,t^1]; H^{4}(\Omega) )} \big) h^2.
	\end{aligned}
\end{equation}
However, direct estimation of the term $\mR^1_{II}$ can not maintain the whole second-order convergence, instead there will be a loss of error accuracy. Therefore, special treatment of this term shall be included in the subsequent convergence analysis.

Next, by discretizing \eqref{model:ks:b} using the Crank-Nicolson BCFD method, it can be deduced that the exact solution $c^1$ satisfies the discretized equation:
 \begin{equation}\label{exact:CN-BCFD:c1} 
 	\begin{aligned}
 	d_t c_{i,j}^{1} 
 	&= \left[ D_x ( d_x (c-\delta(c)) +  \epsilon^x(c) ) \right]_{i,j}^{1/2} 
 	+ \left[ D_y ( d_y (c-\delta(c)) +  \epsilon^y(c) ) \right]_{i,j}^{1/2} \\
 	&\qquad-  c_{i,j}^{1/2} + \rho_{i,j}^{1/2}+  \mF^{1}_{i,j}, 
 	\end{aligned}
 \end{equation}
where the truncation error
 \begin{equation*}\label{exact:PC-BCFD:c1:e1} 
\begin{aligned}
	 \mF^1 & :=\big(d_t c^{1}  -  \p_t c(t^{1/2}) \big)
                    + \big(\Delta c(t^{1/2}) - \left[ D_x\p_x c + D_y\p_y c \right]^{1/2}\big)
                  + \big( c^{1/2} - c(t^{1/2}) \big) \\
               &\qquad   + \big( \rho(t^{1/2}) - \rho^{1/2} \big),
\end{aligned}
\end{equation*}
such that can be estimated by Lemmas \ref{lem:trunerr:t}--\ref{trun:Dp} as
\begin{equation}\label{TruncErr:PC:ec1} 
	\begin{aligned}
	\|\mF^1\|_{\rm M}	
 &\leq K \big(\|c_{ttt}\|_{L^{2}([t^0,t^1] ; L^{\infty}(\Omega))} + \|c_{tt}\|_{L^{2}([t^0,t^1] ;W^{2, \infty}(\Omega))} \big) \tau^{3/2} \\
      &\qquad   + K \|c\|_{L^{\infty}([t^0,t^1];H^{4}(\Omega) )} h^2.
	\end{aligned}
\end{equation}

Finally, applying the same Crank-Nicolson BCFD discretization to \eqref{model:ks:a} as above, it is seen that the exact solution $\rho^1$ also satisfies the following discretized equation:
\begin{equation}\label{exact:PC-BCFD:rho1c} 
	\begin{aligned}
		d_t \rho_{i,j}^{1} 
		&= \left[ D_x ( d_x (\rho-\delta(\rho)) +  \epsilon^x(\rho) ) \right] _{i,j}^{1/2} 
		+ \left[ D_y( d_y (\rho-\delta(\rho)) +  \epsilon^y(\rho) ) \right] _{i,j}^{1/2} \\
		&\qquad -\lambda  \left[D_x ( [\ell_h\rho] [ d_x (c-\delta(c)) +  \epsilon^x(c) ]) \right]_{i,j}^{1/2} \\
		&\qquad -\lambda \left[D_y ( [\ell_h\rho] [ d_y (c-\delta(c)) +  \epsilon^y(c) \right]) ]_{i,j}^{1/2}  + \mS^{1}_{i,j},  
	\end{aligned}
\end{equation}
with truncation error 
\begin{align*}
	 \mS^1 & :=\big(d_t \rho^{1}  - \p_t \rho(t^{1/2}) \big)
                    + \big(\Delta \rho(t^{1/2}) - \left[ D_x\p_x \rho + D_y\p_y \rho\right]^{1/2}\big)\\
                &\qquad  + \lambda \big(  [D_x ( [\ell_h\rho] \p_x c )]^{1/2} - [ \p_x( \rho\p_x c ) ](t^{1/2})\big)\\
	         &\qquad  + \lambda \big(   [ D_y ( [\ell_h\rho]  \p_y c ) ]^{1/2} -  [ \p_y ( \rho \p_y c )](t^{1/2}) \big)\\
           & =\big(d_t \rho^{1}  - \p_t \rho(t^{1/2}) \big)
                    + \big(\Delta \rho(t^{1/2}) - \left[ D_x\p_x \rho + D_y\p_y \rho\right]^{1/2}\big)\\
         &\qquad + \lambda \big( [D_x (\rho \p_x c) - \p_x (\rho \p_x c)]^{1/2} + [D_y (\rho \p_x c) - \p_y (\rho \p_y c)]^{1/2} \big) \\
         &\qquad + \lambda \big( [\p_x (\rho \p_x c)]^{1/2} - [ \p_x( \rho\p_x c ) ](t^{1/2}) + [\p_y (\rho\p_y c)]^{1/2} -  [ \p_y ( \rho \p_y c )](t^{1/2}) \big)\\
          &\qquad +  \lambda \big( [D_x ([\ell_h\rho-\rho]\p_x c)]^{1/2} +  D_y ([\ell_h\rho-\rho] \p_y c)^{1/2} \big) = : \mS^1_{I} + \mS^1_{II},
\end{align*}
where
$
\mS^1_{II} := \lambda \big( [D_x ([\ell_h\rho-\rho]\p_x c)]^{1/2} +  D_y ([\ell_h\rho-\rho] \p_y c)^{1/2} \big),
$
and other terms are denoted as $\mS^1_{I}$. Similar to the estimate \eqref{TruncErr:PC:e11}, we have
\begin{equation}\label{TruncErr:PC:e2} 
	\begin{aligned}
	\|\mS^1_{I}\|_{\rm M}	&\leq K \big( \|\rho_{ttt}\|_{L^{2}([t^0,t^1] ; L^{\infty}(\Omega))} + \|\rho_{tt}\|_{L^{2}([t^0,t^1] ; W^{2, \infty}(\Omega))} \\
    &\qquad + \|\rho_t\|_{L^{4}([t^0,t^1] ; W^{1,\infty}(\Omega))}\|c_t\|_{L^{4}([t^0,t^1] ; W^{2,\infty}(\Omega))} \big) \tau^{3/2} \\
       &\quad  + K \big( \|\rho\|_{L^{\infty}([t^0,t^1] ;  H^{4}(\Omega) )} + \|\rho\|_{L^{\infty}([t^0,t^1] ; W^{3,\infty}(\Omega))} \| c\|_{L^{\infty}([t^0,t^1] ; H^{4}(\Omega))} \big) h^2,
	\end{aligned}
\end{equation}
while $\mS^1_{II}$ will be analyzed in subsection \ref{subsec:err} to avoid loss of error accuracy.

\paragraph{\indent \bf Analysis of the truncation errors for \eqref{2D:CN-BCFD}} In this step, a linearized Crank-Nicolson BCFD method is developed for the discretization of \eqref{model:ks:a}--\eqref{model:ks:b} at $t=t^{n+1/2} (1\leq n\leq N-1)$. This approach ensures the maintenance of second-order convergence. It can be checked that the exact solutions $c^{n+1}$ and $\rho^{n+1}$ satisfy the discretized form:
 \begin{equation}\label{exact:CN-BCFD:c} 
 	\begin{aligned}
 	d_t c_{i,j}^{n+1} 
 	&= \left[ D_x ( d_x (c-\delta(c)) +  \epsilon^x(c) ) \right] _{i,j}^{n+1/2} 
 	+ \left[ D_y ( d_y (c-\delta(c)) +  \epsilon^x(c) ) \right] _{i,j}^{n+1/2} \\
 	&\qquad  -  c_{i,j}^{n+1/2} + \rho_{i,j}^{*,n+1} +  \mF_{i,j}^{n+1}, 
 	\end{aligned}
 \end{equation}
 \begin{equation}\label{exact:CN-BCFD:rho} 
	\begin{aligned}
		d_t \rho_{i,j}^{n+1} 
		&= \left[ D_x ( d_x (\rho-\delta(\rho)) +  \epsilon^x(\rho) ) \right] _{i,j}^{n+1/2} 
		+ \left[ D_y( d_y (\rho-\delta(\rho)) +  \epsilon^y(\rho) ) \right] _{i,j}^{n+1/2} \\
		&\qquad -\lambda  \left[D_x ( [\ell_h\rho] [ d_x (c-\delta(c)) +  \epsilon^x(c) ]) \right]_{i,j}^{n+1/2} \\
		&\qquad -\lambda \left[D_y ( [\ell_h\rho] [ d_y (c-\delta(c)) +  \epsilon^y(c) ]) \right]_{i,j}^{n+1/2}  
                  + \mS_{i,j}^{n+1},  
	\end{aligned}
\end{equation}
for $1\leq n\leq N-1$, where the truncation errors
\begin{align*}
	 \mF^{n+1} & :=\big(d_t c^{n+1}  -  \p_t c(t^{n+1/2}) \big)
                    + \big(\Delta c(t^{n+1/2}) - \left[ D_x\p_x c + D_y\p_y c \right]^{n+1/2}\big)\\
              &\qquad    + \big( c^{n+1/2} - c(t^{n+1/2}) \big) 
                  + \big( \rho(t^{n+1/2}) - \rho^{*,n+1} \big),\\
	 \mS^{n+1} & :=\big(d_t \rho^{n+1}  - \p_t \rho(t^{n+1/2}) \big)
                    + \big(\Delta \rho(t^{1/2}) - \left[ D_x\p_x \rho + D_y\p_y \rho\right]^{n+1/2}\big)\\
                &\qquad  + \lambda \big(  [D_x ( [\ell_h\rho] \p_x c )]^{n+1/2} - [ \p_x( \rho\p_x c ) ](t^{n+1/2}) \big)\\
	          &\qquad  + \lambda \big( [ D_y ( [\ell_h\rho]  \p_y c ) ]^{n+1/2} -  [ \p_y ( \rho \p_y c )](t^{n+1/2}) \big)\\
           & =\big(d_t \rho^{n+1}  - \p_t \rho(t^{n+1/2}) \big)
                   +  \big(\Delta \rho(t^{n+1/2}) - \left[ D_x\p_x \rho + D_y\p_y \rho\right]^{n+1/2}\big)\\
         &\qquad + \lambda \big( [D_x (\rho \p_x c) - \p_x (\rho \p_x c)]^{n+1/2} + [D_y (\rho \p_x c) - \p_y (\rho \p_y c)]^{n+1/2} \big) \\
         &\qquad + \lambda \big( [\p_x (\rho \p_x c)]^{n+1/2} - [ \p_x( \rho\p_x c ) ](t^{n+1/2}) \big)\\
         &\qquad + \lambda  \big( [\p_y (\rho\p_y c]^{n+1/2} -  [ \p_y ( \rho \p_y c )](t^{n+1/2}) \big)\\
          &\qquad +  \lambda \big( [D_x ([\ell_h\rho-\rho]\p_x c)]^{n+1/2} +  D_y ([\ell_h\rho-\rho] \p_y c)^{n+1/2} \big) = : \mS^{n+1}_{I} + \mS^{n+1}_{II},
\end{align*}
with $
\mS^{n+1}_{II} := \lambda \big( [D_x ([\ell_h\rho-\rho]\p_x c)]^{n+1/2} +  D_y ([\ell_h\rho-\rho] \p_y c)^{n+1/2}\big)$ representing the last term of $\mS^{n+1}$, and will be specially analyzed as $\mS^1_{II}$ and $\mR^1_{II}$ to retain the accuracy.

Similarly, it can be concluded from Lemmas \ref{lem:trunerr:t}--\ref{trun:Dp} that
\begin{equation}\label{TruncErr:CN:e1} 
	\begin{aligned}
	\|\mF^{n+1}\|_{\rm M}	
   &\leq K \big(\|c_{ttt}\|_{L^{2}([t^{n},t^{n+1}] ; L^{\infty}(\Omega))} + \|c_{tt}\|_{L^{2}([t^{n},t^{n+1}] ;W^{2, \infty}(\Omega))} \big) \tau^{3/2} \\
    &\qquad   + K \|c\|_{ L^{\infty}([t^{n},t^{n+1}]; H^{4}(\Omega) )} h^2,
	\end{aligned}
\end{equation}
and
\begin{equation}\label{TruncErr:CN:e2} 
	\begin{aligned}
	\|\mS^{n+1}_{I}\|_{\rm M}
 &\leq K \Big( \|\rho_{ttt}\|_{L^{2}([t^{n},t^{n+1}] ; L^{\infty}(\Omega))} + \|\rho_{tt}\|_{L^{2}([t^{n},t^{n+1}] ; W^{2, \infty}(\Omega))} \\
   &\qquad  +  \|\rho_t\|_{L^{4}([t^{n},t^{n+1}] ; W^{1,\infty}(\Omega))}\|c_t\|_{L^{4}([t^{n},t^{n+1}] ; W^{2,\infty}(\Omega))} \Big) \tau^{3/2} \\
       &\quad  + K \big( \|\rho\|_{L^{\infty}([t^{n},t^{n+1}] ; H^{4}(\Omega) )} + \|\rho\|_{L^{\infty}([t^{n},t^{n+1}] ; W^{3,\infty}(\Omega))} \| c\|_{L^{\infty}([t^{n},t^{n+1}] ; H^{4}(\Omega))} \big) h^2.
	\end{aligned}
\end{equation}
 
\subsection{Error estimates for the DeC-MC-BCFD scheme }\label{subsec:err}

In this subsection, we show that the DeC-MC-BCFD scheme \eqref{2D:PC-BCFD}--\eqref{2D:CN-BCFD:IBc} is second-order accurate in both time and space in corresponding discrete norms.
Set
\begin{equation*}
    \bar{e}_\rho^1 = \rho^1 - \bar{U}^1,\quad
\hat{\bar{e}}_\rho^1=\bar{e}_{\rho}^{1} - \delta(\rho^{1}),
\end{equation*}
\begin{equation*}
  e_\rho^n=\rho^n - U^n, ~ e_c^n=c^n - Z^n, ~  \hat{e}_{c}^n=e_{c}^n - \delta(c^{n}), ~   \hat{e}_\rho^n=e_{\rho}^{n} - \delta(\rho^{n}), ~~ n\ge 0.
\end{equation*}
It is easy to check that 
\begin{equation}\label{err:initial}
 e_{\rho,i,j}^0= \hat{e}_{c,i,j}^0=0, ~\hat{e}_{\rho,i,j}^0=\mo(h^2), ~ e_{c,i,j}^0=\mo(h^2), ~~ i=1,\ldots, N_x, j=1,\ldots, N_y.
\end{equation}
\begin{theorem}\label{thm:coverg}
Let $\left\{Z^{m}, U^{m}\right\}$ be the solutions to the DeC-MC-BCFD scheme \eqref{2D:PC-BCFD}--\eqref{2D:CN-BCFD:IBc}. Under the  assumptions \eqref{Regularity}--\eqref{bound} and let $\tau \leq K^* h$, then the following estimates hold for $\tau\le \tau_*:=\min\{\tau_1,\tau_2, \tau_3,\tau_4\}, ~ h\le  h_* :=\min\{h_1,h_2,h_3\} $
\begin{equation}\label{thm:converg:result}
	\| c^{m}-Z^{m}\|_{\rm M} 
			+  \| \nabla c^{m} - \bm{d}Z^{m}\|_{\rm TM} + \| \rho^{m}-U^{m}\|_{\rm M} \leq K_0 (\tau^2 + h^2),\quad 1\leq m \leq N,
	\end{equation}
 where $K_0>0$ is a constant, independent of $h$, $\tau$, and $m$, but depends on $K_c$, $K_{\rho}$, $K_\delta $, $K_\epsilon$, $K_{inv}$, the final time $T$ and the bounds of exact solutions, the constants $\{\tau_i\}_{i=1}^4$ and $\{h_i\}_{i=1}^3$ are defined in the proof of Lemma \ref{lem:erho} and Theorem \ref{thm:coverg}, respectively.
\end{theorem}

To facilitate a clear and comprehensive understanding of the proof of Theorem \ref{thm:coverg}, we divide the proof into two separate lemmas. 

\begin{lemma}\label{lem:ec}
 Under the conditions of Theorem \ref{thm:coverg}, there exists a positive constant $K_c$ independent of $h$ and $\tau$ such that
  \begin{equation}\label{lem:ec:sum}
 	\begin{aligned}
		&\|\hat{e}_{c}^{m}\|_{\rm M}^2 + \| \bm{d} \hat{e}_{c}^{m} +  \bm{\epsilon}(c^{m})\|_{\rm TM}^2 \\
       &\quad\leq  K_c\Big( \tau \sum_{n=1}^{m} \| \bm{d} \hat{e}_{c}^{n}  + \bm{\epsilon}(c^n) \|_{\rm TM}^2 
            + \tau\,\|\bar{e}_\rho^1\|_{\rm M}^2 + \tau\sum_{n=1}^{m-1}\|e_{\rho}^{n}\|_{\rm M}^2  +  \tau^4 + h^4\Big), 
	\end{aligned}
\end{equation}
for $1\leq m \leq N$, where henceforth $\bm{\epsilon}(\cdot) := (\epsilon^x(\cdot),\epsilon^y(\cdot))$. 
\end{lemma}
\begin{proof} Note that \eqref{2D:PC-BCFD:c} and \eqref{2D:CN-BCFD:c} have a similar expression for the discretizations of the continuous equation \eqref{model:ks:b} with respect to $c(\bm{x},t)$ at different time levels, except the different treatments of the $\rho$-term. Thus, we can use a unified framework to estimate  $\|\hat{e}_{c}^{n+1}\|_{\rm M}$ and $\|\bm{d}\hat{e}_{c}^{n+1} + \bm{\epsilon}(c^{n+1})\|_{\rm TM}$.

By subtracting \eqref{2D:CN-BCFD:c} from \eqref{exact:CN-BCFD:c}, and adding $\underline{\mF}_{i,j}^{n+1} := -d_t (\delta(c))_{i,j}^{n+1}$ onto both sides of the resulting equation, we obtain the following error equation
 \begin{equation}\label{err:CN-BCFD:c} 
	\begin{aligned}
		d_t \hat{e}_{c,i,j}^{n+1} 
		& = \left[ D_x (d_x \hat{e}_{c} +  \epsilon^x(c) ) \right] _{i,j}^{n+1/2} 
		+ \left[ D_y ( d_y \hat{e}_{c} +  \epsilon^x(c) ) \right] _{i,j}^{n+1/2}-  e_{c,i,j}^{n+1/2}   + \Gamma_{i,j}^{n+1/2}  +  \hat{\mF}_{i,j}^{n+1},
	\end{aligned}
\end{equation}
where $\hat{\mF}^{n+1} := \mF^{n+1} + \underline{\mF}^{n+1}$, 
 and
 \begin{equation}
     \Gamma^{n+1/2}=\begin{cases}
       {\bar{e}_{\rho}^1}/ {2}, & n=0,\\
       e_{\rho}^{*,n+1}, & n\geq 1.
     \end{cases}
 \end{equation}
Then, taking the discrete inner product with $d_t \hat{e}_{c}^{n+1}$, and  using Lemma \ref{lemma:Dd} we obtain
\begin{equation}\label{est:Dec:inn}
\begin{aligned}
	&\|d_t \hat{e}_c^{n+1}\|_{\rm M}^2 \\
	&\quad =- \big( (d_x \hat{e}_{c}  +  \epsilon^x(c) ) ^{n+1/2},   d_x d_t\hat{e}_{c}^{n+1}\big)_x  
 - \big( (d_y \hat{e}_{c}  +  \epsilon^y(c) ) ^{n+1/2},   d_y d_t\hat{e}_{c}^{n+1}\big)_y \\
	&\quad\qquad-(e_c^{n+1/2},  d_t \hat{e}_{c}^{n+1})_{\rm M} 
	 +(\Gamma^{n+1/2},  d_t \hat{e}_{c}^{n+1})_{\rm M}  + (\hat{\mF}^{n+1},  d_t \hat{e}_{c}^{n+1})_{\rm M}
   =: \sum_{i=1}^5 J_i.
\end{aligned}
\end{equation}

Now, we estimate the right-hand side of \eqref{est:Dec:inn} term by term. For the first term, it follows from Cauchy-Schwarz inequality that
 \begin{equation}\label{est:Dec:rhs1}
 	\begin{aligned}
 		J_1 &= - \big( (d_x \hat{e}_{c}  +  \epsilon^x(c) ) ^{n+1/2}, d_t ( d_x \hat{e}_{c}  +  \epsilon^x(c) ) ^{n+1}\big)_x \\
 	 &\qquad + \big( (d_x\hat{e}_{c}  +  \epsilon^x(c) ) ^{n+1/2}, d_t \epsilon^x(c^{n+1})  \big)_x\\
 		&\leq -\f{1}{2\tau}\left( \| (d_x\hat{e}_{c}  +  \epsilon^x(c) )^{n+1}\|_x^2 - \| (d_x\hat{e}_{c} +  \epsilon^x(c) )^{n}\|_x^2 \right)  \\
 		&\qquad + \f{1}{2}\| (d_x \hat{e}_{c}  +  \epsilon^x(c) )^{n+1/2}\|_x^2 
 		+ \f{1}{2}\| d_t \epsilon^x(c^{n+1})\|_x^2,
 	\end{aligned}
 \end{equation}
 and, similar to \eqref{est:Dec:rhs1}, the second term  can be estimated by
 \begin{equation}\label{est:Dec:rhs2}
 	\begin{aligned}
 		J_2
 		&\leq -\f{1}{2\tau}\left( \| (d_y\hat{e}_{c}  +  \epsilon^y(c) )^{n+1}\|_y^2 - \| (d_y\hat{e}_{c}  +  \epsilon^y(c) )^{n}\|_y^2 \right)  \\
 		&\qquad + \f{1}{2}\| (d_y\hat{e}_{c}  +  \epsilon^y(c) )^{n+1/2}\|_y^2 
 		+ \f{1}{2}\| d_t \epsilon^y (c^{n+1})\|_y^2.
 	\end{aligned}
 \end{equation}

For the third term, noting $e_{c}^n=\hat{e}_{c}^n + \delta(c^{n})$ we have 
\begin{equation}\label{est:Dec:rhs3}
	\begin{aligned}
		J_3 &= -\big( (\hat{e}_{c} + \delta(c) )^{n+1/2}, d_t \hat{e}_{c}^{n+1}\big)_{\rm M}\\
		&\leq - \f{1}{2\tau} \big( \|\hat{e}_{c}^{n+1}\|_{\rm M}^2 - \| \hat{e}_{c}^{n}\|_{\rm M}^2 \big) + \|  \delta(c^{n+1/2})\|_{\rm M}^2 +  \f{1}{4} \| d_t \hat{e}_{c}^{n+1}\|_{\rm M}^2.
	\end{aligned}
\end{equation}

Thanks to the Cauchy-Schwarz inequality and Young's inequality, the last two terms can be bounded by
\begin{equation}\label{est:Dec:rhs4} 
	J_4
	\leq \|\Gamma^{n+1/2}\|_{\rm M}^2 + \f{1}{4}\|d_t\hat{e}_{c}^{n+1}\|_{\rm M}^2.
\end{equation}
\begin{equation}\label{est:Dec:rhs5} 
	J_5
	\leq \|\mF^{n+1}\|_{\rm M}^2 + \|\underline{\mF}^{n+1}\|_{\rm M}^2 + \f{1}{2}\|d_t \hat{e}_{c}^{n+1}\|_{\rm M}^2.
\end{equation}

Now, inserting these estimates \eqref{est:Dec:rhs1}--\eqref{est:Dec:rhs5} into \eqref{est:Dec:inn}, and multiplying the resulting equation by $2\tau$, 
we obtain 
 \begin{equation}\label{est:CN-BCFD:c:befsum} 
	\begin{aligned}
		&\| \hat{e}_{c}^{n+1}\|_{\rm M}^2 + \| \bm{d} \hat{e}_{c}^{n+1} +  \bm{\epsilon}(c^{n+1})\|_{\rm TM}^2\\
			&\quad\leq \| \hat{e}_{c}^{n}\|_{\rm M}^2 + \| \bm{d} \hat{e}_{c}^{n} +  \bm{\epsilon}(c^{n})\|_{\rm TM}^2 
                  + \tau \| d_t \bm{\epsilon}(c^{n+1}) \|_{\rm TM}^2   + 2\tau\|\mF^{n+1}\|_{\rm M}^2 + 2\tau\|\underline{\mF}^{n+1}\|_{\rm M}^2   \\
				&\qquad\qquad  + 2\tau\| \delta(c^{n+1/2})\|_{\rm M}^2  + \tau\| \bm{d} \hat{e}_{c}^{n+1/2} +  \bm{\epsilon}(c^{n+1/2})\|_{\rm TM}^2  
           + 2\tau \|\Gamma^{n+1/2}\|_{\rm M}^2,
			\end{aligned}
	\end{equation}
in which, 
 \begin{equation}\label{est:CN-BCFD:c:sum1} 
	\begin{aligned}
		\|\Gamma^{1/2}\|_{\rm M}\leq  \f{1}{2} \|\bar{e}_\rho^1\|_{\rm M},\quad  
        \|\Gamma^{n+1/2}\|_{\rm M}\leq \|e_{\rho}^{*,n+1}\|_{\rm M} \le \f{3}{2}\|e_{\rho}^{n}\|_{\rm M}+ \f{1}{2} \|e_{\rho}^{n-1}\|_{\rm M},~ n\geq 1.
	\end{aligned}
\end{equation}

Therefore, summing \eqref{est:CN-BCFD:c:befsum} over $n$ from $0$ to $m-1$ for $1\leq m \leq N$, and involving the estimates \eqref{est:CN-BCFD:c:sum1}, we have
 \begin{equation*}
	\begin{aligned}
		&\| \hat{e}_{c}^{m}\|_{\rm M}^2 + \| \bm{d} \hat{e}_{c}^{m} +  \bm{\epsilon}(c^{m})\|_{\rm TM}^2\\
		&\quad\leq \| \hat{e}_{c}^{0}\|_{\rm M}^2 + \| \bm{d} \hat{e}_{c}^{0} +  \bm{\epsilon}(c^{0})\|_{\rm TM}^2  + \tau\sum_{n=1}^{m} \| d_t \bm{\epsilon}(c^{n})\|_{\rm TM}^2  + 2\tau\sum_{n=1}^{m}\|\mF^{n}\|_{\rm M}^2 + 2\tau\sum_{n=1}^{m}\|\underline{\mF}^{n}\|_{\rm M}^2  \\
		&\quad\quad + \tau\sum_{n=0}^{m}\| \delta(c^{n})\|_{\rm M}^2   +  \f{\tau}{2} \|\bar{e}_\rho^1\|_{\rm M}^2  
           + 10\tau \sum_{n=0}^{m-1} \|e_\rho^{n}\|_{\rm M}^2
           + \tau \sum_{n=0}^{m} \| \bm{d} \hat{e}_{c}^{n} +  \bm{\epsilon}(c^n) \|_{\rm TM}^2.
	\end{aligned}
\end{equation*}
Thus, the conclusion \eqref{lem:ec:sum} is obtained by gathering the estimates in Lemma \ref{trun:ap-dp}, initial errors \eqref{err:initial}, and truncation errors \eqref{TruncErr:PC:ec1} and \eqref{TruncErr:CN:e1}.
\end{proof}

\begin{remark} Lemma \ref{lem:ec} shows that the errors $\|\hat{e}_{c}^{m}\|_{\rm M}$ and $\|\bm{d}\hat{e}_{c}^{m} + \bm{\epsilon}(c^{m})\|_{\rm TM}$  are related to the errors $\|\bar{e}_{\rho}^{1}\|_{\rm M}$ and $\|e_{\rho}^{n}\|_{\rm M}$ for $1 \leq n \leq m-1 $. These will be estimated in the following lemma.
\end{remark}

\begin{lemma}\label{lem:erho} 
	Under the conditions of Theorem \ref{thm:coverg}, there exists a positive constant $K_\rho$ independent of $h$ and $\tau$ such that
	\begin{equation}\label{lem:erho:sum1}
		\| \bar{e}_{\rho}^{1}\|_{\rm M}^2 \leq  K_{\rho} \big( \tau^3 + h^4 \big),
	\end{equation}
 \begin{equation}\label{lem:erho:sum} 
		\begin{aligned}
			 \| \hat{e}_{\rho}^{m}\|_{\rm M}^2 
			&\leq  K_{\rho} \Big(\tau\sum_{n=0}^{m}\|\hat{e}_{\rho}^{n}\|_{\rm M}^2  
			 	+   \tau\sum_{n=0}^{m} \| \bm{d} \hat{e}_{c}^{n} +  \bm{\epsilon}(c^{n})\|_{\rm TM}^2 \\
              &\qquad \qquad + \tau \max_{0 \le n \le m} \|\bm{d} Z^{n}\|_\infty^2 \, \sum_{n=0}^{m}  \| e_\rho^{n} \|_{\rm M}^2+  \tau^4 + h^4 \Big), ~1\leq m \leq N. 
		\end{aligned}
	\end{equation}

\end{lemma}
\begin{proof}
The proof is split into two steps, due to the two different discretizations of the continuous equation \eqref{model:ks:a} about $\rho(\bm{x},t)$ at different time levels, see equations \eqref{2D:PC-BCFD:rho}, \eqref{2D:PC-BCFD:rho:e2} and \eqref{2D:CN-BCFD:rho}.

\paragraph{\bf Step I. Estimate for $\|\bar{e}_\rho^1\|_{\rm M}$} 
By subtracting \eqref{2D:PC-BCFD:rho} from \eqref{exact:PC-BCFD:rho1} and noting the zero initial error \eqref{err:initial} for $e_\rho^0$, we obtain the following error equation
\begin{equation}\label{err:PC-BCFD:rho1} 
	\begin{aligned}
    \f{\hat{\bar{e}}_{\rho,i,j}^{1}+ \delta_{i,j}(\rho^{1}) }{\tau} 
    &= \big[ D_x ( d_x \hat{\bar{e}}_\rho +  \epsilon^x(\rho) ) \big] _{i,j}^{1} 
		+ \big[ D_y ( d_y \hat{\bar{e}}_\rho +  \epsilon^y(\rho) ) \big] _{i,j}^{1}  -  \lambda \Lambda_{i,j}^1 +  \mR^1_{i,j}, 
	\end{aligned}
\end{equation}
where  
\begin{equation*}\label{err:PC-BCFD:lbda} 
	\begin{aligned}
		\Lambda^1
		&:=    D_x \Big( [\ell_h\rho]^1 [ d_x Z^0 +  \epsilon^x(c^o ) ]
		    - [\ell_h \bar{U}]^1 [d_x Z]^0 \Big)   \\
		&\qquad+  D_y \Big( [\ell_h\rho]^1 [ d_y Z^0 +  \epsilon^y(c^o ) ]
              - [\ell_h \bar{U}]^1 [d_y Z]^0 \Big).
	\end{aligned}
\end{equation*}
Then, taking discrete inner product with $\tau \hat{\bar{e}}_\rho^1$, using Lemma \ref{lemma:Dd} we have
 \begin{equation}\label{inn:PC-BCFD:rho1} 
	\begin{aligned}
		\|\hat{\bar{e}}_\rho^{1}\|_{\rm M}^2  
       &= -\tau\left(  d_x \hat{\bar{e}}_\rho^{1} +  \epsilon^x(\rho^{1}), d_x \hat{\bar{e}}_\rho^{1} \right)_x
         -\tau\left(  d_y \hat{\bar{e}}_\rho^{1} +  \epsilon^y(\rho^{1}), d_y \hat{\bar{e}}_\rho^{1} \right)_y \\
		&\qquad -  \lambda \tau \left( \Lambda^1, \hat{\bar{e}}_\rho^{1} \right)_{\rm M}
		 +  \tau ( \mR^1,\hat{\bar{e}}_\rho^{1})_{\rm M} -(\delta(\rho^{1}),\hat{\bar{e}}_\rho^{1})_{\rm M}=:\sum_{i=1}^{5} I_i. 
	\end{aligned}
\end{equation}

Now, we estimate the right-hand side of \eqref{inn:PC-BCFD:rho1} term by term. Firstly, the first two terms can be analyzed via Cauchy-Schwarz inequality that
 \begin{equation}\label{est:pc-BCFD:rho1:rhs12} 
	\begin{aligned}
		I_1+I_2 	& \leq -\f{\tau}{2}\|  d_x \hat{\bar{e}}_\rho^{1} \|_x^2 -\f{\tau}{2}\| d_y \hat{\bar{e}}_\rho^{1}\|_y^2
  + \f{\tau}{2}\|\epsilon^x(\rho^{1})\|_x^2	+ \f{\tau}{2}\|\epsilon^y(\rho^{1})\|_y^2\\
          &= -\f{\tau}{2}  \|\bm{d}\hat{\bar{e}}_\rho^1\|_{\rm TM}^2 + \f{\tau}{2}\|\bm{\epsilon}(\rho^{1})\|_{\rm TM}^2.
	\end{aligned}
\end{equation}

Secondly, we pay attention to the third term of the right-hand side of \eqref{inn:PC-BCFD:rho1}. Following Lemma \ref{lemma:Dd}, we  derive an equivalent form
\begin{equation}\label{est:pc-erho:rhs3}
	\begin{aligned} 
I_3 & = -  \lambda\tau \left( D_x ( [\ell_h\rho]^1 [ d_x Z^0 +  \epsilon^x(c^o) ] 
 - [\ell_h \bar{U}]^1 [d_x Z]^0),\hat{\bar{e}}_\rho^1 \right)_{\rm M}\\
&\qquad -  \lambda\tau \left( D_y ( [\ell_h\rho]^1 [ d_y Z^0 +  \epsilon^y(c^o) ] - [\ell_h\bar{U}]^1 [d_y Z]^0),\hat{\bar{e}}_\rho^1 \right)_{\rm M} \\
& = \lambda\tau \left(   [\ell_h\rho]^1 [ d_x Z^0 +  \epsilon^x(c^o) ] 
   - [\ell_h\bar{U}]^1 [d_x Z]^0, d_x\hat{\bar{e}}_\rho^1 \right)_x\\
&\qquad +  \lambda\tau \left(   [\ell_h\rho]^1 [ d_y Z^0 +  \epsilon^y(c^o) ] 
   - [\ell_h\bar{U}]^1 [d_y Z]^0 ,d_y\hat{\bar{e}}_\rho^1 \right)_y
=: I_{31} + I_{32},
 \end{aligned}
\end{equation}
in which, the first term of \eqref{est:pc-erho:rhs3} can be bounded by
\begin{equation}\label{est:pc-erho:rhs3:1} 
	\begin{aligned}
     I_{31}
	&\leq \lambda\tau \| [\ell_h\rho]^1 [ d_x Z^0 +  \epsilon^x(c^o)] - [\ell_h\bar{U}]^1 [d_x Z]^0 \|_x  \|d_x\hat{\bar{e}}_\rho^1\|_x\\
	&\leq \lambda^2 \tau \| [\ell_h \bar{e}_\rho]^1 [ d_x Z]^0 +  [\ell_h \rho]^1 [\epsilon^x(c^o)] \|_x^2  +\f{\tau}{4}\|d_x\hat{\bar{e}}_\rho^1\|_x^2\\
   &\leq K_1 \tau\left(  \|  \epsilon^x(c^o) \|_x^2	+ \| \bar{e}_\rho^1 \|_{\rm M}^2\right) 
		+\f{\tau}{4}\|d_x\hat{\bar{e}}_\rho^1\|_x^2,
		\end{aligned}
\end{equation}
where $K_1 := 2\lambda^2\max\left\{ \|\rho^1\|_\infty^2, \|\bm{d} Z^0\|_\infty^2\right\}$.
The second term of \eqref{est:pc-erho:rhs3} yields a similar result to \eqref{est:pc-erho:rhs3:1} that
\begin{equation}\label{est:pc-erho:rhs3:2} 
   \begin{aligned}
	I_{32} \leq K_1\tau \left(  \|   \epsilon^y(c^o) \|_y^2	+ \| \bar{e}_\rho^1 \|_{\rm M}^2\right) 
		+\f{\tau}{4}\|d_y\hat{\bar{e}}_\rho^1\|_y^2,
   \end{aligned}
\end{equation}
and therefore, following  \eqref{est:pc-erho:rhs3:1}--\eqref{est:pc-erho:rhs3:2} and the fact that 
$
        \| \bar{e}_{\rho}^1\|_{\rm M} \le \| \hat{\bar{e}}_\rho^{1}\|_{\rm M}+ \| \delta(\rho^{1})\|_{\rm M}
$, the third term of the right-hand side of \eqref{inn:PC-BCFD:rho1} can be estimated as
\begin{equation}\label{est:pc-BCFD:rho1:rhs3}
	\begin{aligned} 
		I_3 & \leq  K_1\tau\, \|  \bm{\epsilon}(c^o)\|_{\rm TM}^2 
                    +\f{\tau}{4}  \|\bm{d}\hat{\bar{e}}_\rho^1\|_{\rm TM}^2 
             + 4K_1\tau\,\| \delta(\rho^{1})\|_{\rm M}^2		
             + 4K_1\tau\| \hat{\bar{e}}_\rho^1 \|_{\rm M}^2.
	\end{aligned}
\end{equation}

Next, for the fourth term on the right-hand side of \eqref{inn:PC-BCFD:rho1}, the definition of $\mR^1$ in \eqref{exact:PC-BCFD:rho1:e1} implies that
\begin{equation}\label{est:pc-BCFD:rho1:rhs4bf}
    I_4 = \tau ( \mR_{I}^1 ,\hat{\bar{e}}_\rho^{1})_{\rm M} 
        + \tau ( \mR_{II}^1,\hat{\bar{e}}_\rho^{1})_{\rm M} =: I_{41} + I_{42}, 
\end{equation}
in which the first term can be estimated by applying the Cauchy-Schwarz inequality as
\begin{equation}\label{est:pc-BCFD:rho1:rhs4:1}
    I_{41} \leq 2 \tau^2 \|\mR_{I}^1\|_{\rm M}^2 + \f{1}{8}\|\hat{\bar{e}}_\rho^1\|_{\rm M}^2.
\end{equation}
Moreover, it follows from Lemmas \ref{lemma:Dd} and \ref{lem:interp} that $I_{42}$ can be estimated as
\begin{equation}\label{est:pc-BCFD:rho1:rhs4:2}
\begin{aligned}
    I_{42} &= - \lambda\tau \big( [\ell_h\rho-\rho]^1 \p_x c^o, d_x \hat{\bar{e}}_\rho^{1} )_x 
             - \lambda\tau \big( [\ell_h\rho-\rho]^1 \p_y c^o, d_y \hat{\bar{e}}_\rho^{1} )_y\\
          &\leq \lambda^2\tau K \|\rho\|_{L^{\infty}([t^{0},t^{1}] ;W^{2,\infty}(\Omega))}^2\|c\|_{L^{\infty}([t^{0},t^{1}] ;W^{1,\infty} (\Omega))}^2 h^4  +\f{\tau}{4}  \|\bm{d}\hat{\bar{e}}_\rho^1\|_{\rm TM}^2,
\end{aligned}
\end{equation}
and thus, the fourth term on the right-hand side of \eqref{inn:PC-BCFD:rho1} can be estimated as
\begin{equation}\label{est:pc-BCFD:rho1:rhs4}
\begin{aligned}
    I_4 &\leq 2 \tau^2 \|\mR_{I}^1\|_{\rm M}^2 + \lambda^2\tau K \|\rho^1\|_{L^{\infty}([t^{0},t^{1}] ;W^{2,\infty}(\Omega) )}^2 \|c\|_{L^{\infty}([t^{0},t^{1}] ;W^{1,\infty} (\Omega) )}^2 h^4 \\
   &\qquad + \f{1}{8} \| \hat{\bar{e}}_\rho^1 \|_{\rm M}^2 + \f{\tau}{4}  \|\bm{d}\hat{\bar{e}}_\rho^1\|_{\rm TM}^2.
\end{aligned}
\end{equation}

Finally, thanks to the Cauchy-Schwarz inequality and Young's inequality, the last term on the right-hand side of \eqref{inn:PC-BCFD:rho1} can be estimated as 
\begin{equation}\label{est:pc-BCFD:rho1:rhs5}
		I_5\leq  2 \|\delta(\rho^{1}) \|_{\rm M}^2+ \f{1}{8}\|\hat{\bar{e}}_\rho^1\|_{\rm M}^2.
\end{equation}

Inserting these estimates \eqref{est:pc-BCFD:rho1:rhs12}, \eqref{est:pc-BCFD:rho1:rhs3} and \eqref{est:pc-BCFD:rho1:rhs4}--\eqref{est:pc-BCFD:rho1:rhs5}
into \eqref{inn:PC-BCFD:rho1}, for a sufficiently small chosen $\tau_1$, for example, $\tau_1 :=1/(16K_1)$, we obtain from Lemma \ref{trun:ap-dp} and estimate \eqref{TruncErr:PC:e11} that 
\begin{equation}
    \begin{aligned}
\| \hat{\bar{e}}_\rho^{1}\|_{\rm M}^2 
        &\leq K\left( \tau \| \bm{\epsilon}(c^o)  \|_{\rm TM}^2 + \tau \|\bm \epsilon(\rho^{1}) \|_{\rm TM}^2
            +  \| \delta(\rho^{1})\|_{\rm M}^2 + \tau^2 \|\mR_{I}^1\|_{\rm M}^2 
             + \tau h^4 \right)\\  
            & \leq  K\left( \tau^3 + h^4 \right), \quad \text{for}~ \tau\le \tau_1.
    \end{aligned}
\end{equation}
which further implies
$
      \| \bar{e}_{\rho}^1\|_{\rm M}^2  \le 2\| \hat{\bar{e}}_\rho^{1}\|_{\rm M}^2 + 2\| \delta(\rho^{1})\|_{\rm M}^2  
        \leq  K_{\rho}\big( \tau^3 + h^4 \big)
$ by the triangle inequality.

\paragraph{\bf Step II. Estimate for $\|\hat{e}_\rho^{n+1}\|_{\rm M}$ $(0 \leq n \leq m-1 )$} 
By subtracting \eqref{exact:CN-BCFD:rho} from \eqref{2D:CN-BCFD:rho}, and adding $\underline{\mS}^{n+1}_{i,j} := -d_t (\delta(\rho))_{i,j}^{n+1}$ onto both sides of the resulting equation, we obtain the following error equation
\begin{equation}\label{err:CN-BCFD:rho} 
	\begin{aligned}
		d_t \hat{e}_{\rho,i,j}^{n+1} 
		& = \big[ D_x ( d_x \hat{e}_{\rho} +  \epsilon^x(\rho) ) \big] _{i,j}^{n+1/2} 
		+ \big[ D_y ( d_y \hat{e}_{\rho} +  \epsilon^y(\rho) ) \big] _{i,j}^{n+1/2}
  -  \lambda \Phi_{i,j}^{n+1/2} +  \hat{\mS}_{i,j}^{n+1}, 
	\end{aligned}
\end{equation}
where $\hat{\mS}^{n+1}:=\mS^{n+1} + \underline{\mS}^{n+1}$,
 and for $0\leq l \leq N$,
 \begin{equation*}
 	\begin{aligned}
 		\Phi_{i,j}^{l} 
 		&:=  \big[D_x ( [\ell_h\rho] [ d_x (c-\delta(c)) +  \epsilon^x(c) ] 
 		- [\ell_hU] [d_x Z]) \big]_{i,j}^{l} \\
 		&\quad+ \big[D_y ( [\ell_h\rho] [ d_y (c-\delta(c)) +  \epsilon^y(c) ] - [\ell_hU] [d_y Z]) \big]_{i,j}^{l}.
 	\end{aligned}
 \end{equation*}

Then, taking discrete inner product with $\hat{e}_{\rho}^{n+1/2}$, and using Lemma \ref{lemma:Dd} we have
\begin{equation}\label{inn:CN-BCFD:rho} 
   \begin{aligned}
	&\f{1}{2\tau} \big( \| \hat{e}_{\rho}^{n+1}\|_{\rm M}^2 - \| \hat{e}_{\rho}^{n}\|_{\rm M}^2\big)
            = \big( d_t \hat{e}_{\rho}^{n+1}, \hat{e}_{\rho}^{n+1/2}\big)_{\rm M} \\
               & \quad  = -\big( ( d_x \hat{e}_{\rho} +  \epsilon^x(\rho))^{n+1/2}, d_x\hat{e}_{\rho}^{n+1/2} \big)_{\rm M}  
                          - \big( ( d_y \hat{e}_{\rho} +  \epsilon^x(\rho) )^{n+1/2}, d_y\hat{e}_{\rho}^{n+1/2} \big)_{\rm M} \\
				&\qquad -  \lambda \big( \Phi^{n+1/2},\hat{e}_{\rho}^{n+1/2} \big)_{\rm M}
				         +  \big( \hat{\mS}^{n+1},\hat{e}_{\rho}^{n+1/2}\big)_{\rm M}. 
			\end{aligned}
	\end{equation}

Note that the estimates for the right-hand side terms in \eqref{inn:CN-BCFD:rho} are very similar to those of Lemma \ref{lem:ec},
we only pay special attention to the third term, denoted by $H_3$. By Lemma \ref{lemma:Dd}, and similar to \eqref{est:pc-erho:rhs3}, we obtain
\begin{equation}\label{est:erho:rhs3}
	\begin{aligned} 
   H_3 &=  \lambda \big( ( [\ell_h\rho]  \left[ d_x (c-\delta(c)) +  \epsilon^x(c) \right]  - [\ell_hU] [d_x Z])^{n+1/2},d_x \hat{e}_{\rho}^{n+1/2} \big)_{\rm M}\\
 &\qquad +  \lambda \big(  ( [\ell_h\rho] \left[ d_y (c-\delta(c)) +  \epsilon^y(c) \right] - [\ell_hU] [d_y Z]) ^{n+1/2}, d_y \hat{e}_{\rho}^{n+1/2} \big)_{\rm M} \\
 &=  \lambda \big( ( [\ell_h \rho]  [ d_x \hat{e}_c +  \epsilon^x(c) ]  + [\ell_h e_\rho] [d_x Z])^{n+1/2},d_x \hat{e}_{\rho}^{n+1/2} \big)_{\rm M}\\
 &\qquad +  \lambda \big(  ( [\ell_h\rho] [ d_y \hat{e}_c +  \epsilon^y(c) ] + [\ell_h e_\rho] [d_y Z]) ^{n+1/2}, d_y \hat{e}_{\rho}^{n+1/2} \big)_{\rm M}, 
 \end{aligned}
\end{equation}
which can be estimated similarly to \eqref{est:pc-BCFD:rho1:rhs3} as
\begin{equation}\label{est:erho:rhs3f}
	H_3  \leq  K  \sum_{q=0}^{1}  \left( \|  \bm{d} \hat{e}_c^{n+q} +  \bm{\epsilon}(c^{n+q}) \|_{\rm TM}^2
		+  \|\bm{d} Z^{n+q}\|_\infty^2 \|e_\rho^{n+q} \|_{\rm M}^2\right) 
		 + \f{1}{4}\|\bm{d}\hat{e}_{\rho}^{n+1/2}\|_{\rm TM}^2.
\end{equation}

Now, inserting \eqref{est:erho:rhs3f} into \eqref{inn:CN-BCFD:rho}, similar to  the estimate \eqref{est:CN-BCFD:c:befsum} in Lemma \ref{lem:ec}, we can obtain
\begin{equation}\label{est:CN-BCFD:rho:befsum}
\begin{aligned}
       \| \hat{e}_{\rho}^{n+1}\|_{\rm M}^2 
    &\leq \| \hat{e}_{\rho}^{n}\|_{\rm M}^2 + \tau\|\bm{\epsilon}(\rho^{n+1/2})\|_{\rm TM}^2 + 2\tau \|\mS_{I}^{n+1}\|_{\rm M}^2+ K\tau h^4 + 2\tau \|\underline{\mS}^{n+1}\|_{\rm M}^2 \\
    &\quad + K \tau \sum_{q=0}^{1}  \left( \|  \bm{d} \hat{e}_c^{n+q} +  \bm{\epsilon}(c^{n+q}) \|_{\rm TM}^2
		+  \|\bm{d} Z^{n+q}\|_\infty^2 \|e_\rho^{n+q} \|_{\rm M}^2\right) + \tau\|\hat{e}_{\rho}^{n+1/2}\|_{\rm M}^2.
\end{aligned}
\end{equation}
 Therefore, summing \eqref{est:CN-BCFD:rho:befsum} over $n$ from $0$ to $m-1$ for $1\leq m \leq N$, we have
\begin{equation}\label{result:CN-rho}
	\begin{aligned}
		\| \hat{e}_{\rho}^{m}\|_{\rm M}^2  
		&\leq \| \hat{e}_{\rho}^{0}\|_{\rm M}^2 + K\tau\,\sum_{n=1}^{m}\|\bm{\epsilon}(\rho^{n})\|_{\rm TM}^2 + K\tau\,\sum_{n=1}^{m}\|\mS_{I}^{n}\|_{\rm M}^2 + K h^4 + K\tau\,\sum_{n=1}^{m}\|\underline{\mS}^{n}\|_{\rm M}^2 \\
		&\qquad +   K\tau\,\sum_{n=0}^{m}\left( \|\hat{e}_{\rho}^{n}\|_{\rm M}^2+ \|  \bm{d} \hat{e}_c^{n} +  \bm{\epsilon}(c^{n}) \|_{\rm TM}^2 +  \|\bm{d} Z^{n}\|_\infty^2 \|e_\rho^{n} \|_{\rm M}^2  \right), 
	\end{aligned}
\end{equation}
which implies the conclusion \eqref{lem:erho:sum} by combining the estimates in Lemma \ref{trun:ap-dp}, \eqref{TruncErr:PC:e2}, \eqref{TruncErr:CN:e2} and \eqref{err:initial}.
\end{proof}

\begin{remark}  We remark that in Lemma \ref{lem:erho}, the term $\| \bm{d} \hat{e}_{c}^{n} +  \bm{\epsilon}(c^{n})\|_{\rm TM}$ is also related to the conclusion of Lemma \ref{lem:ec}. Thus, a combination of Lemmas \ref{lem:ec} and \ref{lem:erho} may deduce the final conclusion of Theorem \ref{thm:coverg}. While we still need to prove that the term $\|\bm{d} Z^{n}\|_\infty$ in \eqref{lem:erho:sum}  is uniformly bounded for all $n$ using the mathematical induction method.
\end{remark}

We are now in position to prove our main results.
\paragraph{\bf Proof of Theorem \ref{thm:coverg}} 
To complete the proof, we use the mathematical induction method to prove that $\|\bm{d} Z^{n}\|_\infty$ is bounded, and thus the conclusion \eqref{thm:converg:result} holds. Noting the assumption \eqref{bound} indicates that 
\begin{equation}\label{inf:Z0}
    \|\bm{d} Z^0\|_\infty = \| \bm{d} c^0 - \bm{d}\delta(c^0) \|_{\infty} \leq K_* + K_\delta h  \le K_* + 1, \quad \text{for}~ h \le h_1:= 1/K_\delta.
\end{equation}

First, at the first time level, i.e., $m=1$, thanks to \eqref{lem:ec:sum} in Lemma \ref{lem:ec} and \eqref{lem:erho:sum1} in Lemma \ref{lem:erho}, we obtain, for a sufficiently small chosen $\tau_2$, for example, $\tau_2 :=1/(2K_c)$,
  \begin{equation}\label{reslt:thm:t1c}
		\| \hat{e}_{c}^{1}\|_{\rm M}^2 
		+ \| \bm{d} \hat{e}_{c}^{1} +  \bm{\epsilon}(c^{1})\|_{\rm TM}^2  \leq (2 K_c+K_\rho) ( \tau^4 + h^4),  \quad \text{for}~ \tau\le \tau_2,
\end{equation}
which implies that
\begin{equation*}
		 \| \bm{d} \hat{e}_{c}^{1} +  \bm{\epsilon}(c^{1})\|_{\rm TM}  
     \leq \sqrt{2K_c+K_\rho} ( \tau^2 + h^2).
\end{equation*}
 Then, applying the triangle inequality, the inverse estimate with constant $K_{inv}$, Lemma \ref{trun:ap-dp}, and assumption \eqref{bound}, we have
\begin{equation}\label{inf:dZ}
	\begin{aligned}
	\|\bm{d} Z^1\|_\infty &\leq  \| \bm{d} \hat{e}_c^{1}  +  \bm{\epsilon}(c^1) \|_\infty
	+\|\bm{d} \delta(c^1) \|_\infty + \| \bm{\epsilon}(c^1) \|_{\infty} + \| \bm{d} c^1 \|_{\infty}\\
 	&\leq K_{inv} h^{-1}\| \bm{d} \hat{e}_{c}^{1} +  \bm{\epsilon}(c^{1})\|_{\rm TM} + K_{\delta} h +K_\epsilon h^2 + \| \nabla c^1 \|_{\infty}\\
 &\leq  (K_{inv}\sqrt{2K_c+K_\rho}+K_\delta + K_\epsilon)  (h+ h^{-1}\tau^{2}) + \| \nabla c^1 \|_{\infty}.
 	\end{aligned}
\end{equation}

Let $\tau \leq K^* h$ for some positive constant $K^*$, and let $h_2$ be a  small enough positive constant such that
$
(K_{inv}\sqrt{2K_c+K_\rho}+K_\delta + K_\epsilon)  \left(1 + (K^*)^2\right) h_2 \leq 1.
$
Then, for $h \in (0, h_2]$, we derive from \eqref{inf:dZ} and assumption \eqref{bound} that
\begin{equation}\label{inf:Z1}
			\|\bm{d} Z^1\|_\infty  \leq K_* + 1.
\end{equation}
Therefore, inserting \eqref{inf:Z0} and \eqref{inf:Z1} into \eqref{lem:erho:sum} for $m=1$ in Lemma \ref{lem:erho}, and utilizing the conclusion \eqref{reslt:thm:t1c}, 
we have
\begin{equation}\label{est:erho:1}
   \| \hat{e}_{\rho}^{1}\|_{\rm M}^2 \leq K_2 (\tau^4 + h^4), \quad \text{for}~ \tau\le \tau_3:=1/(K_2),
\end{equation}
where $K_2 := 2\max\{ K_\rho(2K_c + K_\rho + K_\epsilon^2 + K_\delta^2 + 1) T, K_\rho (2(K_*+1)^2+1)T\}$. Thus, combinations of \eqref{est:erho:1} and \eqref{reslt:thm:t1c} imply that
  \begin{equation}\label{reslt:thm:t1}
		\| \hat{e}_{c}^{1}\|_{\rm M}^2 
		+ \| \bm{d} \hat{e}_{c}^{1} +  \bm{\epsilon}(c^{1})\|_{\rm TM}^2 + \| \hat{e}_{\rho}^{1}\|_{\rm M}^2 \leq K_3 ( \tau^4 + h^4),  \quad \text{for}~ \tau\le \min \{\tau_2, \tau_3\},
\end{equation}
where $K_3 :=  2K_c + K_\rho + K_2$.

Now, suppose that $\|\bm{d} Z^{l}\|_{\infty} \leq K_* + 1$ holds true for all $l\le m$ with some $m \ge 2$. Then, adding the estimates \eqref{lem:ec:sum} in Lemma \ref{lem:ec} and \eqref{lem:erho:sum} in Lemma \ref{lem:erho} together, and inserting the assumption into the resulting equation,
we have 
  \begin{equation} \label{est:ec:l-0}
	\begin{aligned}
		&\| \hat{e}_c^{l}\|_{\rm M}^2 	
             +\| \bm{d} \hat{e}_{c}^{l} +  \bm{\epsilon}(c^{l})\|_{\rm TM}^2 + \| \hat{e}_{\rho}^{l}\|_{\rm M}^2\\
		&\quad\leq K_4 \Big(\tau\sum_{n=0}^{l} \| \bm{d} \hat{e}_{c}^{n} +  \bm{\epsilon}(c^{n})\|_{\rm TM}^2 
		+   \tau\sum_{n=0}^{l}\|\hat{e}_\rho^{n}\|_{\rm M}^2 +   ( \tau^4 + h^4)\Big),
	\end{aligned}
\end{equation}
where $\|e_\rho^{n}\|_{\rm M} \leq \|\hat{e}_\rho^{n}\|_{\rm M}+ \|\delta(\rho^{n})\|_{\rm M}$ is applied, and the constant 
$K_4 :=\max\{2K_c + K_\rho + 2K_\rho (K_*+1)^2), 2K_c + K_\rho +  2( K_\rho (K_*+1)^2 + K_c) TK_\delta^2 \}$.
An application of the discrete Gr\"{o}nwall's inequality to \eqref{est:ec:l-0} directly yields
  \begin{equation} \label{est:ec:l-1}
	\begin{aligned}
		&\| \hat{e}_\rho^{l}\|_{\rm M}^2 \leq   K_5 \big( \tau^4 + h^4 \big),  \quad \text{for}~ \tau\le \tau_4 :=1/(2K_4),
	\end{aligned}
\end{equation}
 where $K_5 := 2 K_4 {\rm e}^{2K_4T}$, and the truth that $K_4\tau \le 1/2$ for $\tau\le \tau_4$ is applied.

Next, we prove that  $\|\bm{d} Z^{m+1}\|_{\infty} \leq K_*+1$ holds true. 
Bring the result \eqref{est:ec:l-1} into \eqref{lem:ec:sum}, and then, applying the discrete Gr\"{o}nwall's inequality we get
  \begin{equation}\label{est:ec:l}
		\| \hat{e}_c^{m+1}\|_{\rm M}^2 
		+\| \bm{d} \hat{e}_{c}^{m+1} +  \bm{\epsilon}(c^{m+1})\|_{\rm TM}^2\leq  K_6 ( \tau^4 + h^4 ),  \quad \text{for}~ \tau\le \tau_2.
\end{equation}
where the constant $K_6 := 2 K_c (2 T K_5 + 2 T K_\delta^2 + 2){\rm e}^{2K_c T}$. 
Thus, following the same procedure of \eqref{inf:dZ}, we have
\begin{equation}\label{inf:dZm+1}
		\|\bm{d} Z^{m+1}\|_\infty 
    \leq (K_{inv}\sqrt{K_6}+K_\delta + K_\epsilon)(h + h^{-1}\tau^2) + \| \nabla c^{m+1} \|_{\infty},
\end{equation}
and then, for $\tau \leq K^* h$, there exists a small enough positive constant $h_3$ such that
$
(K_{inv}\sqrt{K_6}+K_\delta + K_\epsilon)\left(1 + (K^*)^2 \right) h_3 \leq 1.
$
Therefore, for $h \in (0, h_3]$, we derive from \eqref{inf:dZm+1} and assumption \eqref{bound} that
\begin{equation}\label{inf:reslt:dZm+1}
	\|\bm{d} Z^{m+1}\|_\infty  \leq K_* + 1,
\end{equation}
for $\tau \le \tau_*:= \min\{\tau_1,\tau_2, \tau_3, \tau_4\}$ and $h \le h_*:= \min\{h_1, h_2, h_3\}$. This completes the induction.

Finally, we insert \eqref{inf:reslt:dZm+1} into \eqref{lem:erho:sum}, with the same procedure of \eqref{est:ec:l-0}--\eqref{est:ec:l-1}, and applying the discrete Gr\"{o}nwall's inequality, we obtain  
\begin{equation}\label{est:erho:l}
	\| \hat{e}_{\rho}^{m+1} \|_{\rm M}^2 \leq K_5 ( \tau^4 + h^4 ),  \quad \text{for}~ \tau\le \tau_*.
\end{equation}
Therefore, it follows from \eqref{est:ec:l} and \eqref{est:erho:l} that
  \begin{equation} \label{reslt:thm:tm+1}
	\begin{aligned}
		&\| \hat{e}_c^{m+1}\|_{\rm M}^2 	+\| \bm{d} \hat{e}_{c}^{m+1} +  \bm{\epsilon}(c^{m+1})\|_{\rm TM}^2  + \| \hat{e}_\rho^{m+1}\|_{\rm M}^2 \leq  K_7 ( \tau^4 + h^4 ),  \quad \text{for}~ \tau\le \tau_*,
	\end{aligned}
\end{equation}
where $K_7 := K_5 + K_6$.
 
To gain the final result \eqref{thm:converg:result}, we note that $\|e_\rho^{n}\|_{\rm M} \leq \|\hat{e}_\rho^{n}\|_{\rm M}+ \|\delta(\rho^{n})\|_{\rm M}$, $\|e_c^{n}\|_{\rm M} \leq \|\hat{e}_c^{n}\|_{\rm M}+ \|\delta(c^{n})\|_{\rm M}$ and $\|\bm{d} \hat{e}_c^{n} +  \bm{\epsilon}(c^{n})\|_{\rm TM} = \|\nabla c^{n} - \bm{d}Z^{n}\|_{\rm TM}$. Thus, the estimate \eqref{reslt:thm:tm+1} combined with Lemma \ref{trun:ap-dp} proves Theorem \ref{thm:coverg}. 
\qed

It remains to prove the existence and uniqueness of the solutions $\{Z^{n}, U^{n}\}_{n=1}^{N}$ to the DeC-MC-BCFD scheme \eqref{2D:PC-BCFD}--\eqref{2D:CN-BCFD:IBc}. Although the scheme is linear, the corresponding coefficient matrix of the resulting linear algebraic system of the variable $U^{n}$ is directly related to $Z^{n}$. Therefore, the proof shall be based on the bounded inductive assumption $\|\bm{d} Z^{n}\|_{\infty} \leq K_* + 1$ in the convergence analysis, and should be carried out via the time-marching method. Given the complexity and length of the preceding proof, this proof is presented separately in the following theorem to ensure clarity and thoroughness.

\begin{theorem}\label{thm:existence}  Assume the conditions in Theorem \ref{thm:coverg} hold. If $\tau < \f{4}{\lambda^2(K_*+1)^2}$, there exist unique solutions $\{Z^{n}, U^{n}\}_{n=1}^{N}$ to the DeC-MC-BCFD scheme \eqref{2D:PC-BCFD}--\eqref{2D:CN-BCFD:IBc}.
\end{theorem}
\begin{proof} Noting that \eqref{2D:PC-BCFD}--\eqref{2D:CN-BCFD:IBc} form linear square systems in a finite-dimensional space for each pair $\{Z^{n}, U^{n}\}$, the uniqueness of the solutions also imply the existence. Therefore, we only pay attention to the proof of uniqueness of $\{Z^{n}, U^{n}\}$ that satisfying \eqref{2D:PC-BCFD}--\eqref{2D:CN-BCFD:IBc}. 
To this aim, let $\{W^{n}, V^{n}\}_{n=0}^{N}$ and $\{Z^{n}, U^{n}\}_{n=0}^{N}$ be two solution pairs of \eqref{2D:PC-BCFD}--\eqref{2D:CN-BCFD:IBc} with $W^0=Z^0$ and $V^0=U^0$, and also assume the first time level prediction solutions corresponding to \eqref{2D:PC-BCFD:rho} are denoted by $ \bar{V}^1$ and $\bar{U}^1$ respectively. 

Let $\theta^n=Z^n-W^n$, $\mu^n=U^n-V^n$, and $\bar{\mu}^1=\bar{U}^1-\bar{V}^1$. It is clear that $\theta^0=\mu^0=0$. We shall prove by the induction method that $\theta^n=\mu^n\equiv 0$ for all $n=1, \ldots, N$. 

\paragraph{\bf Step I: Uniqueness of the approximate solutions $( Z^1, U^1)$}
We observe from the prediction-correction BCFD scheme \eqref{2D:PC-BCFD} that $(\bar{\mu}^1, \theta^1, \mu^1)$ are solutions of
 \begin{subequations}\label{2D:PC-BCFD:Unique}
	\begin{numcases}{}
		\f{1}{\tau} \bar{\mu}_{i,j}^{1} = [D_x(d_x \bar{\mu})]_{i,j}^{1} + [D_y(d_y \bar{\mu})]_{i,j}^{1} \nonumber\\
		\qquad\qquad\qquad -\lambda \left( \left[D_x ([\ell_h \bar{\mu}]^1 [d_x Z]^0)\right]_{i,j} + \left[D_y ([\ell_h \bar{\mu}]^1 [d_y Z]^0)\right]_{i,j}\right) ,  \label{2D:PC-BCFD:Unique:rho}	\\
  \f{1}{\tau} \theta_{i,j}^{1} = \f{1}{2}[D_x(d_x \theta)]_{i,j}^{1} + \f{1}{2}[D_y(d_y \theta)]_{i,j}^{1}
 	-  \f{1}{2}\theta_{i,j}^{1} + \f{1}{2}\bar{\mu}_{i,j}^{1}, \label{2D:PC-BCFD:Unique:c}\\
  		\f{1}{\tau} \mu_{i,j}^{1} = \f{1}{2}[D_x(d_x \mu)]_{i,j}^{1} + \f{1}{2}[D_y(d_y \mu)]_{i,j}^{1} \nonumber\\
		\qquad\qquad\qquad -\f{\lambda}{2} \left( \left[D_x ([\ell_h U] [d_x Z])\right]_{i,j}^{1} + \left[D_y ([\ell_h U] [d_y Z])\right]_{i,j}^{1}\right) \nonumber\\
        \qquad\qquad\qquad + \f{\lambda}{2} \left( \left[D_x ([\ell_h V] [d_x W])\right]_{i,j}^{1} + \left[D_y ([\ell_h V] [d_y W])\right]_{i,j}^{1}\right).  \label{2D:PC-BCFD:Unique:d}
  \end{numcases}
\end{subequations}
with homogeneous boundary and initial conditions as \eqref{2D:CN-BCFD:IBc}.

First, taking discrete inner product of \eqref{2D:PC-BCFD:Unique:rho} with $\tau \bar{\mu}^{1}$, with the help of Lemma \ref{lemma:Dd}, Cauchy-Schwarz inequality and the fact $\|\bm{d} Z^{0}\|_\infty \le K_*+1$ , we obtain
\begin{equation*}
\begin{aligned}
        \|\bar{\mu}^{1}\|_{\rm M}^2  + \tau \| \bm{d}\bar{\mu}^{1} \|_{\rm TM}^2 
        & =  \lambda \tau \left( [\ell_h \bar{\mu}]^1 [d_x Z]^0, d_x \bar{\mu}^1 \right)_x 
        + \lambda \tau \left( [\ell_h \bar{\mu}]^1 [d_y Z]^0, d_y \bar{\mu}^1\right)_y\\
        & \leq \f{\tau \lambda^2\|\bm{d} Z^0\|_\infty^2}{4}   \|\bar{\mu}^{1}\|_{\rm M}^2
        + \tau\| \bm{d}\bar{\mu}^{1} \|_{\rm TM}^2
        \leq \f{\tau \lambda^2 (K_*+1)^2}{4}   \|\bar{\mu}^{1}\|_{\rm M}^2
        + \tau\| \bm{d}\bar{\mu}^{1} \|_{\rm TM}^2,
\end{aligned}
\end{equation*}
which implies that
\begin{equation*}
    \Big(1- \f{\tau \lambda^2 (K_*+1)^2}{4} \Big) \|\bar{\mu}^{1}\|_{\rm M}^2 \leq 0.
\end{equation*}
Then, for $\tau < \f{4}{\lambda^2(K_*+1)^2}$, we get
\begin{equation}\label{2D:PC-BCFD:Unique:e1}
    \|\bar{\mu}^{1}\|_{\rm M}^2 \leq 0 \Longrightarrow \bar{\mu}^{1} = 0 \Longrightarrow \bar{U}^1 =\bar{V}^1.
\end{equation}

Next, taking discrete inner product of \eqref{2D:PC-BCFD:Unique:c} with $\tau \theta^{1}$, and again applying Lemma \ref{lemma:Dd} and 
conclusion \eqref{2D:PC-BCFD:Unique:e1}, we directly have
\begin{equation}\label{uniqu:Z1}
        (1+\f{\tau}{2}) \|\theta^{1}\|_{\rm M}^2  + \f{\tau}{2} \| \bm{d}\theta^{1} \|_{\rm TM}^2 
        =  0 \Longrightarrow \theta^{1} = 0 \Longrightarrow Z^1 =W^1.
\end{equation}

Finally, taking discrete inner product of \eqref{2D:PC-BCFD:Unique:d} with  $\tau  \mu^{1}$, with the help of Lemma \ref{lemma:Dd}, and \eqref{2D:PC-BCFD:Unique:e1}--\eqref{uniqu:Z1}, 
we similarly obtain
\begin{equation*}
\begin{aligned}
        \|\mu^{1}\|_{\rm M}^2  + \f{\tau}{2} \| \bm{d}\mu^{1} \|_{\rm TM}^2 
        &=  \f{\lambda \tau}{2} \left( [\ell_h \mu]^1 [d_x Z]^1, d_x \mu^1 \right)_x 
        + \f{\lambda \tau}{2} \left( [\ell_h \mu]^1 [d_y Z]^1, d_y \mu^1\right)_y \\
        &\leq \f{\tau \lambda^2\|\bm{d} Z^1\|_\infty^2}{8}  \|\mu^{1}\|_{\rm M}^2
        + \f{\tau}{2} \| \bm{d}\mu^{1} \|_{\rm TM}^2\\
        &\le \f{\tau \lambda^2(K_*+1)^2}{8}  \|\mu^{1}\|_{\rm M}^2
        + \f{\tau}{2} \| \bm{d}\mu^{1} \|_{\rm TM}^2,
\end{aligned}
\end{equation*}
which implies that
\begin{equation}\label{uniqu:U1}
    \Big(1-\f{\tau \lambda^2(K_*+1)^2}{8}\Big) \|\mu^{1}\|_{\rm M}^2 \leq 0 \Longrightarrow \mu^{1} = 0,
\end{equation}
for $\tau < \f{8}{\lambda^2(K_*+1)^2}$.
Therefore, $U^1 =V^1$.
Thus, the approximate solutions $(\bar{U}^1, Z^1, U^1)$ of the prediction-correction BCFD scheme \eqref{2D:PC-BCFD} are unique for enough small $\tau$.

\paragraph{\bf Step II: Uniqueness of the approximate solutions $\{Z^{n}, U^{n}\}_{n=1}^{N}$}
By induction, we assume that $\theta^n = \mu^n = 0$ and we aim to show that $\theta^{n+1} = \mu^{n+1} = 0$. It follows that $( \theta^{n+1}, \mu^{n+1})$ are solutions of
 \begin{subequations}\label{2D:CN-BCFD:Unique}
	\begin{numcases}{}
  \f{1}{\tau} \theta_{i,j}^{n+1} = \f{1}{2}[D_x(d_x \theta)]_{i,j}^{n+1} + \f{1}{2}[D_y(d_y \theta)]_{i,j}^{n+1}
 	-  \f{1}{2}\theta_{i,j}^{n+1}, \label{2D:CN-BCFD:Unique:c}\\
  		\f{1}{\tau} \mu_{i,j}^{n+1} = \f{1}{2}[D_x(d_x \mu)]_{i,j}^{n+1} + \f{1}{2}[D_y(d_y \mu)]_{i,j}^{n+1} \nonumber\\
		\qquad\qquad\qquad -\f{\lambda}{2} \left( \left[D_x ([\ell_hU] [d_x Z])\right]_{i,j}^{n+1} + \left[D_y ([\ell_h U] [d_y Z])\right]_{i,j}^{n+1}\right) \nonumber\\
        \qquad\qquad\qquad +\f{\lambda}{2} \left( \left[D_x ([\ell_h V] [d_x W])\right]_{i,j}^{n+1} + \left[D_y ([\ell_h V [d_y W])\right]_{i,j}^{n+1}\right).  \label{2D:CN-BCFD:Unique:rho}
  \end{numcases}
\end{subequations}

It is apparent that the two equations in \eqref{2D:CN-BCFD:Unique} are identical to \eqref{2D:PC-BCFD:Unique:c}--\eqref{2D:PC-BCFD:Unique:d}, then by taking discrete inner products  with $\tau \theta^{n+1}$ and $\tau  \mu^{n+1}$, respectively, we obtain  similar results to \eqref{uniqu:Z1} and \eqref{uniqu:U1} that
\begin{equation}\label{uniqu:Zn}
        (1+\f{\tau}{2}) \|\theta^{n+1}\|_{\rm M}^2  + \f{\tau}{2} \| \bm{d}\theta^{n+1} \|_{\rm TM}^2 
        =  0 
        \Longrightarrow \theta^{n+1} = 0 \Longrightarrow Z^{n+1}=W^{n+1},
\end{equation}
and
\begin{equation}\label{uniqu:Un}
    \Big(1-\f{\tau \lambda^2(K_*+1)^2}{8}\Big) \|\mu^{n+1}\|_{\rm M}^2 \leq 0
    \Longrightarrow \mu^{n+1} = 0 \Longrightarrow U^{n+1}=V^{n+1},
\end{equation}
for $\tau < \f{8}{\lambda^2(K_*+1)^2}$.

In summary,  for $\tau < \f{4}{\lambda^2 (K_*+1)^2}$,  we have $U^{n}=V^{n}$ and $Z^{n}=W^{n}$  $n=1,\ldots, N$. This completes the proof of uniqueness of solutions, hence the existence of solutions.
\end{proof}

\section{Numerical results}\label{sec:num}
In this section, some numerical experiments using the DeC-MC-BCFD scheme \eqref{2D:PC-BCFD}--\eqref{2D:CN-BCFD:IBc} are carried out. Accuracy of the scheme under both uniform and non-uniform spatial grids is demonstrated in subsection \ref{subsec:accuracy}. Meanwhile, various scenarios for the blow-up phenomenon are simulated with different input initial values in subsection \ref{subsec:blow}.

Unless otherwise specified, the non-uniform spatial grids utilized in this paper are generated as follows.
First, a uniform partition $x_{{\rm fix},i+1/2}=a^x+ ih_{\rm fix}$ for $i=0,\ldots, N_x$ with equal grid size $ h_{\rm fix}=(b^x-a^x)/N_x$ is constructed. Then, through a slight random adjustment of the grid size using Matlab inline code \textit{rand}, we define the non-uniform grid points as follows:
\begin{equation}\label{non-uni}
	x_{i+1/2} = x_{{\rm fix},i+1/2} + \beta \cdot h_{\rm fix}\cdot (-1+2\, rand(i)),  \quad i=1, \ldots, N_x-1,
\end{equation}
where $\beta$ represents a small mesh parameter that can regulate random disturbance within a specific range. The $y$-direction non-uniform grid points $y_{j+1/2}~(j=0, \ldots, N_y)$ are defined in a similar manner. For all the tests below, we set $N_x = N_y = M$ and  take $\lambda=1$.

\subsection{Accuracy test}\label{subsec:accuracy}
We start by checking the accuracy of the DeC-MC-BCFD scheme \eqref{2D:PC-BCFD}--\eqref{2D:CN-BCFD:IBc} for the Keller--Segel system \eqref{model:ks} subject to homogeneous Neumann boundary conditions.

\begin{example}\label{exam:s1}
In this example, the computational domain is set to be $\Omega=(0,1) \times(0,1)$ and the exact solutions are chosen as
\begin{align*}
	\rho(x, y, t)=(x^2-x )^2 (y^2-y )^2t, \quad	c(x, y, t)=(x^2-x )^2 (y^2-y )^2t. 
\end{align*}
\end{example}

\begin{figure}[!ht]
	\centering
	\includegraphics[width=0.325\linewidth]{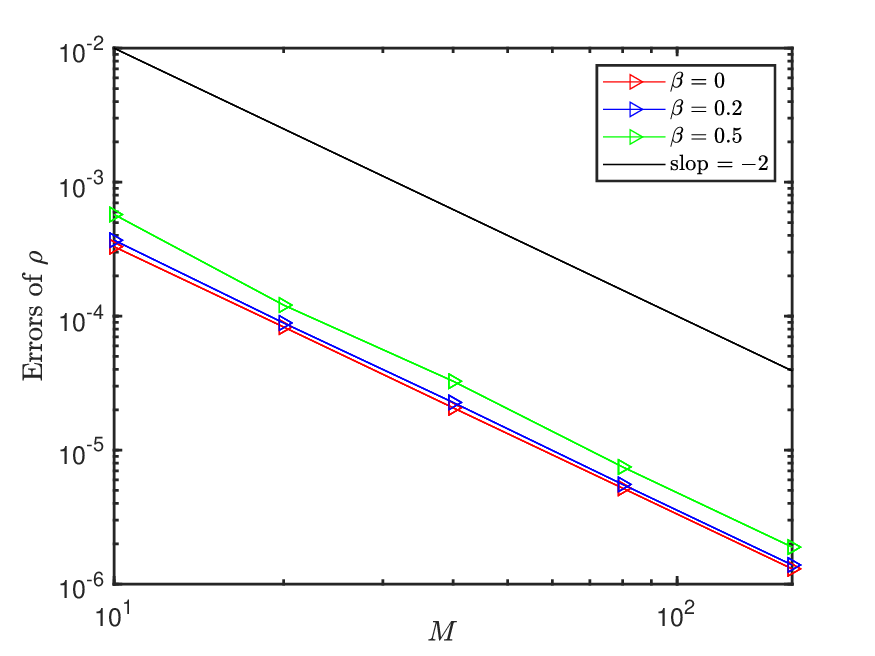}
	\includegraphics[width=0.325\linewidth]{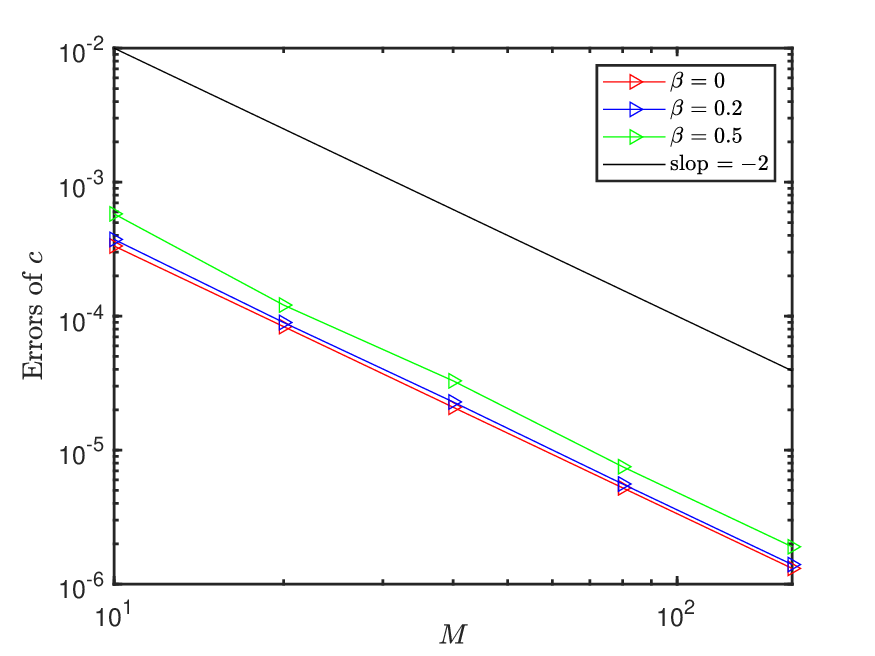}
    \includegraphics[width=0.325\linewidth]{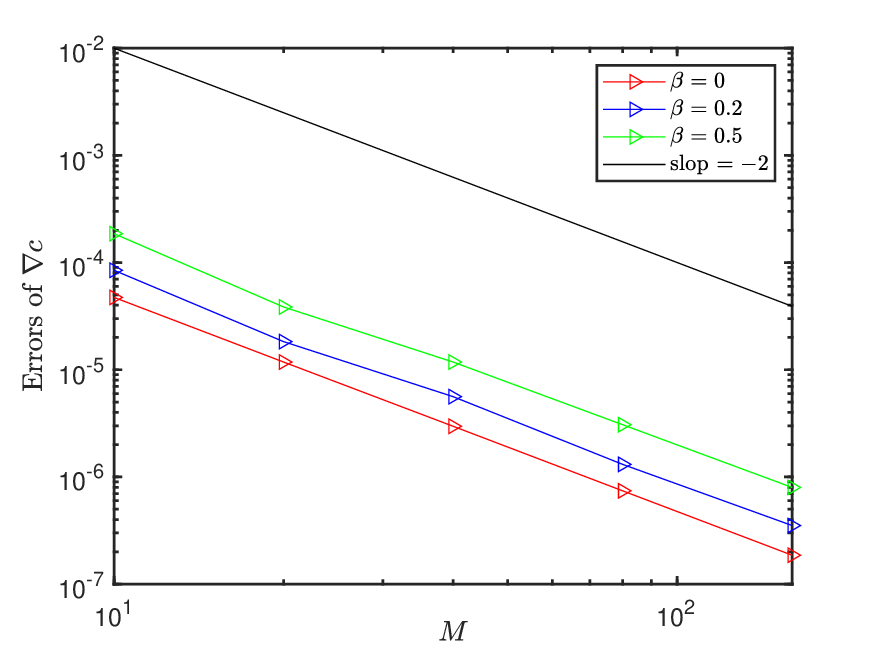}
	\caption{ $L^2$ errors of $\rho$, $c$ and $\nabla c$ for the DeC-MC-BCFD scheme for Example \ref{exam:s1}.}
	\label{fig:err}
\end{figure}

We test numerically the spatial and temporal accuracy by setting the grid sizes $\tau=h_{\rm fix}$, in which the discrete $L^2$ and $H^1$ errors are measured at $T=1$.
Table \ref{tab:accuracy} presents the errors and convergence orders for the cell distribution $\rho$, the chemoattractant concentration $c$ and its gradient $\nabla c$, utilizing the DeC-MC-BCFD scheme \eqref{2D:PC-BCFD}--\eqref{2D:CN-BCFD:IBc}. These results are further illustrated in Figure \ref{fig:err}.
In conclusion, the following observations can be drawn:
\begin{enumerate}[(i)]
    \item Our proposed scheme achieves second-order convergence in both time and space, demonstrated on both uniform and non-uniform spatial grids. These findings align closely with the conclusions drawn in Theorem \ref{thm:coverg}.
    \item Despite increasing random disturbances in the spatial grids, i.e., with larger mesh parameter $\beta$, errors escalate while still consistently maintain the same order of magnitude and second-order accuracy.
\end{enumerate}

\begin{table}[htbp]
	\centering \caption{ $L^2$ and $H^1$ errors for the DeC-MC-BCFD scheme for Example \ref{exam:s1}. } \label{tab:accuracy}
	\setlength{\tabcolsep}{1.5mm}
\begin{tabular}{c| c c c c c c c }
	\toprule
	Mesh  & $M$ &  $\|\rho-U\|_{\rm M}$&  Order&  $\|c-Z\|_{\rm M}$& Order &  $\|\nabla c-\bm{d} Z\|_{\rm TM}$& Order  \\
	\midrule
	\multirow{5}*{\makecell{ $\beta = 0$}}
	&10 & 3.30e-04 & --    & 3.34e-04 & --    & 4.73e-05 &--\\
	&20 & 8.30e-05 & 1.99  & 8.36e-05 & 2.00  & 1.18e-05 & 1.99\\
	&40 & 2.07e-05 & 2.00  & 2.09e-05 & 2.00  & 2.97e-06 & 2.00\\
	&80 & 5.20e-06 & 2.00  & 5.23e-06 & 2.00  & 7.42e-07 & 2.00\\
	&160 & 1.30e-06 & 2.00 & 1.31e-06 & 2.00  & 1.86e-07 & 2.00\\
 \midrule
 \multirow{5}*{\makecell{ $\beta = 0.2$}}
	&10  & 3.69e-04 & --   & 3.74e-04 & --   & 8.47e-05 & --\\
	&20  & 8.90e-05 & 2.05 & 8.97e-05 & 2.06 & 1.83e-05 & 2.05\\
	&40  & 2.27e-05 & 1.97 & 2.29e-05 & 1.97 & 5.61e-06 & 1.97\\
	&80  & 5.55e-06 & 2.03 & 5.59e-06 & 2.03 & 1.31e-06 & 2.03\\
	&160 & 1.39e-06 & 1.99 & 1.40e-06 & 1.99 & 3.52e-07 & 1.99\\
		\midrule
  \multirow{5}*{\makecell{ $\beta = 0.5$}}
		&10  & 5.74e-04 & --   & 5.80e-04 & --    & 1.86e-04 & --\\
		&20  & 1.21e-04 & 2.25 & 1.21e-04 & 2.26  & 3.85e-05 & 2.25\\
		&40  & 3.27e-05 & 1.88 & 3.29e-05 & 1.88  & 1.18e-05 & 1.88\\
		&80  & 7.49e-06 & 2.13 & 7.52e-06 & 2.13  & 3.07e-06 & 2.13\\
		&160 & 1.89e-06 & 1.99 & 1.90e-06 & 1.99  & 7.98e-07 & 1.99\\
	\bottomrule
\end{tabular}
\end{table}

\subsection{Investigation of the blow-up phenomenon}\label{subsec:blow}
It is well known that the solutions to the Keller--Segel system \eqref{model:ks} exist globally if $Mass_\rho (0) < 8 \pi$, and undergo finite-time blow-up if $Mass_\rho(0) > 8 \pi$ \cite{CC'08, LWZ'18}. In this subsection, we employ the proposed scheme \eqref{2D:PC-BCFD}--\eqref{2D:CN-BCFD:IBc} to simulate various blow-up phenomena.

\begin{example}\label{exam:less8pi}
	 (Global existence with $Mass_\rho (0) \approx 24.67<8 \pi$.) According to the existence condition of the Keller--Segel system, i.e., the total mass of cells is strictly less than $8 \pi$, we set the initial cell density and chemoattractant concentration as follows:
	\begin{align*}
		&\rho(x, y, 0)=50 \exp \left(-5\big((x-0.5)^2+(y-0.5)^2\big)\right),\\
		&c(x, y, 0)=25 \exp \left(-5\big((x-0.5)^2+(y-0.5)^2\big)/2\right) ,
	\end{align*}
	where the total mass of cells $Mass_\rho (0) \approx 24.67$ in $\Omega=(0,1) \times(0,1)$. 
\end{example}

\begin{figure}[!ht]
	\centering
	\includegraphics[width=0.325\linewidth]{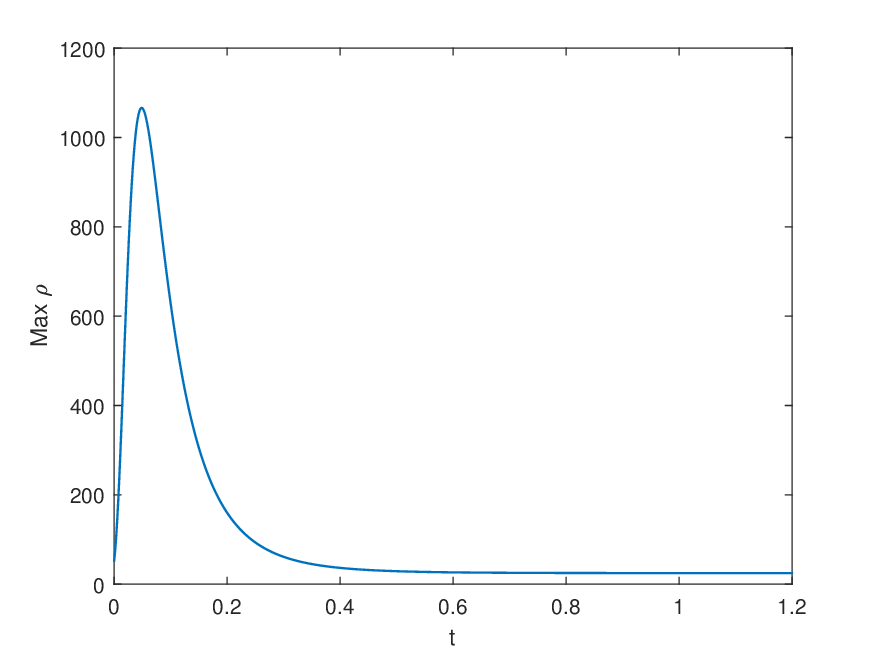}
	\includegraphics[width=0.325\linewidth]{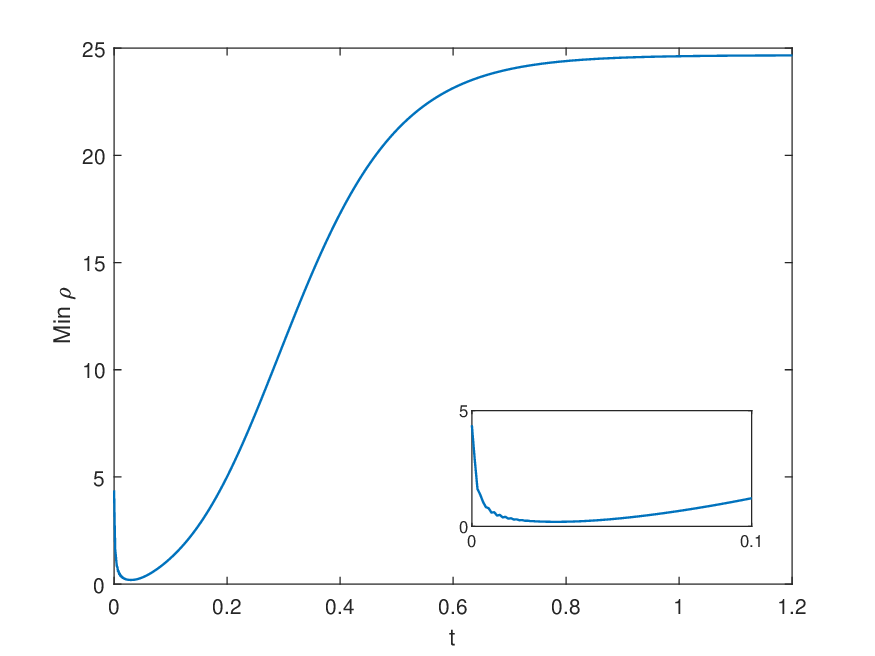}
	\includegraphics[width=0.325\linewidth]{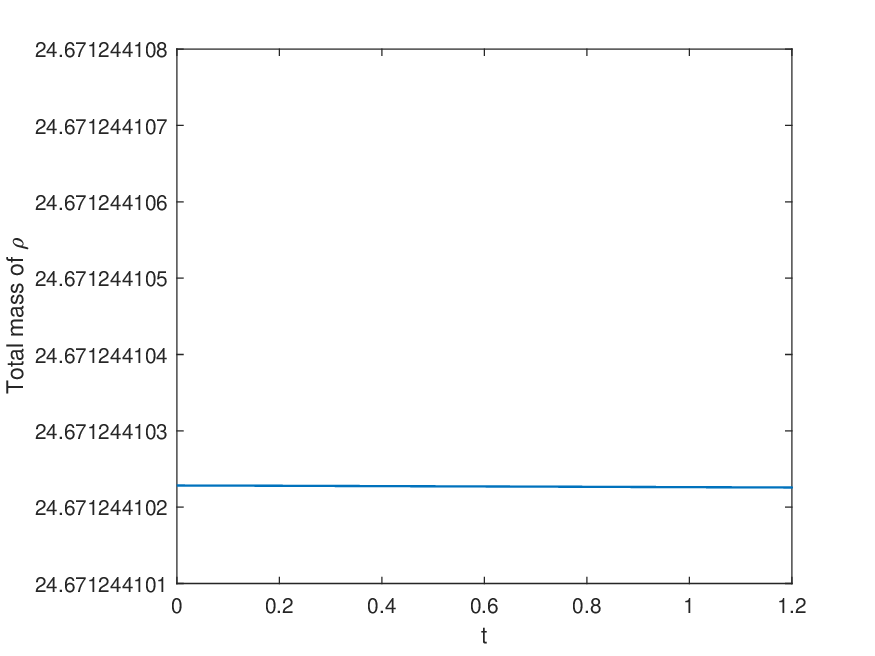}
	\caption{ Evolutions of the maximum, minimum  and total mass of $\rho$ on uniform grids for Example \ref{exam:less8pi}. }
	\label{fig:uni_bdl8}
\end{figure}
\begin{figure}[!ht]
	\centering
	\includegraphics[width=0.325\linewidth]{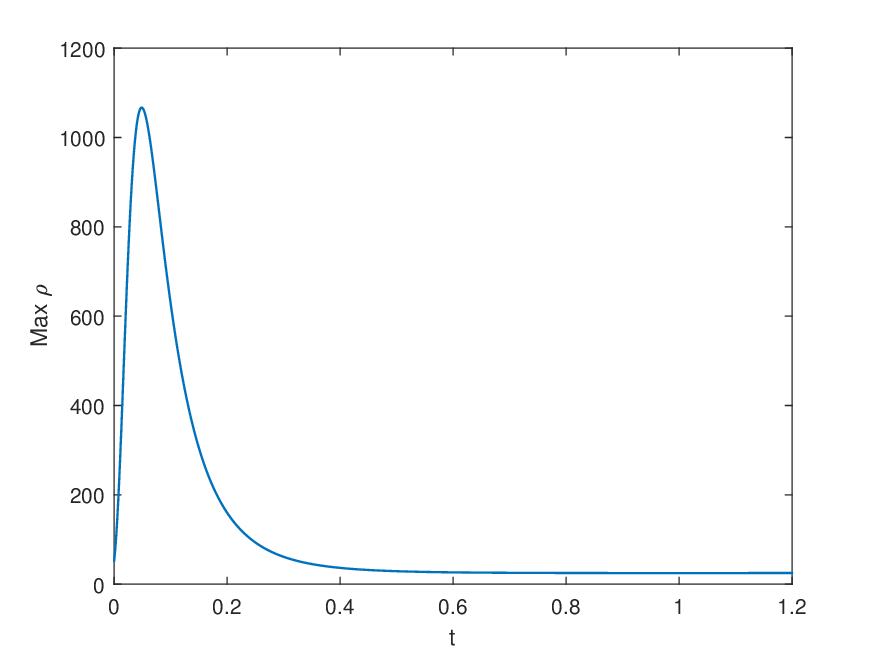}
	\includegraphics[width=0.325\linewidth]{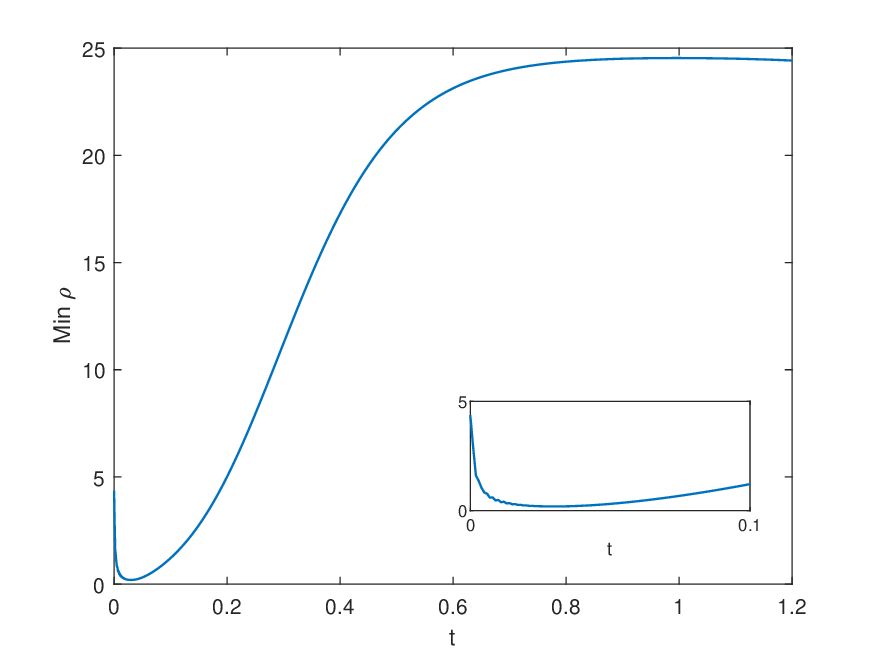}
	\includegraphics[width=0.325\linewidth]{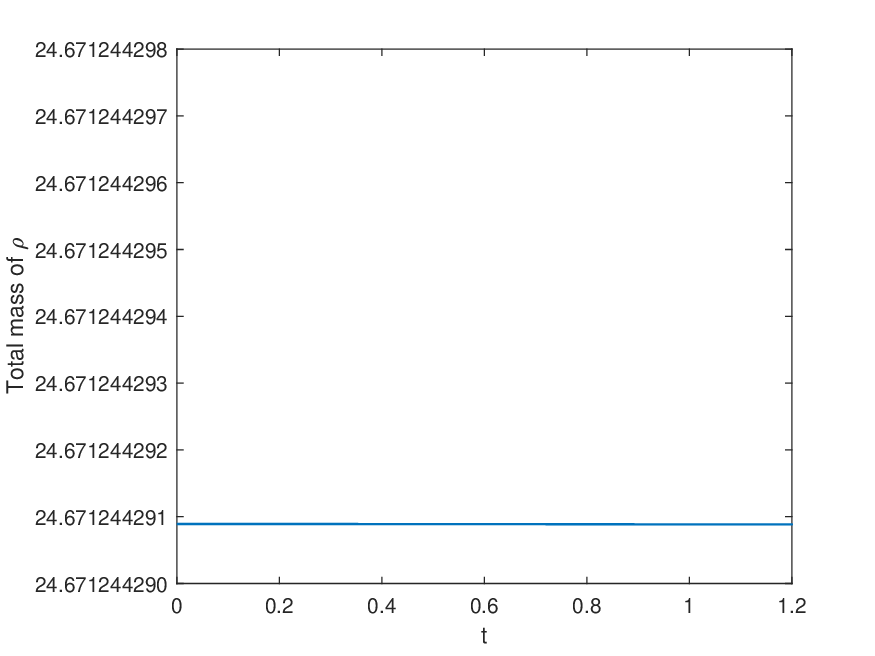}
	\caption{Evolutions of the maximum, minimum and total mass of $\rho$ on non-uniform  grids for Example \ref{exam:less8pi}. }
	\label{fig:nonuni_bdl8}
\end{figure}
\begin{figure}[!ht]
	\centering
	\includegraphics[width=0.325\linewidth]{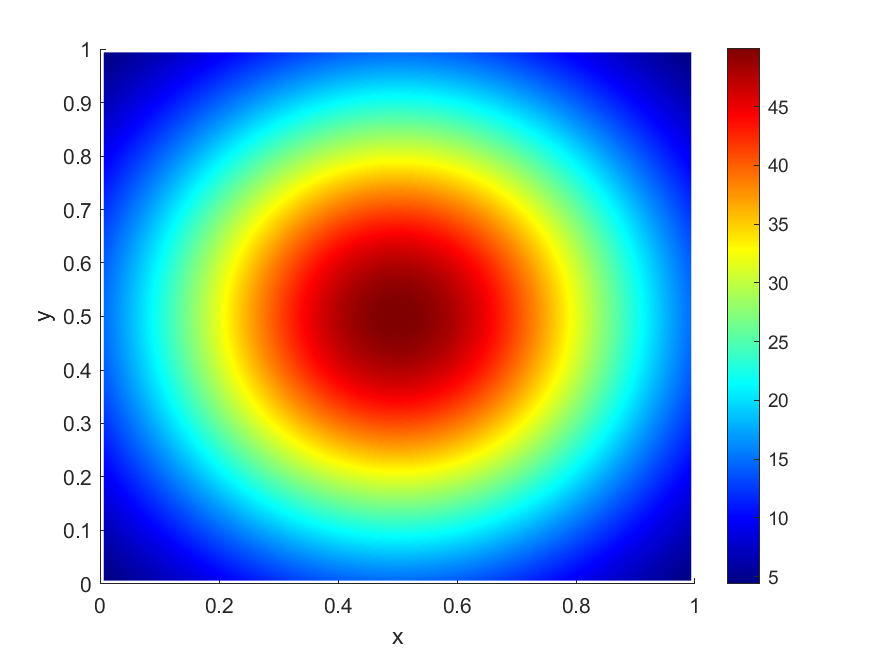}
	\includegraphics[width=0.325\linewidth]{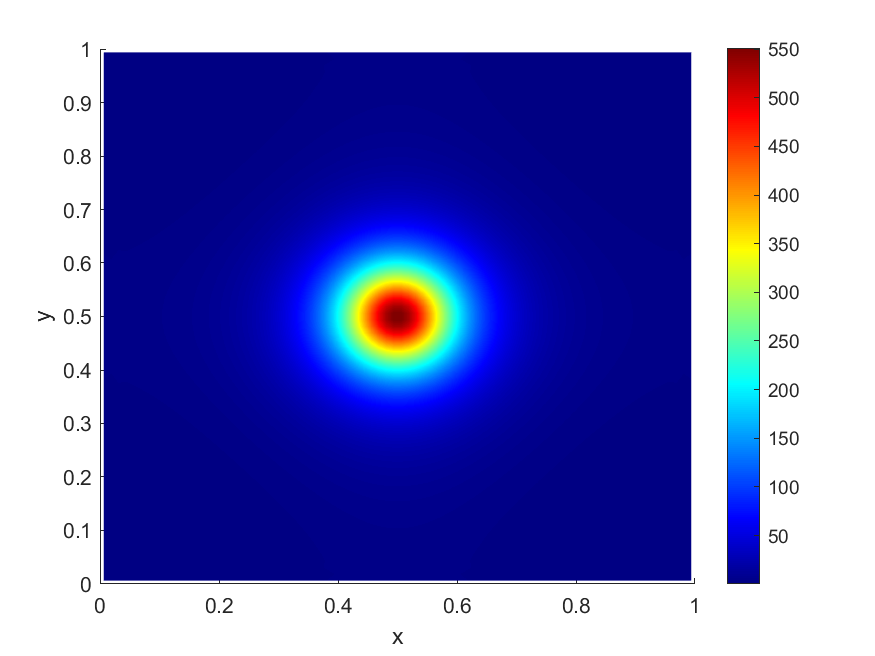}
	\includegraphics[width=0.325\linewidth]{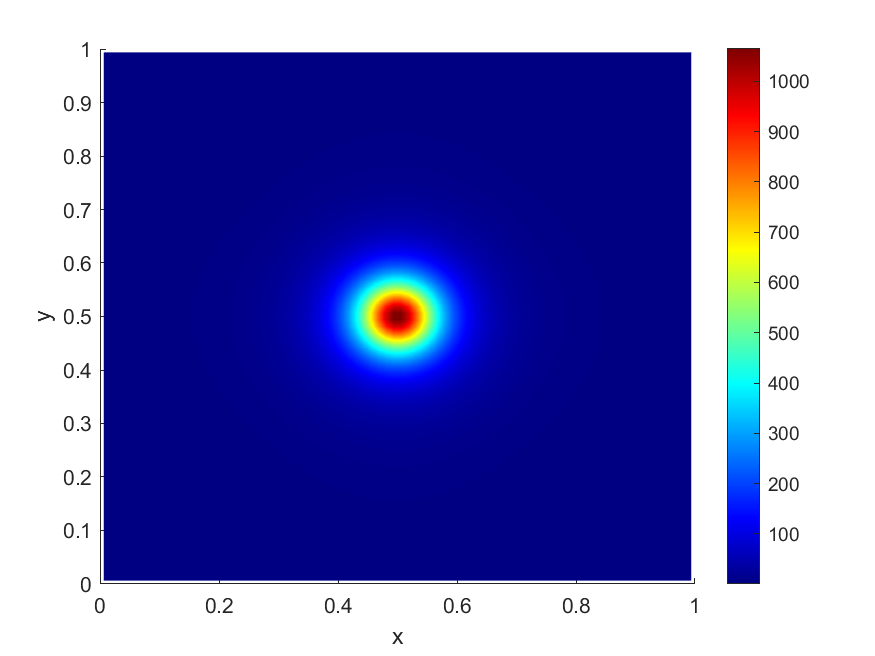}
    \includegraphics[width=0.325\linewidth]{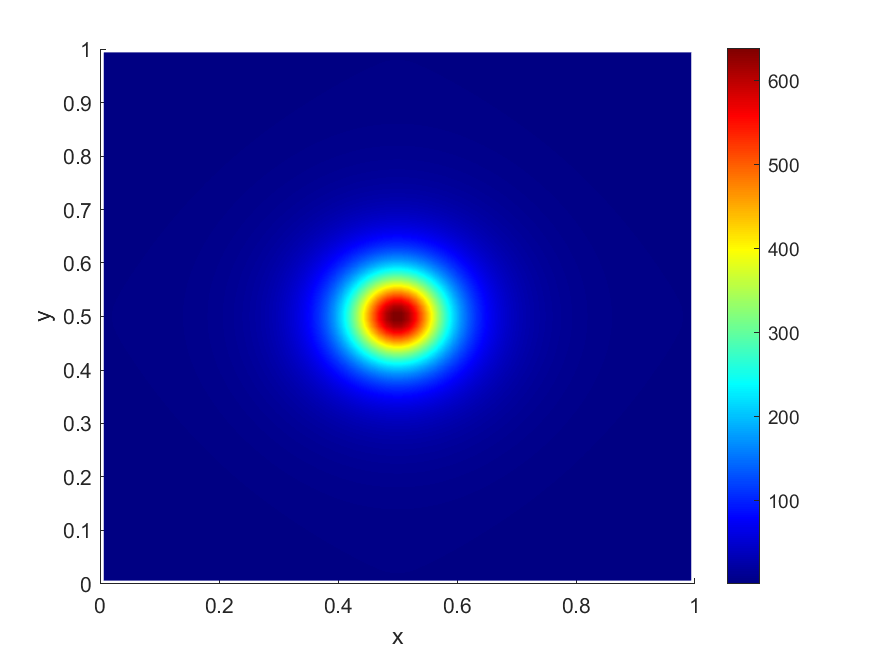}
	\includegraphics[width=0.325\linewidth]{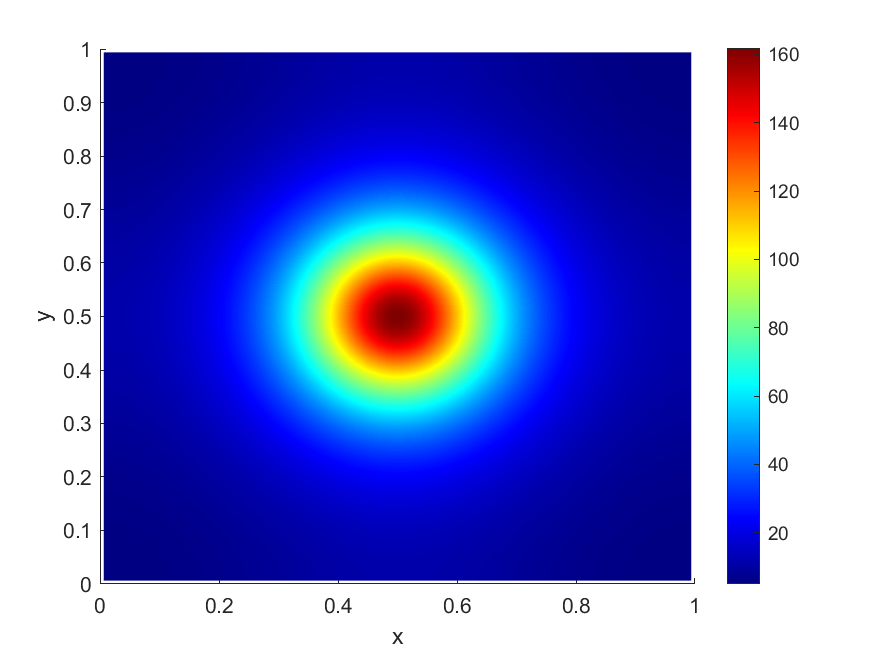}
	\includegraphics[width=0.325\linewidth]{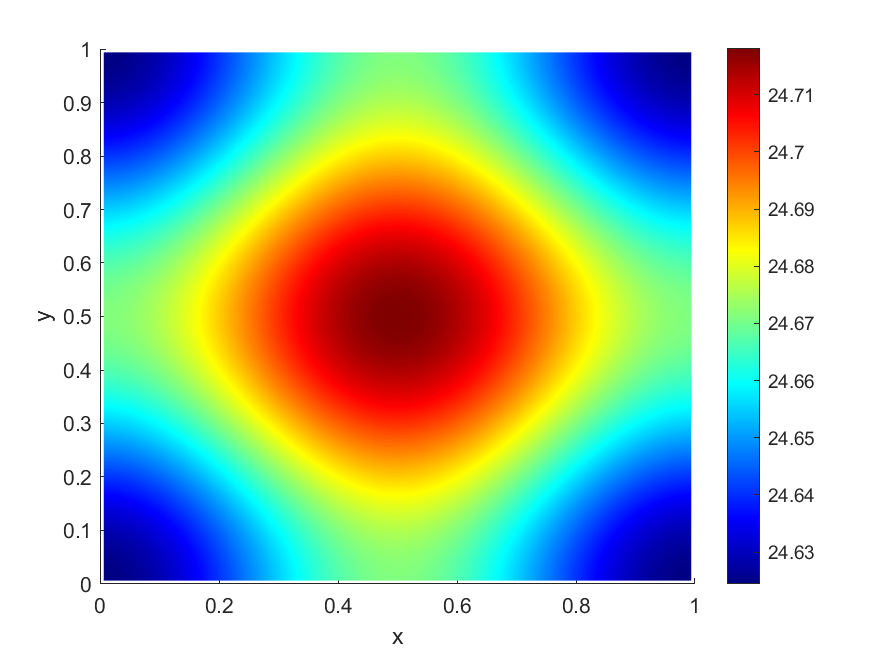}
	\caption{Contour plots of $\rho$ at time instants $t = 0,\ 2.0\times10^{-2},\ 4.9\times10^{-2},\ 0.1,\ 0.2,\ 1$ (from left to right) for Example \ref{exam:less8pi}. }
	\label{fig:evolution_bdl8_rho}
\end{figure}
\begin{figure}[!ht]
	\centering
	\includegraphics[width=0.325\linewidth]{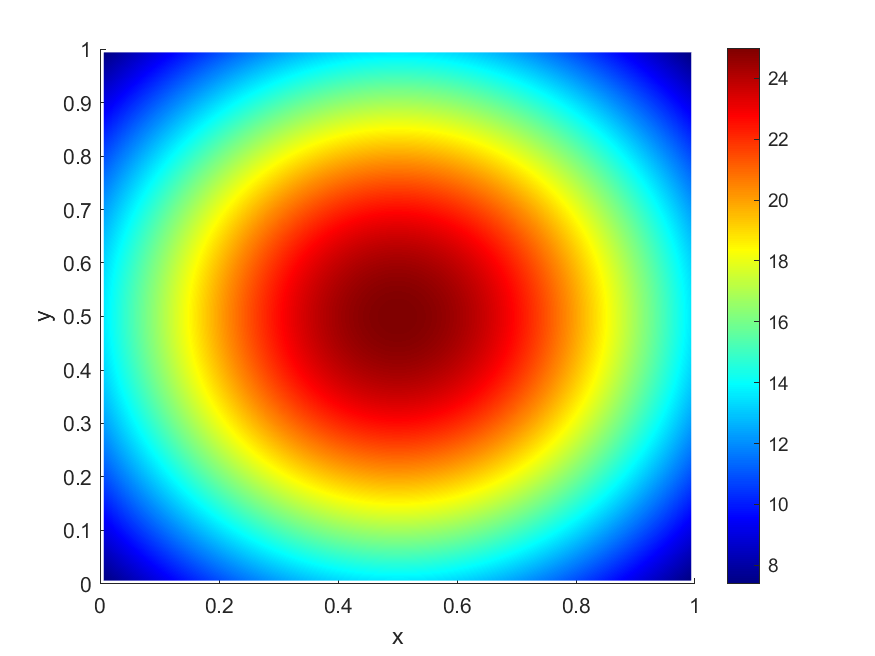}
	\includegraphics[width=0.325\linewidth]{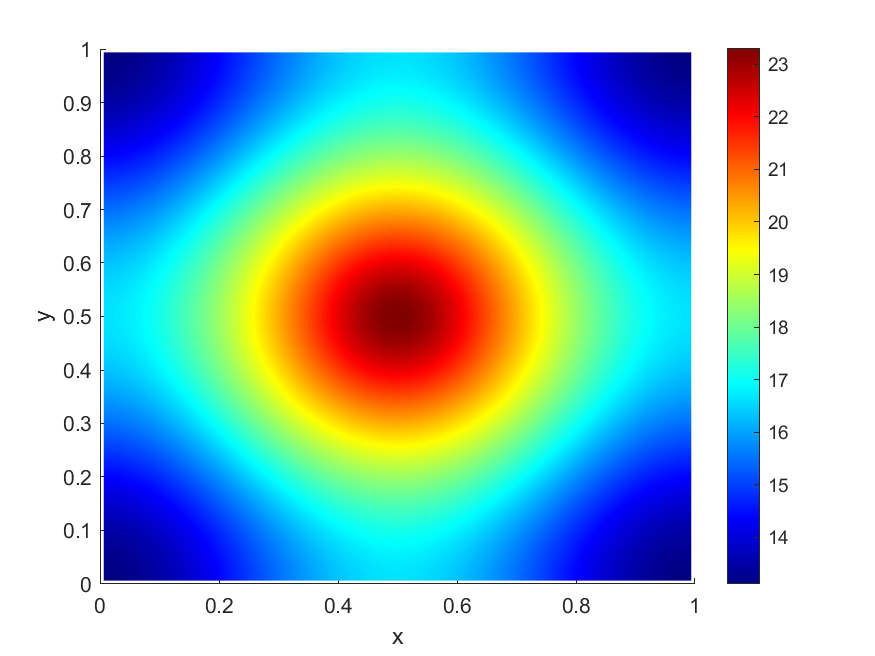}
	\includegraphics[width=0.325\linewidth]{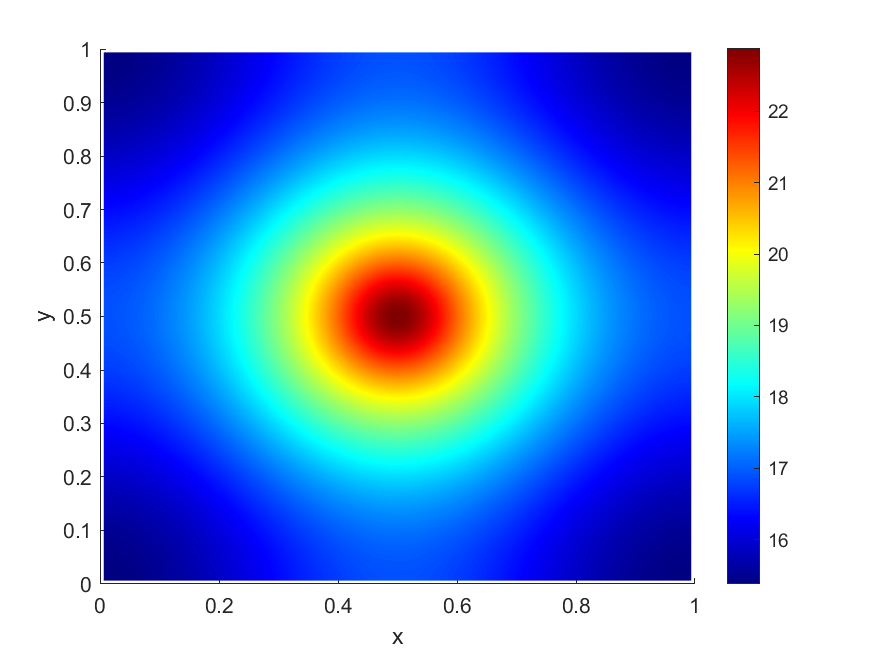}
        \includegraphics[width=0.325\linewidth]{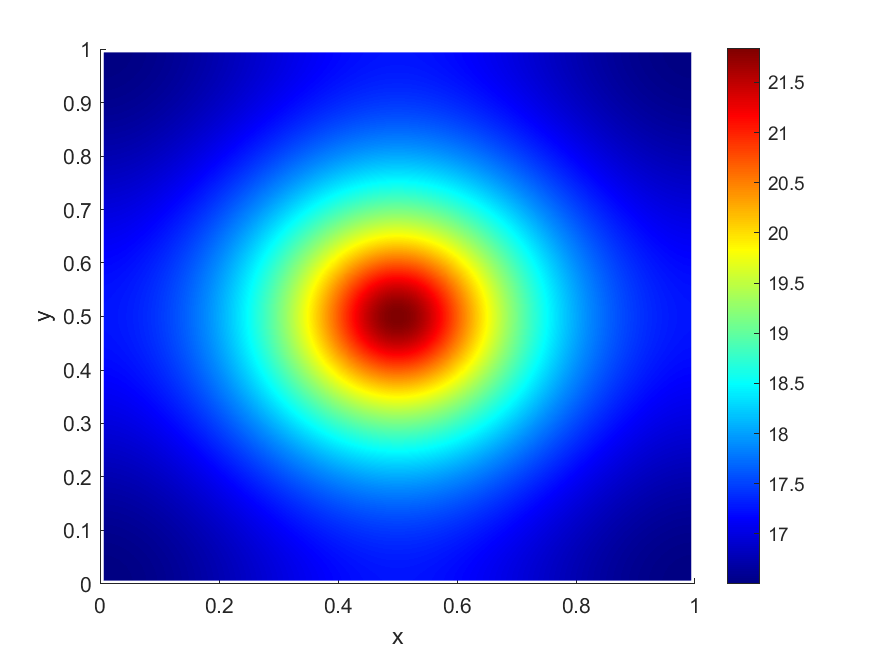}
	\includegraphics[width=0.325\linewidth]{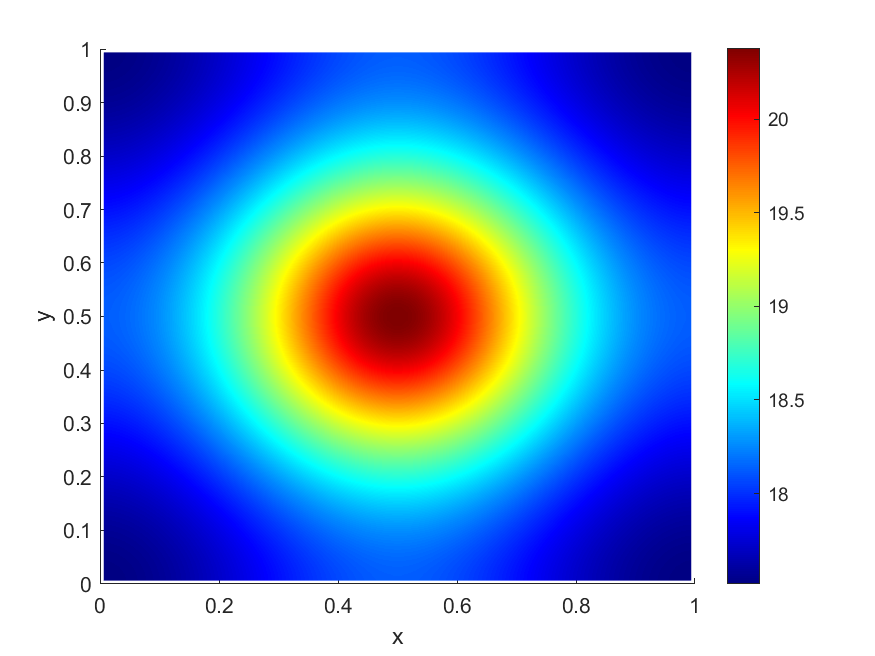}
	\includegraphics[width=0.325\linewidth]{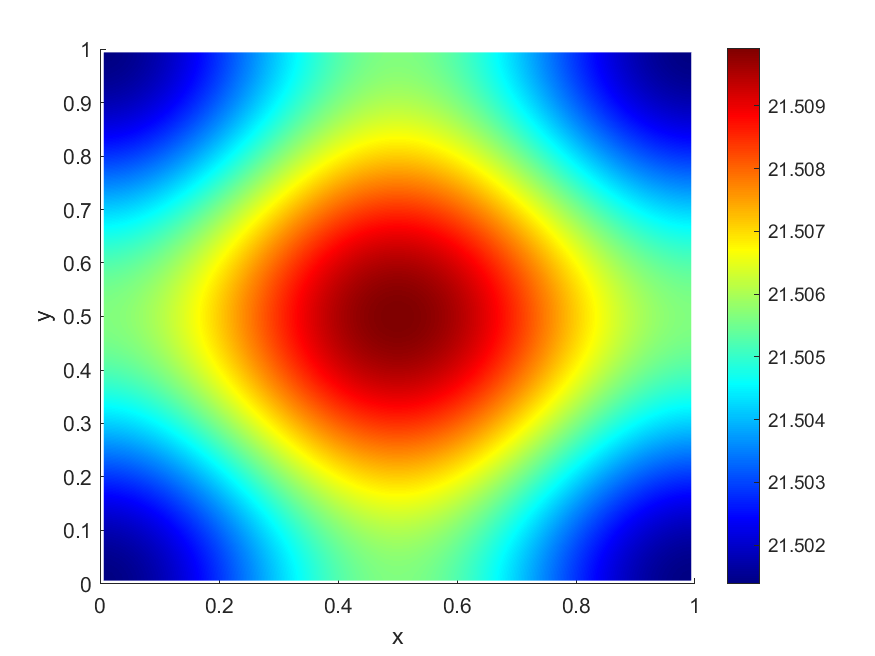}
	\caption{Contour plots of $c$ at time instants $t = 0,\ 2.0\times10^{-2},\ 4.9\times10^{-2},\ 0.1,\ 0.2,\ 1$ (from left to right) for Example \ref{exam:less8pi}. }
	\label{fig:evolution_bdl8_c}
\end{figure}

In this simulation, we set $M = 80$ and $\tau = 10^{-3}$. The evolutions of the maximum $\rho$, minimum $\rho$ and total mass of $\rho$ under both uniform and non-uniform spatial grids with $\beta=0.1$ are depicted in Figures \ref{fig:uni_bdl8}--\ref{fig:nonuni_bdl8}. Figures \ref{fig:evolution_bdl8_rho}--\ref{fig:evolution_bdl8_c} show the contour plots of $\rho$ and $c$ at different time instants $t = 0,\ 2.0\times10^{-2},\ 4.9\times10^{-2},\ 0.1,\ 0.2,\ 1$, respectively. Key observations include:
\begin{enumerate}[(i)]
    \item The maximum value of $\rho$ initially increases and then gradually decreases to its steady state, whether on uniform or non-uniform spatial grids.
    \item It shows that  the cell density $\rho$ remains non-negative on both sets of grids, consistent with the positivity-preserving requirement of the solution to the Keller--Segel system.
    \item The total mass of $\rho$ is conserved in the scenario that the solutions to the Keller--Segel system globally exist, demonstrating the  mass conservation law as proved in Theorem \ref{thm:MassConserve}.
    \item Despite small random perturbations in the spatial grids, the simulation remains robust and comparable to that observed on uniform grids. This underscores the superiority of the proposed scheme in maintaining accuracy on general non-uniform grids. In fact, the ability of preserving second-order accuracy on non-uniform grids opens the possibility of taking advantage of adaptive grids to further improve the accuracy of the solution.
\end{enumerate}

\begin{example}\label{exam:grate8pi}
	(Blow-up with $Mass_\rho (0) \approx 27.23>8 \pi$.) In order to explore the chemotactic blow-up, we increase the total mass of cell density to $Mass_\rho (0) \approx 27.23$ by changing the initial cell density and chemoattractant concentration as
	\begin{align*}
		\rho(x, y, 0)=130 \exp \left(-15 (x^2+y^2 )\right),\		c(x, y, 0)=13 \exp \left(-2 (x^2+y^2 )\right),
	\end{align*}
	in the domain $\Omega=(-1,1) \times(-1,1)$. 
\end{example}

\begin{figure}[!ht]
	\centering
	\includegraphics[width=0.4\linewidth]{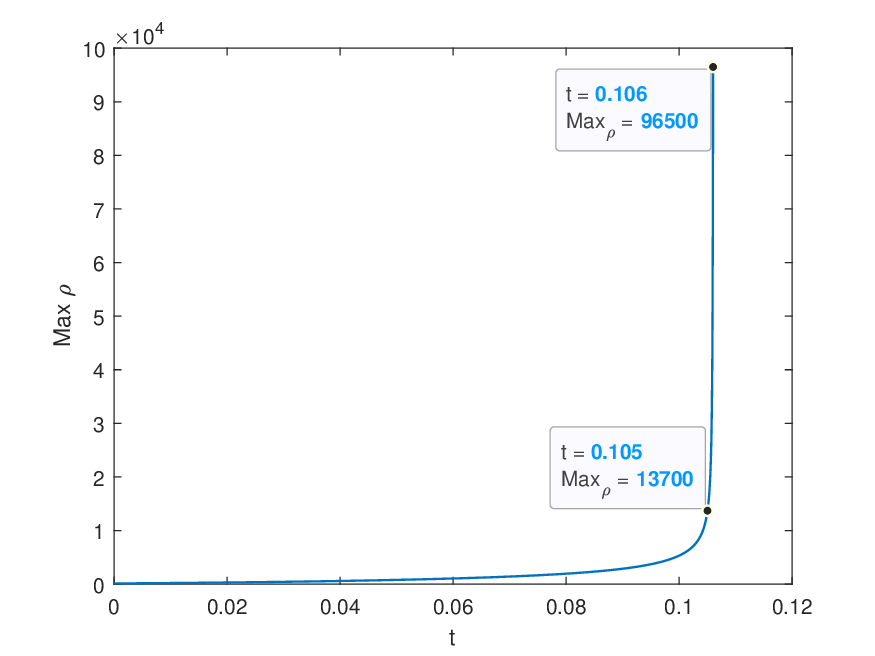}
	\includegraphics[width=0.4\linewidth]{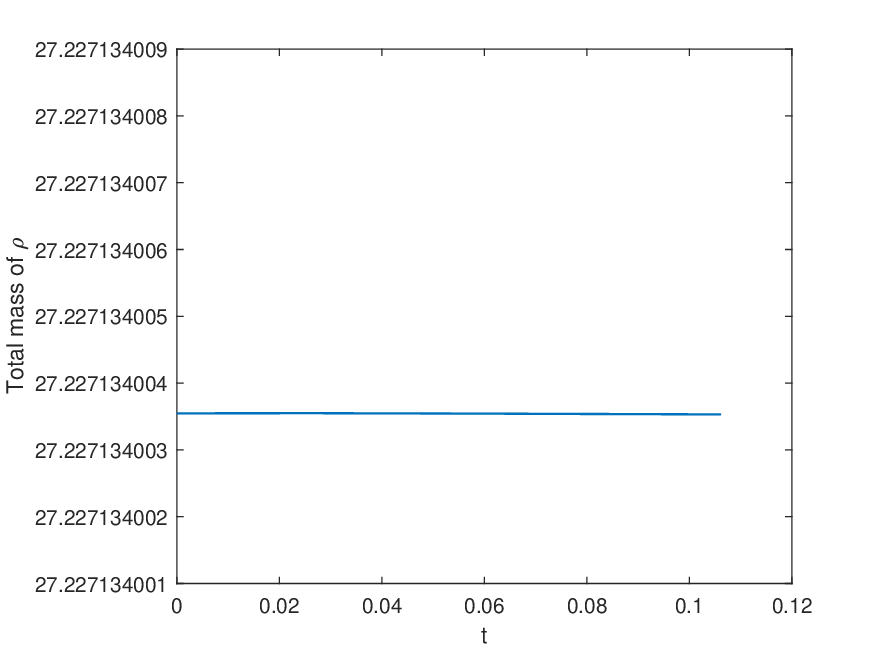}
	\caption{Evolutions of the maximum and total mass of $\rho$ on uniform grids for Example \ref{exam:grate8pi}.}
	\label{fig:uni_bug8}
\end{figure}
To ensure the reliability of the numerical results, we set $M = 80$ and  an even small time stepsize $\tau=10^{-5}$. 
The evolutions of the maximum value and total mass of $\rho$ using uniform spatial grids are illustrated in Figure \ref{fig:uni_bug8}. We observe that the maximum value of $\rho$ is increasing with respect to time and eventually exhibits finite-time blow-up, while still maintaining the principle of mass conservation.

At this stage, we refrain from employing general non-uniform spatial grids with random disturbances for numerical simulations due to their potential adverse impact on accurately capturing the blow-up phenomena. In the subsequent examples, we plan to employ specific non-uniform spatial grids to accurately model the rapid blow-up phenomenon.

\begin{example}\label{exam:bu:center} (Blow-up at the center of unit square.)
 We consider the Keller--Segel model \eqref{model:ks} in a square domain $\Omega=(0,1) \times(0,1)$. The initial conditions are given as
\begin{equation*}
	\begin{aligned}
		& \rho(x, y, 0)=1000 \exp\left( -100 ((x-0.5)^2+(y-0.5)^2 )\right) , \\
		& c(x, y, 0)=500 \exp\left( -50 ((x-0.5)^2+(y-0.5)^2 )\right)  .
	\end{aligned}
\end{equation*}	
\end{example}

\begin{figure}[!ht]
	\centering
	\includegraphics[width=0.4\linewidth]{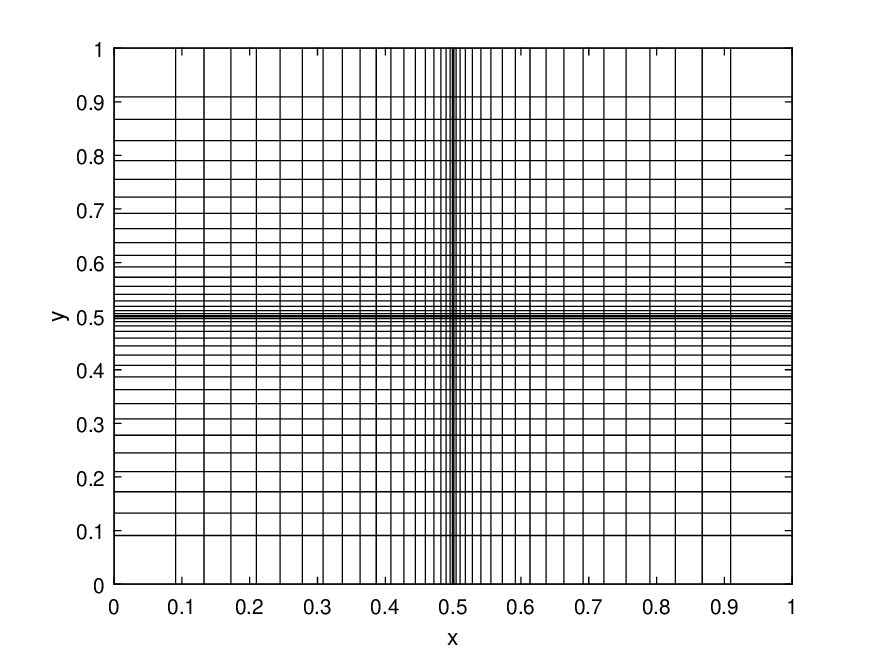}
		\caption{ Middle refinement grids.}
	\label{fig:gridmid}
\end{figure}

According to Refs. \cite{AGR'23, CE'18}, the cell density of \eqref{model:ks} with the above given initial conditions has a peak at the center $(0.5, 0.5)$, where the blow-up phenomena for $\rho$ occurs rapidly in finite time. Thus, in this simulation, the following specified non-uniform grids are applied:
\begin{equation}\label{grid:mid}
\begin{cases}
	x_{N_{x}/2+i+1/2}=\f{1}{2} + \frac{i^2}{2\left(N_x / 2+1\right)^2}, & i=0,1, \ldots, N_x / 2+1, \\ 
    x_{N_{x}/2-i+1/2}=\f{1}{2} -\frac{i^2}{2\left(N_x / 2+1\right)^2}, & i=1,2, \ldots, N_x / 2+1,
\end{cases}
\end{equation}
and $y_{j+1/2}~(j=0,\ldots, N_y)$ are similarly defined. Here, for simplicity, $N_x=N_y=M$ are chosen as even numbers. The resulting mesh with $M=40$ is shown in Figure \ref{fig:gridmid}.

\begin{figure}[!ht]
	\centering
	\includegraphics[width=0.4\linewidth]{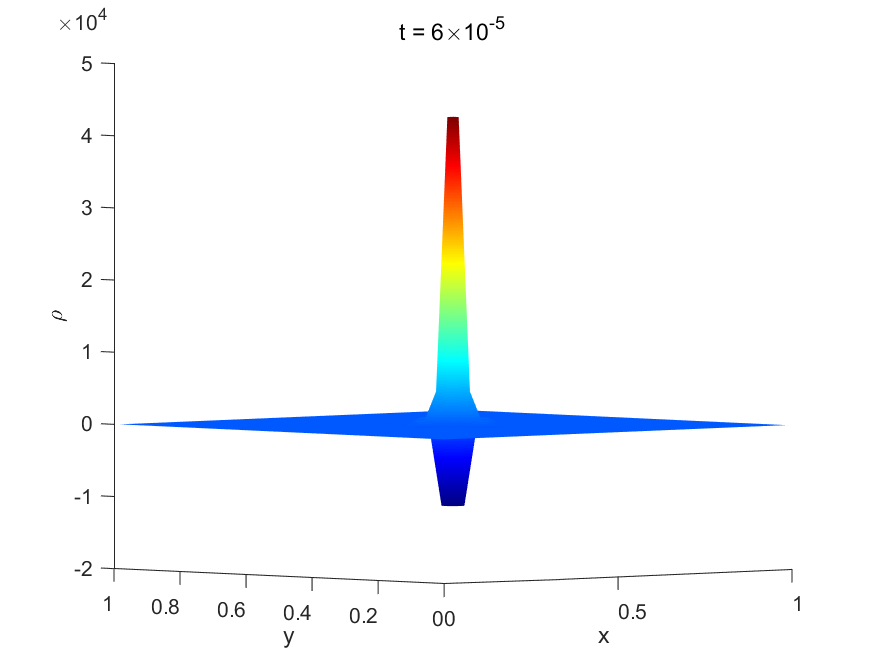}
	\includegraphics[width=0.4\linewidth]{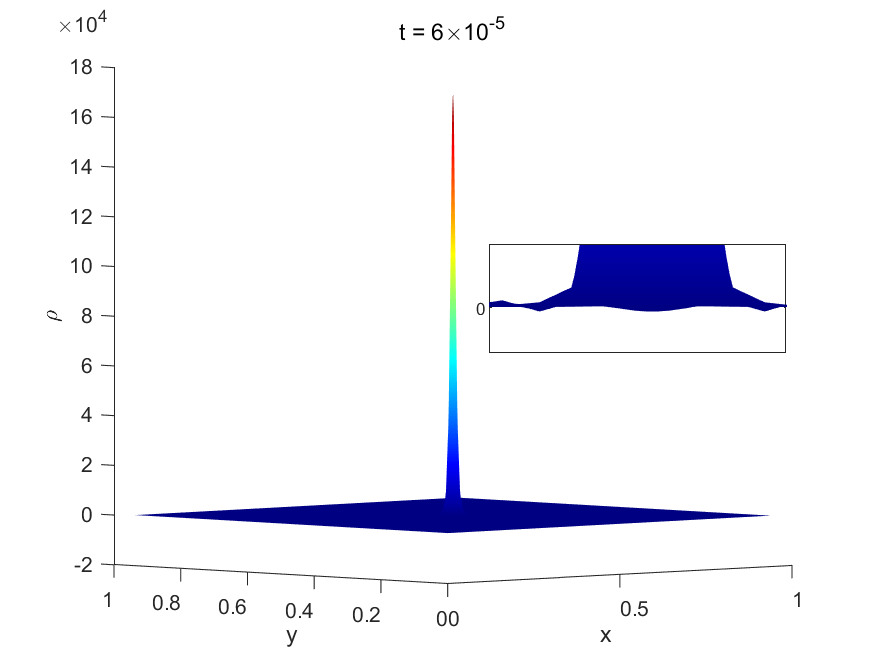}
	\includegraphics[width=0.4\linewidth]{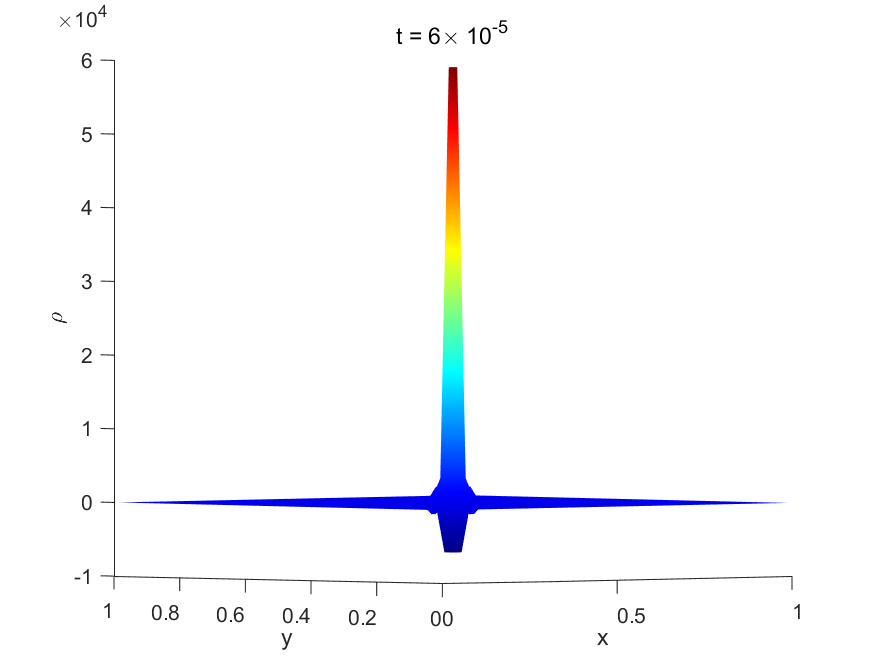}
	\includegraphics[width=0.4\linewidth]{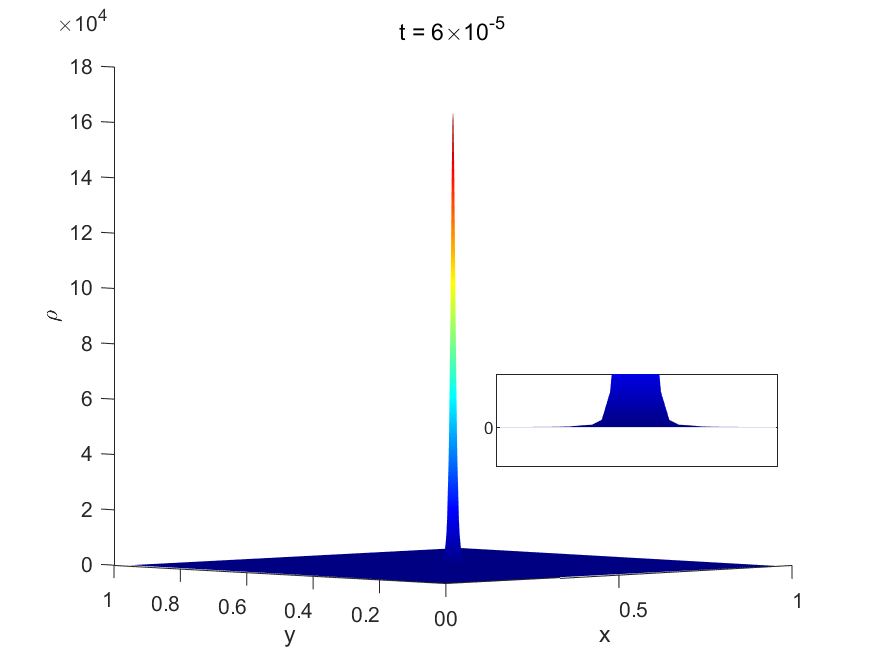}
	\includegraphics[width=0.4\linewidth]{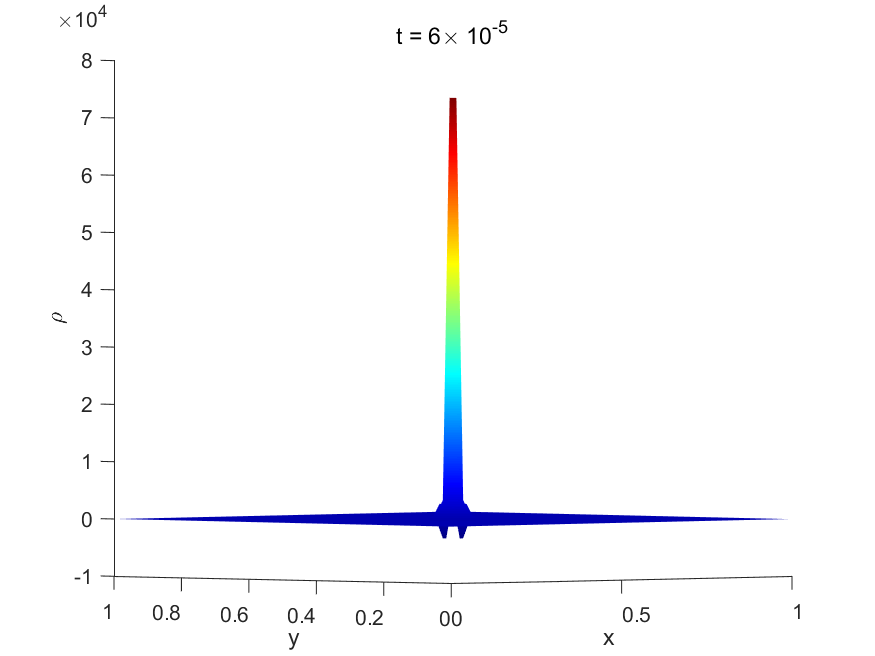}
	\includegraphics[width=0.4\linewidth]{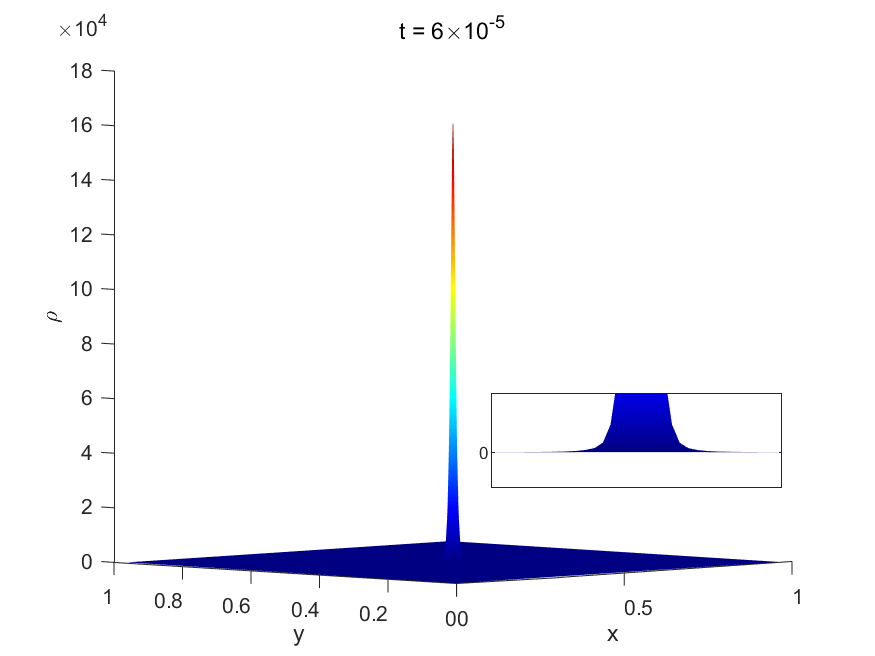}
	\includegraphics[width=0.4\linewidth]{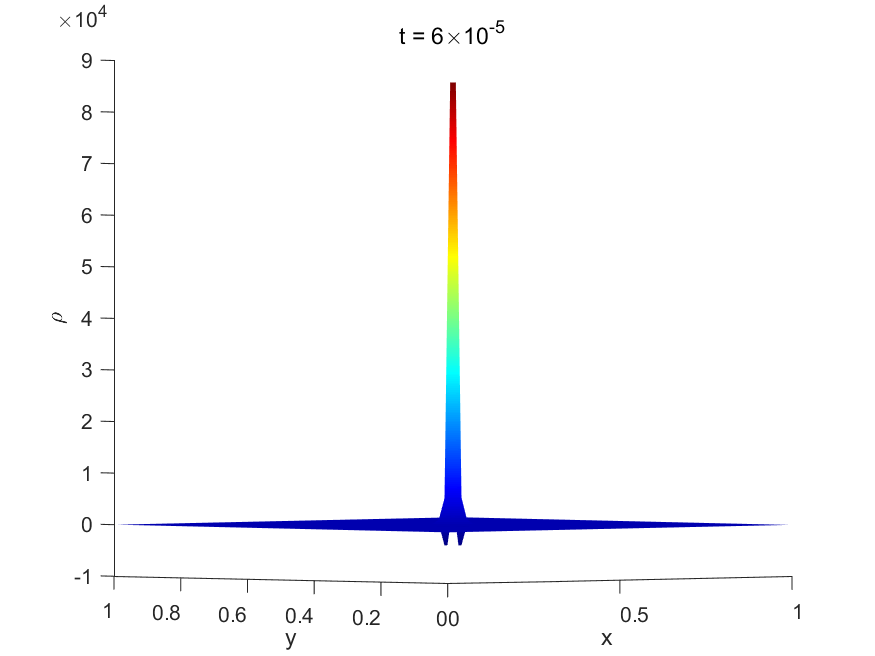}
	\includegraphics[width=0.4\linewidth]{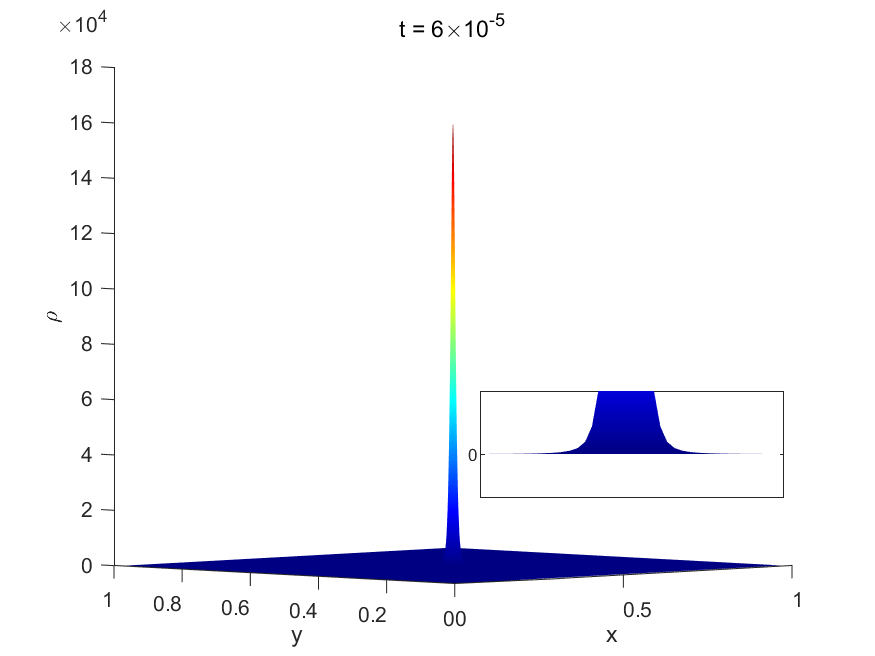}
	\caption{ Blow-up phenomena at $t = 6 \times 10^{-5}$ on uniform grids (left) and non-uniform grids (right) using the DeC-MC-BCFD scheme, each row corresponds to $M = 60,\ 80,\ 100,\ 120$.}
	\label{fig:blowupnon0}
\end{figure}
As the cell density shall blow up at a very short finite time, we set $\tau=10^{-6}$. Blow-up phenomena at $t = 6 \times 10^{-5}$ are tested using both uniform and non-uniform grids \eqref{grid:mid} for the DeC-MC-BCFD scheme \eqref{2D:PC-BCFD}--\eqref{2D:CN-BCFD:IBc} with various $M = 60, 80, 100, 120$, as depicted in Figure \ref{fig:blowupnon0}. It is observed that
\begin{enumerate}[(i)]
    \item Compared to simulation on the uniform spatial grids, that on the non-uniform grids shows a marked improvement in the representation of blow-up phenomena. In particular, numerical results on the non-uniform grids with $M = 60$ surpass those on uniform grids even with $M = 120$. Thus, the specified non-uniform grids can  offer enhanced resolution in specific regions, which demonstrate its better performance in achieving more accurate and reliable simulations. 

    \item The peak values of the cell density $\rho$ on specified non-uniform grids with $M = 60, 80, 100, 120$ exhibit remarkable similarity. This observation indicates that fewer non-uniform grid points can effectively simulate the blow-up phenomena, highlighting the efficiency of the DeC-MC-BCFD scheme \eqref{2D:PC-BCFD}--\eqref{2D:CN-BCFD:IBc} on non-uniform grids. Therefore, the proposed scheme is crucial for accurately capturing the blow-up dynamics, particularly in complex real-world applications where computational resources are often limited. 
 
    \item Well-chosen non-uniform grids can reduce the non-positivity violations. The solutions to the Keller--Segel system obtained using the specific non-uniform grids with $M = 80, 100, 120$ can be observed to be non-negative, while negative solutions appear when using uniform grids even with $M = 120$. However, this is not proved theoretically and requires further discussions.
\end{enumerate} 

\begin{example}\label{exam:bu:corner} (Blow-up at the corner of the rectangular region.)
In this example, we consider the Keller--Segel chemotaxis system \eqref{model:ks} in a rectangular region $(-0.5,0.5) \times (-0.5,0.5)$ with the following initial conditions:
\begin{equation*}
		 \rho(x, y, 0)=1000 \exp\left(-100 ((x-0.15)^2+(y-0.15)^2 )\right), ~~	c(x, y, 0)=0.
\end{equation*}	
\end{example}

\begin{figure}[!ht]
	\centering
	\includegraphics[width=0.4\linewidth]{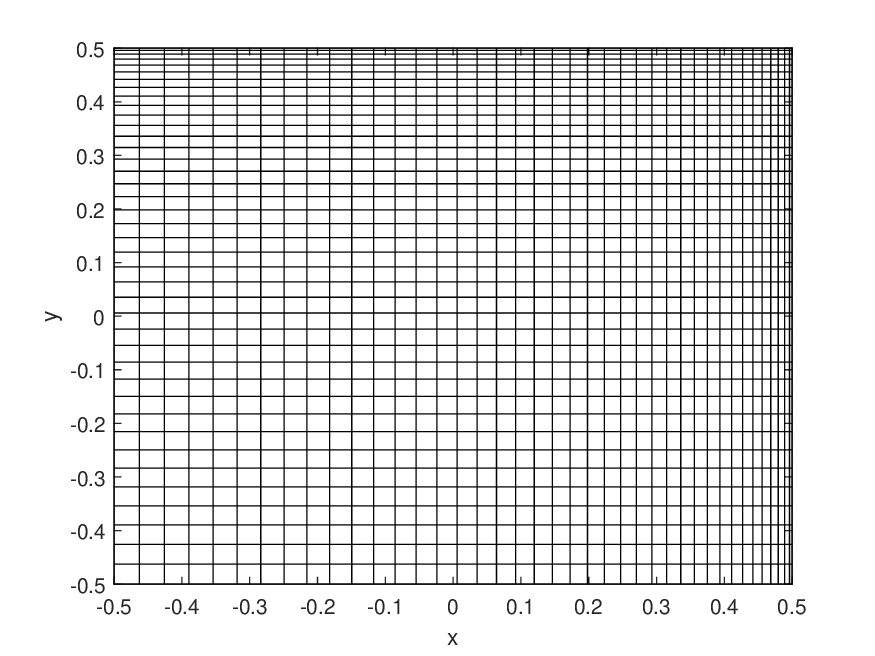}	\caption{ Corner refinement grids.}
	\label{fig:gridcor}
\end{figure}

As stated in Ref. \cite{EK'12}, the cell density of \eqref{model:ks} with the above given initial conditions will blow up in finite time at the corner (0.5, 0.5) of the rectangular region. Therefore, in this simulation, we choose the following specified non-uniform grids:
\begin{equation}\label{grid:cor}
\begin{cases}
	x_{M-i+1/2}= \f{1}{2} -(\f{i}{M})^{3/2}, \quad i=0,1, \ldots, M,\\
 y_{M-j+1/2}= \f{1}{2} -(\f{j}{M})^{3/2}, \quad j=0,1, \ldots, M.
 \end{cases}
\end{equation}
The resulting mesh with $M=40$ is shown in Figure \ref{fig:gridcor}.

\begin{figure}[!ht]
	\centering
	\includegraphics[width=0.32\linewidth]{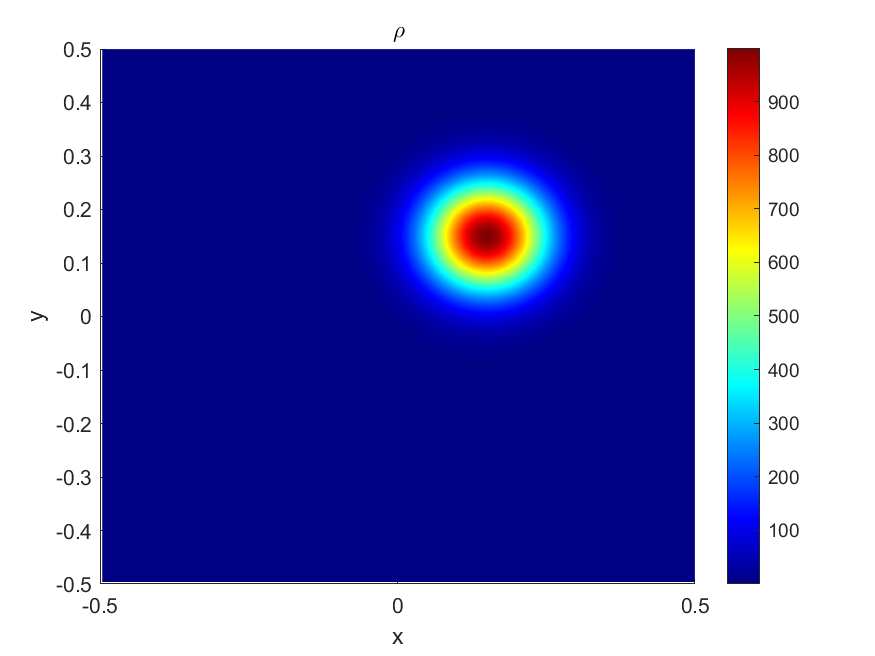}
	\includegraphics[width=0.32\linewidth]{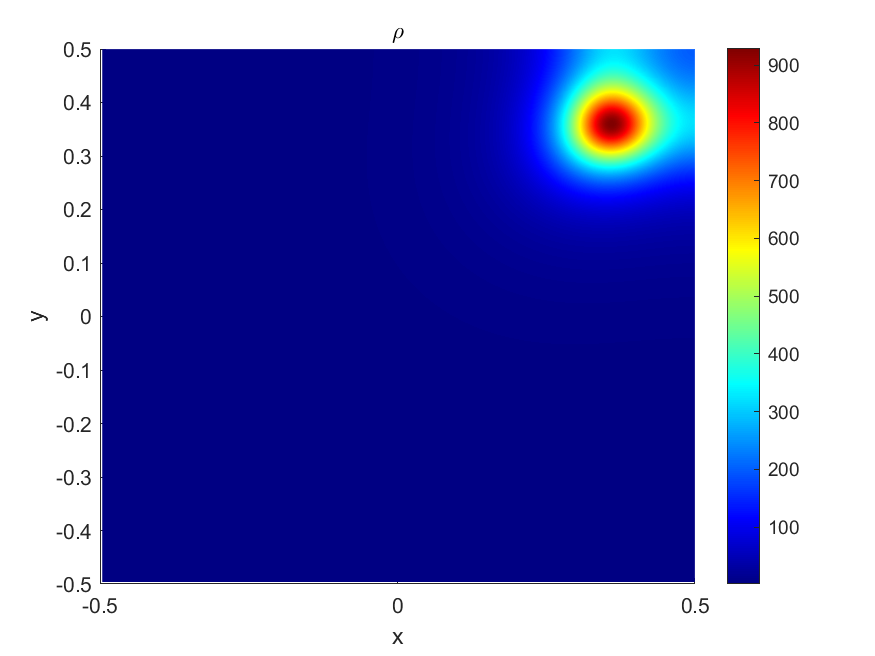}
	\includegraphics[width=0.32\linewidth]{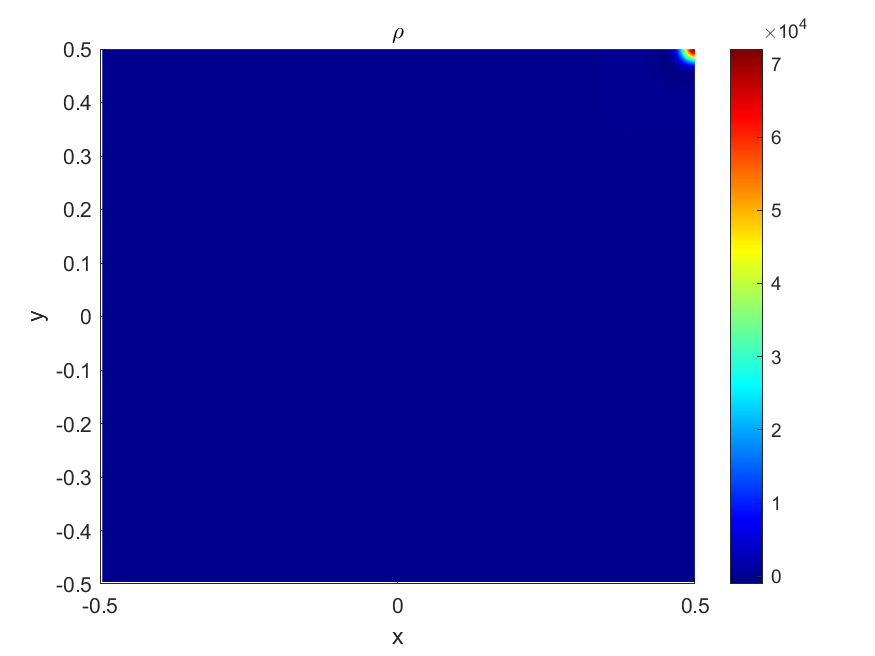}
	\includegraphics[width=0.32\linewidth]{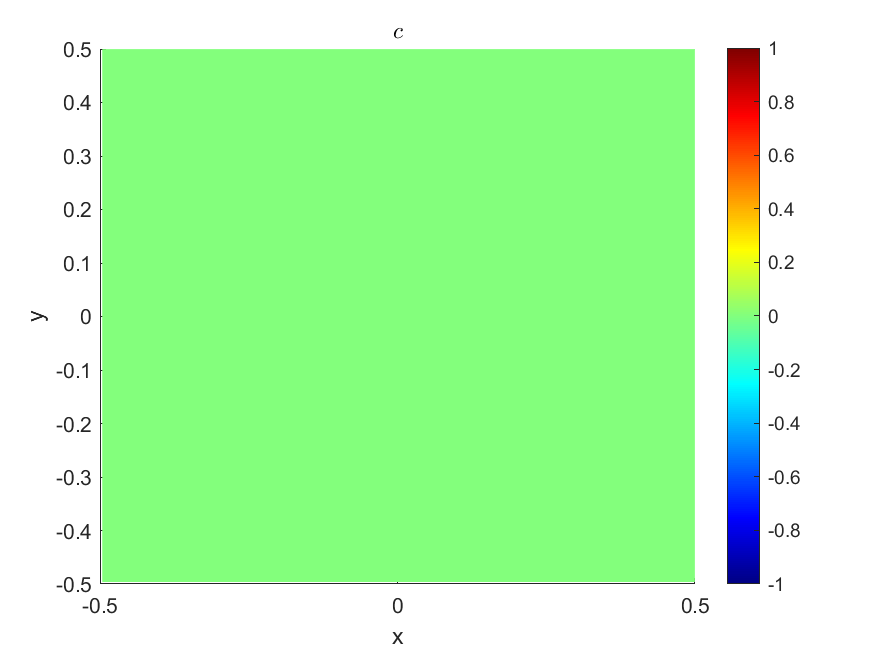}
	\includegraphics[width=0.32\linewidth]{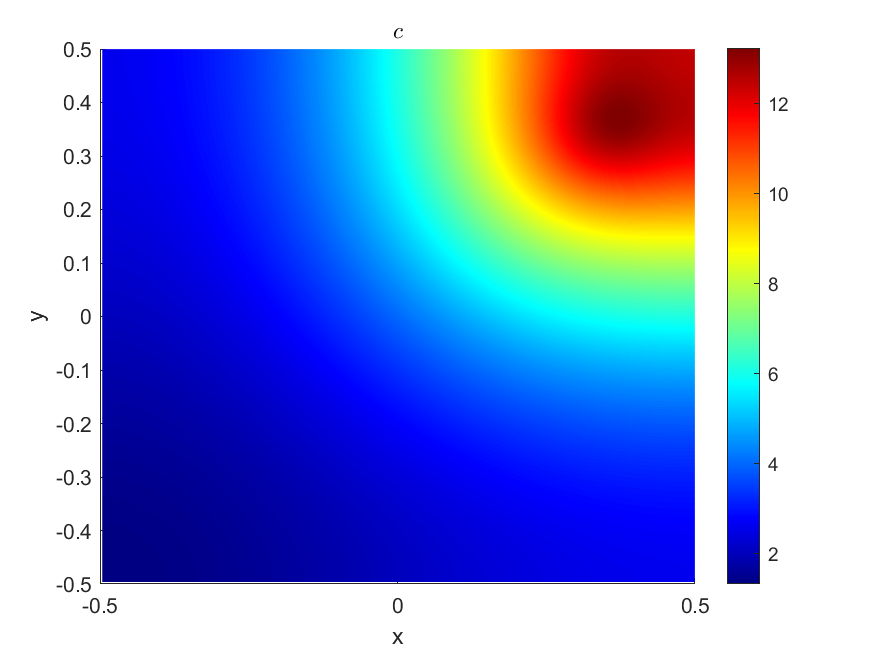}
	\includegraphics[width=0.32\linewidth]{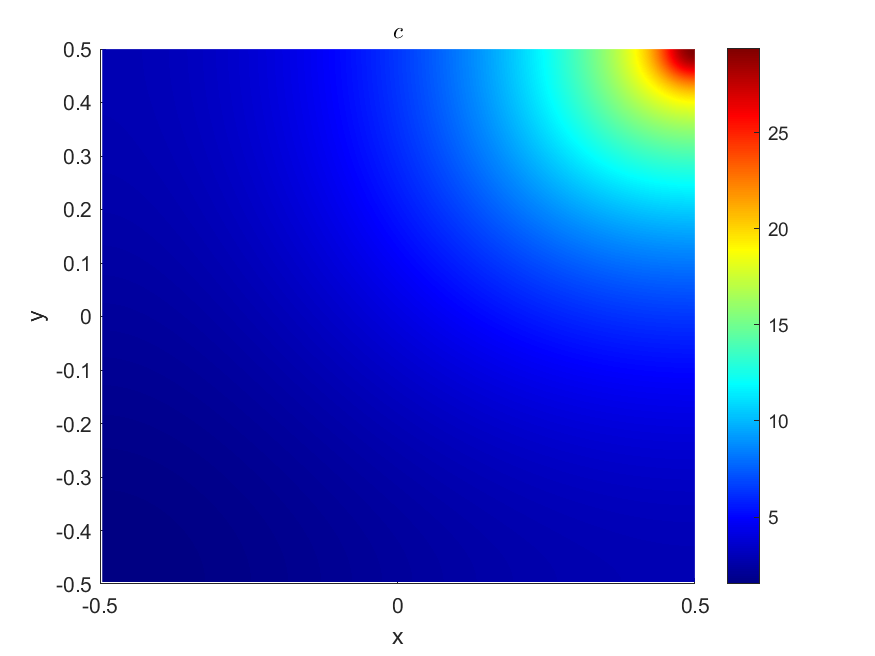}
	\caption{ Contour plots of $\rho$ (up) and $c$ (down) on non-uniform grids at time instants $t=0, 0.15, 0.163$ (from left to right) with $M = 200$ for Example \ref{exam:bu:corner}.}
	\label{fig:blowup—cor-evolution}
\end{figure}

Using a fixed time stepsize $\tau = 10^{-3}$ and $M=200$, Figure \ref{fig:blowup—cor-evolution} depicts the contours of $\rho$ and $c$ at different time instants $t=0$, 0.15, 0.163, respectively. It shows that both the maximum values of $\rho$ and $c$ gradually move toward the corner (0.5, 0.5) of the rectangular domain. As they are approaching the corner, the cell density undergoes a rapid and huge increase, and finally getting blow-up within a finite time. Meanwhile, the chemoattractant concentration $c$ also exhibits a pronounced peak at the corner (0.5, 0.5) of the domain. Our simulation results for this scenario indicate that the blow-up occurs around $t = 0.163$, which is closely aligned with the estimated blow-up time reported in \cite{EK'12}.

\begin{example}\label{exam:bu:3D} (Blow-up in 3D Keller--Segel model.) In the last example, we also consider the modeling of the three-dimensional (3D) Keller--Segel model \eqref{model:ks} in a cubic domain $\Omega=(0,1) \times(0,1)\times(0,1)$. The initial conditions are given as
		\begin{equation*}
			\begin{aligned}
				& \rho(x, y, z, 0)=1000 \exp\left( -100 ((x-0.5)^2+(y-0.5)^2+(z-0.5)^2 )\right) , \\
				& c(x, y, z, 0)=500 \exp\left( -50 ((x-0.5)^2+(y-0.5)^2+(z-0.5)^2 )\right)  .
			\end{aligned}
		\end{equation*}	
	\end{example}

\begin{figure}[!htbp]
	\centering
	\includegraphics[width=0.32\linewidth]{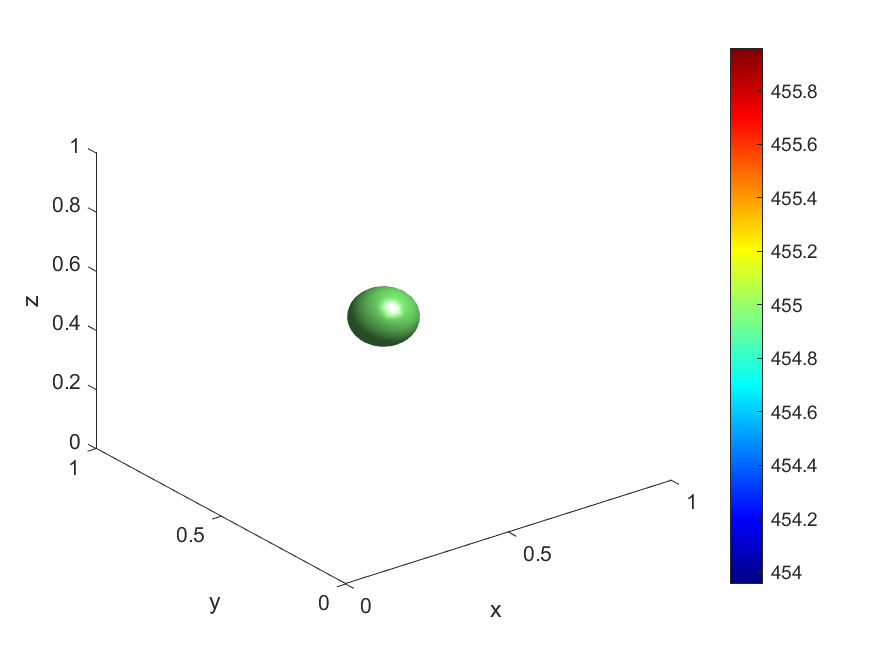}
	\includegraphics[width=0.32\linewidth]{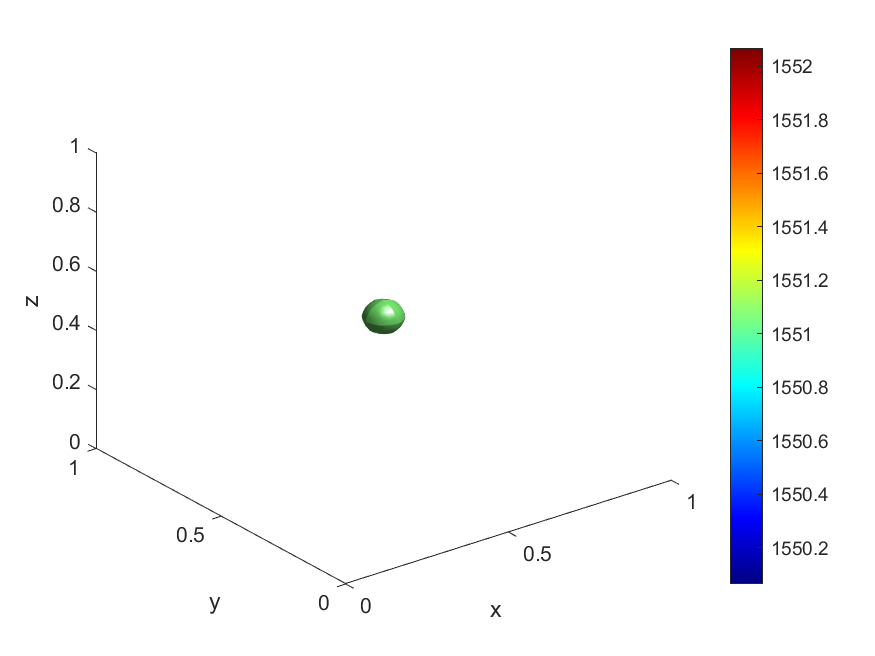}
	\includegraphics[width=0.32\linewidth]{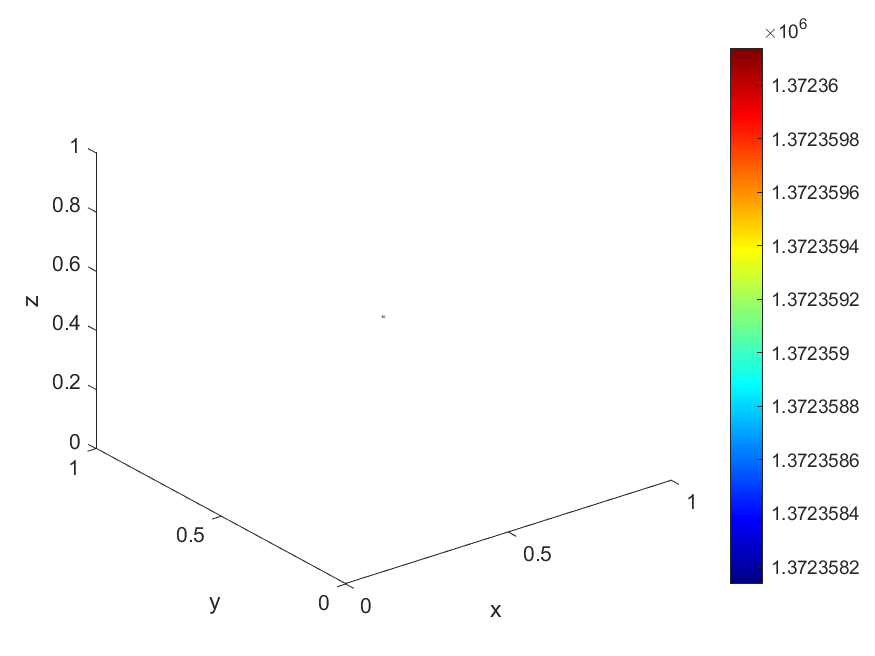}
	\includegraphics[width=0.32\linewidth]{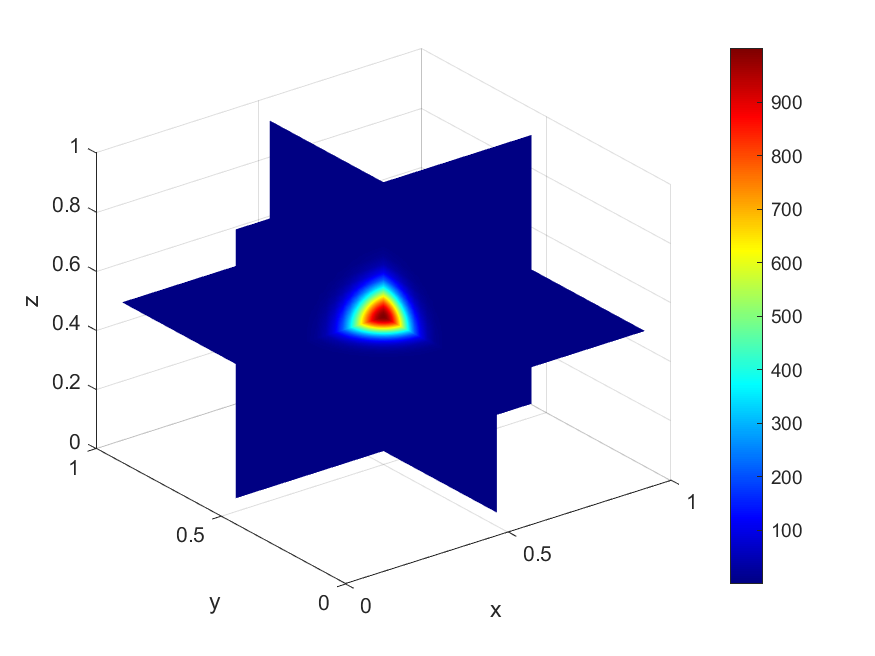}
	\includegraphics[width=0.32\linewidth]{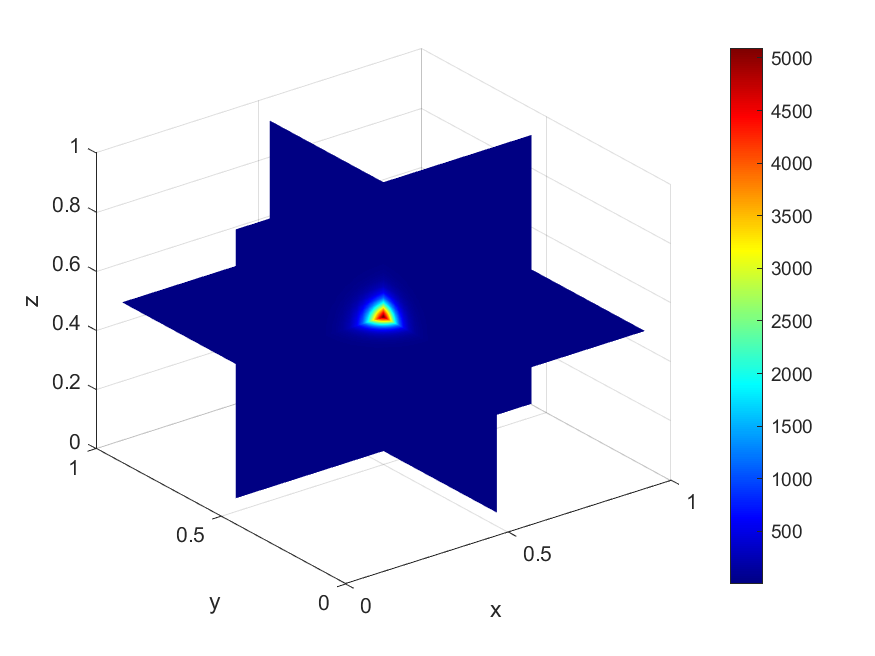}
	\includegraphics[width=0.32\linewidth]{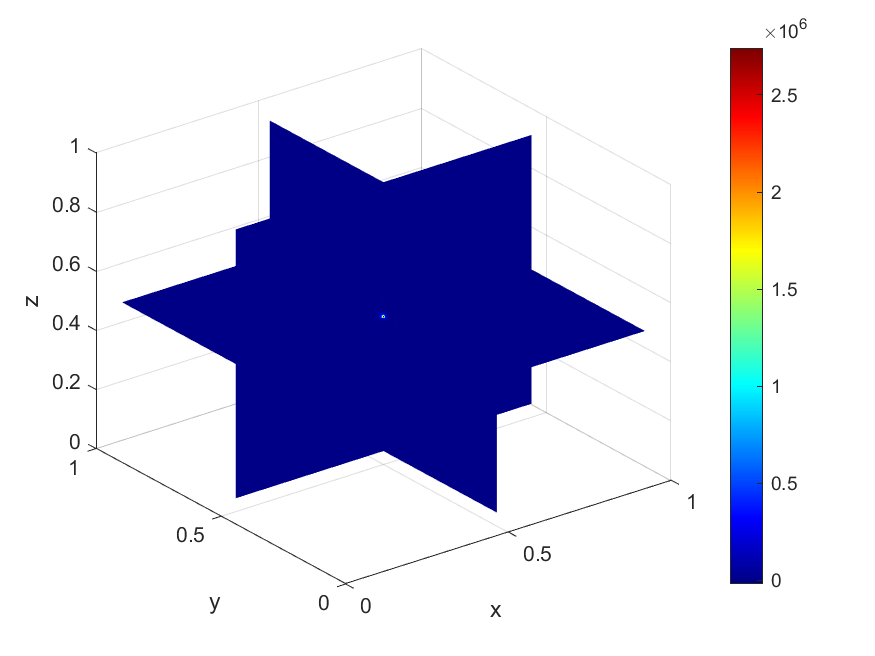}
    \caption{ Snapshots of $\rho$ (up) and their corresponding slices (down) at time instants $t=0, 1\times 10^{-5}, 5\times 10^{-5}$ (from left to right) with $M = 40$ for Example \ref{exam:bu:3D}. }
\label{fig:3D}
\end{figure}

In this example, similar middle refinement grids, as illustrated in \eqref{grid:mid} and Figure \ref{fig:gridmid}, are employed. We set $\tau=10^{-6}$. Snapshots of the density field $\rho$ and their corresponding slices at  time instants $t = 0,\, 1 \times 10^{-5},\, 5 \times 10^{-5}$ are presented in Figure \ref{fig:3D}, with a grid resolution of $M = 40$. 
The spherical isosurface, shown in the first row of Figure \ref{fig:3D}, reflects the symmetry of the numerical solution in each space dimension. The second row of Figure \ref{fig:3D} displays three extracted data slices from the planes $x = 0.5$, $y = 0.5$, and $z = 0.5$ in the three-dimensional space, respectively. 
Additionally, Figure \ref{fig:3D} reveals a blow-up phenomenon occurring in each space dimension of the 3D Keller--Segel system. The observed results are similar to those of 2D model. However, the computational complexity is significantly higher in the 3D scenario. Therefore, it is necessary to develop much more efficient numerical algorithms that avoid solving variable-coefficient algebraic system at each time step.
\section{Conclusion}\label{sec:conclusion}
This paper introduces a decoupled, linearized, and mass-conservative BCFD scheme for modeling of the classical Keller--Segel system, see  \eqref{2D:PC-BCFD}--\eqref{2D:CN-BCFD:IBc}. We demonstrate several key features of the developed scheme:
\begin{itemize}
    \item The scheme ensures mass conservation for the cell density at the discrete level, see Theorem \ref{thm:MassConserve} for details.
    
    \item Optimal second-order error estimates are rigorously established for the cell density $\rho$ in the discrete $L^{\infty}(L^2)$-norm and for the chemoattractant concentration $c$ in the $L^{\infty}(H^1)$-norm on non-uniform spatial grids, as proved in Theorem \ref{thm:coverg}, where the mathematical induction method and the discrete energy method are applied. Moreover, the existence and uniqueness of the solutions are also demonstrated in Theorem \ref{thm:existence}.
    
    \item By adopting mesh refinement strategies in areas with large gradient deformation (see \eqref{grid:mid}--\eqref{grid:cor}), the scheme can effectively simulate the blow-up phenomenon, leading to more accurate and insightful results.  
\end{itemize}

Numerical experiments underscore the robustness and efficiency of the DeC-MC-BCFD scheme in capturing the complex dynamics of the Keller--Segel system. This is particularly evident in scenarios involving center and corner blow-up phenomena. The Keller--Segel model also has some other intrinsic physical properties, such as positivity-preserving and energy dissipation law (see Refs. \cite{AGR'23, HS'21, HZ'23, LWZ'18, SX'20}).
In the future, our aim is to develop some finite difference schemes that can preserve these properties for the Keller--Segel chemotaxis system, even for the coupled chemotaxis–fluid models (see Refs. \cite{SN'19, HFXW'21, FHW'21}).

\section*{CRediT authorship contribution statement}
			\textbf{Jie Xu}: Methodology, Formal analysis, Software, Writing-Original draft.
			\textbf{Hongfei Fu}: Conceptualization, Supervision, Writing-Reviewing and Editing, Methodology, Funding acquisition.
		
\section*{Declaration of competing interest} The authors declare that they have no competing interests.
			
\section*{Data availability}   
A free github repository contains the Matlab code employed, which can be accessed through the following link: 
\href{https://github.com/HFu20/Keller-Segel/tree/Dec-MC-BCFD}{https://github.com/HFu20/Keller-Segel/tree/Dec-MC-BCFD}.

\section*{Acknowledgements}
This work was supported in part by the National Natural Science Foundation of China (No. 12131014), by the Natural Science Foundation of Shandong Province (No. ZR2024MA023), by the Fundamental Research Funds for the Central Universities (Nos. 202264006, 202261099) and by the OUC Scientific Research Program for Young Talented Professionals.

\bibliographystyle{elsarticle-num}
\bibliography{reference.bib}

\begin{thebibliography}{10}
\expandafter\ifx\csname url\endcsname\relax
  \def\url#1{\texttt{#1}}\fi
\expandafter\ifx\csname urlprefix\endcsname\relax\def\urlprefix{URL }\fi
\expandafter\ifx\csname href\endcsname\relax
  \def\href#1#2{#2} \def\path#1{#1}\fi

\bibitem{KS'70}
E.~Keller, L.~Segel, Initiation of slide mold aggregation viewed as an
  instability, J. Theor. Biol. 26 (1970) 399--415.

\bibitem{KS'71}
E.~Keller, L.~Segel, Model for chemotaxis, J. Theor. Biol. 30 (1971) 225--234.

\bibitem{CC'08}
V.~Calvez, L.~Corrias, {The parabolic-parabolic {K}eller--{S}egel model in
  $R^2$}, Commun. Math. Sci. 6 (2008) 417--447.

\bibitem{LWZ'18}
J.~Liu, L.~Wang, Z.~Zhou, Positivity-preserving and asymptotic preserving
  method for 2d {K}eller--{S}egel equations, Math. Comp. 87 (2018) 1165--1189.

\bibitem{LSY'17}
X.~Li, C.~Shu, Y.~Yang, Local discontinuous {G}alerkin method for the
  {K}eller--{S}egel chemotaxis model, J. Sci. Comput. 73 (2017) 943--967.

\bibitem{XFH'19}
X.~Xiao, X.~Feng, Y.~He, Numerical simulations for the chemotaxis models on
  surfaces via a novel characteristic finite element method, Comput. Math.
  Appl. 78 (2019) 20--34.

\bibitem{SN'19}
M.~Sulman, T.~Nguyen, A positivity preserving moving mesh finite element method
  for the {K}eller--{S}egel chemotaxis model, J. Sci. Comput. 80 (2019)
  649--666.

\bibitem{SX'20}
J.~Shen, J.~Xu, Unconditionally bound preserving and energy dissipative schemes
  for a class of {K}eller--{S}egel equations, SIAM J. Numer. Anal. 58 (2020)
  1674--1695.

\bibitem{HS'21}
F.~Huang, J.~Shen, Bound/positivity preserving and energy stable scalar
  auxiliary variable schemes for dissipative systems: applications to
  {K}eller--{S}egel and {P}oisson--{N}ernst--{P}lanck equations, SIAM J. Sci.
  Comput. 43 (2021) A1832--A1857.

\bibitem{TSL'19}
R.~Tyson, L.~Stern, R.~LeVeque, Fractional step methods applied to a chemotaxis
  model, J. Math. Biol. 41 (2000) 455--475.

\bibitem{M'03}
D.~Manoussaki, A mechanochemical model of angiogenesis and vasculogenesis,
  ESAIM: Math. Model. Numer. Anal. 37 (2003) 581--599.

\bibitem{CE'18}
A.~Chertock, Y.~Epshteyn, H.~Hu, A.~Kurganov, High-order positivity-preserving
  hybrid finite-volume-finite-difference methods for chemotaxis systems, Adv.
  Comput. Math. 44 (2018) 327350.

\bibitem{S'09}
N.~Saito, Conservative numerical schemes for the {K}eller--{S}egel system and
  numerical results, RIMS Kôkyûroku Bessatsu 15 (2009) 125--146.

\bibitem{EK'12}
Y.~Epshteyn, Upwind-difffference potentials method for
  {P}atlak--{K}eller--{S}egel chemotaxis model, J. Sci. Comput. 53 (2012)
  689--713.

\bibitem{HZ'23}
J.~Huang, X.~Zhang, Positivity-preserving and energy-dissipative finite
  difference schemes for the {F}okker--{P}lanck and {K}eller--{S}egel
  equations, IMA J. Numer. Anal 43 (2023) 1450--1484.

\bibitem{BGG'20}
J.~Benito, A.~García, L.~Gavete, M.~Negreanu, F.~Ureña, A.~Vargas, Solving a
  fully parabolic chemotaxis system with periodic asymptotic behavior using
  generalized finite difference method, Appl. Numer. Math. 157 (2020) 356--371.

\bibitem{DA'19}
M.~Dehghan, M.~Abbaszadeh, The simulation of some chemotactic bacteria patterns
  in liquid medium which arises in tumor growth with blow-up phenomena via a
  generalized smoothed particle hydrodynamics ({GSPH}) method, Eng. Comput. 35
  (2019) 875--892.

\bibitem{AWY'97}
T.~Arbogast, M.~Wheeler, I.~Yotov, Mixed finite elements for elliptic problems
  with tensor coefficients as cell-centered finite differences, SIAM J. Numer.
  Anal. 34 (1997) 828--852.

\bibitem{RP'12}
H.~Rui, H.~Pan, A block-centered finite difference method for the
  {D}arcy--{F}orchheimer model, SIAM J. Numer. Anal. 50 (2012) 2612--2631.

\bibitem{XXF'22}
J.~Xu, S.~Xie, H.~Fu, A two-grid block-centered finite difference method for
  the nonlinear regularized long wave equation, Appl. Numer. Math. 171 (2022)
  128--148.

\bibitem{WXF'24}
X.~Wang, J.~Xu, H.~Fu, A linearlized mass-conservative fourth-order
  block-centered finite difference method for the semilinear sobolev equation
  with variable coefficients, Commun. Nonlinear Sci. Numer. Simul. 130 (2024)
  107778.

\bibitem{SXLF'21}
Y.~Shi, S.~Xie, D.~Liang, K.~Fu, High order compact block-centered finite
  difference schemes for elliptic and parabolic problems, J. Sci. Comput. 87
  (2021) 86.

\bibitem{LSR'19}
X.~Li, J.~Shen, H.~Rui, Energy stability and convergence of {SAV}
  block-centered finite difference method for gradient flows, Math. Comp. 88
  (2019) 2047--2068.

\bibitem{RL'15}
H.~Rui, W.~Liu, A two-grid block-centered finite difference method for
  {D}arcy--{F}orchheimer flow in porous media, SIAM J. Numeri. Anal. 53 (2015)
  1941--1962.

\bibitem{WW'88}
A.~Weiser, M.~Wheeler, On convergence of block-centered finite differences for
  elliptic problems, SIAM J. Numer. Anal. 25 (1988) 351--375.

\bibitem{RP'13}
H.~Rui, H.~Pan, Block-centered finite difference methods for parabolic equation
  with time-dependent coefficient, Jpn. J. Ind. Appl. Math. 30 (2013) 681--699.

\bibitem{SZG'23}
Z.~Sun, Q.~Zhang, G.~Gao, Finite difference methods for nonlinear evolution
  equations, De Gruyter, Berlin/Boston, 2023.

\bibitem{BH'71}
J.~Bramble, S.~Hilbert, Bounds for a class of linear functionals with
  application to hermite interpolation, Numer. Math. 16 (1971) 362369.

\bibitem{BGM'07}
G.~Berikelashvili, M.~Gupta, M.~Mirianashvili, Convergence of fourth order
  compact difference schemes for three-dimensional convection-diffusion
  equations, SIAM J. Numer. Anal. 45 (2007) 443--455.

\bibitem{DWW'98}
C.~Dawson, M.~Wheeler, C.~Woodward, A two-grid finite difference scheme for
  nonlinear parabolic equations, SIAM J. Numer. Anal. 35 (1998) 435--452.

\bibitem{AGR'23}
D.~Acosta-Soba, F.~Guillén-González, J.~Rodríguez-Galván, An
  unconditionally energy stable and positive upwind {DG} scheme for the
  {K}eller--{S}egel model, J. Sci. Comput. 97 (2023) 18.

\bibitem{HFXW'21}
X.~Huang, X.~Feng, X.~Xiao, K.~Wang, Fully decoupled, linear and
  positivity-preserving scheme for the chemotaxis--{S}tokes equations, Comput.
  Methods Appl. Mech. Engrg. 383 (2021) 113909.

\bibitem{FHW'21}
X.~Feng, X.~Huang, K.~Wang, Error estimate of unconditionally stable and
  decoupled linear positivity-preserving {FEM} for the chemotaxis--{S}tokes
  equations, SIAM J. Numer. Anal. 59 (2021) 3052--3076.

\end{thebibliography}

\end{document}